# On the numerical solution of the Van der Pol equation using collocation

Authored by J. H. M. Darbyshire,

acknowledging D. Roscoe for his suggestions and review.

Submitted, April 2015, to

*The Open University*

in partial fulfilment for the award of the degree of

MASTER OF MATHEMATICS.

# Abstract

This study introduces the reader to the theory of approximating the solution(s) of a non-linear, second order, ordinary differential equation (ODE) with piecewise polynomial functions by using the collocation method[3]. It then focuses on the application of the method to generate collocation approximations for the Van der Pol equation with initial conditions and with varying equation parameter, $\mu$. It is shown that the collocation method in its original form is impractical for generating these approximations due to the numerical costs of producing such an approximation, particularly as the stiffness of the equation increases with parameter, $\mu$. An adaptation to the method, termed segmented collocation, is proposed, reliant upon on the specific structure and convergence for initial value problems, and shown theoretically and numerically to be capable of generating equivalent approximations for much superior costs. Further study is proposed; to compare the segmented method to other numerical methods traditionally employed with equations of this type such as the Runge-Kutta family; and to develop and test the method in order to ascertain its complete usefulness and competitive standing in the context of approximating second order, ODEs with initial conditions.

# Contents

# Introduction

## 1.1 Introduction

In this chapter we present the necessary mathematics for the complete understanding of the methods used to generate piecewise polynomial collocation approximations to the solution of second order, non-linear, ordinary differential equations (ODEs). This theory is presented in a practical sense and structured with successive sections typically building on previous ones, as might be expected.

We also introduce the reader to Balthazar Van der Pol and the equation that bears his name, which forms the basis of this study, that being the application of the aforementioned collocation method to generate approximations to second order, ODEs. The Van der Pol is an excellent example of a typical, and well studied, equation of its type as a non-linear oscillatory equation with a stable limit-cycle, making it highly suitable for purpose here.

## 1.2 Introduction to theory

### 1.2.1 The space of polynomial functions

We state without proof that the set of all polynomial functions of order $k$ with real coefficients forms a vector space. Indeed the monomial basis;

$$\{1, x, x^2, \dots, x^{k-1}\} \ ,$$

spans the space, and the dimension of the space is observed to be $k$. We denote this space by, $\mathbf{\Pi}_{<k}$.

### 1.2.2 The space of piecewise polynomial (pp) functions

In a domain, $[a..b]$, let $g(x)$ be an arbitrary continuous function and consider a partition;

$$\{\xi_1, \xi_2, \dots, \xi_{l+1}\} \quad \text{with,} \quad a = \xi_1 < \xi_2 < \dots < \xi_l < \xi_{l+1} = b,$$

such that,

$$[a..b] = \bigcup_{i=1}^{l} [\xi_i..\xi_{i+1}] \, .$$

For each subdomain $[\xi_i..\xi_{i+1}]$ we denote by $f_i(x)$, any polynomial function of order $k$, and highlight:

$$f_i(x) = \begin{cases} f_i(x) : x \in [\xi_i..\xi_{i+1}] \\ 0 \qquad : x \notin [\xi_i..\xi_{i+1}] \end{cases} .$$

Over the domain, $[a..b]$, the function $f(x) = \{f_i(x) : x \in [\xi_i..\xi_{i+1}]\}$ is termed the pp function of order $k$. It is observed that under this definition for, $i = 2, \dots, l$, that $f(\xi_i)$ has dual values, where it is continuous from the left and from the right. It is clear that two pp functions are equivalent if and only if their piecewise functions, $f_i$, agree for all $x$.

Restricting the pp to be, at the interior breakpoints, either left continuous or right continuous then it is simpler to state without proof that the set of all pp functions with break sequence, $\boldsymbol{\xi} = (\xi_i)_1^{l+1}$ forms a vector space. Indeed for a right continuous pp function the piecewise monomial basis:

$$\Big\{\{M_i(1), M_i(x-\xi_i), \dots, M_i(x-\xi_i)^{k-1}\} : i = 1, \dots, l\Big\}$$

where the operator, $M_i f(x) := \begin{cases} f(x) & \text{for } x \in [\xi_i..\xi_{i+1}) \\ 0 & \text{otherwise} \end{cases}$, is observed to span the space, and the dimension is observed to be $kl$ (note that the definition of the operator is extended to include the rightmost point $x = b$). We denote the space of all pp functions of order $k$ with break sequence $\boldsymbol{\xi}$ by, $\boldsymbol{\Pi}_{<k,\boldsymbol{\xi}}$ .

In the context of this dissertation, and as is typical, we are interested in finding the pp approximation, $f$, to $g$ which attempts to minimise,

$$\|f - g\|_\infty := \max_{a \le x \le b} |f(x) - g(x)|.$$

### 1.2.3 The space of pp functions with imposed continuity conditions

In an attempt to find $f$, as a pp approximation to $g$, we seek to impose the conditions that $f$, and potentially its derivatives $D^m f$, should be continuous at the interior breakpoints of $\boldsymbol{\xi}$. It is clear that such a space of pps is contained within $\boldsymbol{\Pi}_{<k,\boldsymbol{\xi}}$ and we state without proof that this is indeed a subspace. Each continuity condition imposed reduces the dimension of the subspace from $kl$ by one, and it is customary to denote, by $\boldsymbol{\nu}$, the integer-vector of length $(l-1)$, which indicates the number of continuity conditions imposed on $f$ at each interior breakpoint. Note that it is traditional to assume that any derivative, $D^m f$, should not be considered continuous at $\xi_i$ unless the function value, $f(\xi_i)$, and any lower derivatives, $D^{<m} f(\xi_i)$ are imposed to be continuous in addition. For example if $\nu_{i-1} = m$ one would assume $f(\xi_i), Df(\xi_i), \dots, D^{m-1} f(\xi_i)$ harbour the imposed continuity conditions at $\xi_i$. Note, additionally, that the number of continuity conditions imposed at each interior breakpoint must be less than $k$, otherwise if it equals $k$, whilst permissible, for any practical consideration that breakpoint might simply not exist. The resultant subspace is denoted by $\boldsymbol{\Pi}_{<k,\boldsymbol{\xi},\boldsymbol{\nu}}$ and has dimension, $n = kl - \sum_{i=1}^{l-1} \nu_i$.

As an example suppose $g = e^x$ and we seek the pp approximation, $f$, of order four, with breakpoint sequence, $\boldsymbol{\xi} = \{1, 2, 3, 4\}$, and with $\boldsymbol{\nu} = \{1, 3\}$. The dimension of the space, $n = 8$, which indicates that $f$ can only be determined with knowledge of eight distinct values of $g$ or its derivatives, $D^m g$. More than that, an appropriate distribution of this information in consideration with $\boldsymbol{\xi}$ is also required according to the Schoenberg-Whitney theorem (see section 1.2.8), otherwise, some piecewise polynomials might be over-specified while others might be under-specifed.

### 1.2.4 Construction of a knot sequence, t

The combination of information provided by the order, $k$, of the pp, $\boldsymbol{\xi}$ and $\boldsymbol{\nu}$ to define the space of pps, $\boldsymbol{\Pi}_{<k,\boldsymbol{\xi},\boldsymbol{\nu}}$, can be attributed to a single sequence, typically termed knot sequence, and is denoted by, $\mathbf{t}$.

The construction of the knot sequence here is that identical to the description which forms part of the Curry and Schoenberg theorem, and the steps are as follows;

(i) let the first $k$ knots equal $\xi_1$: $\{t_j : j = 1, \dots, k\} := \xi_1$,

(ii) each interior breakpoint, $\xi_i$, has a number of knots equal to $k - \nu_{i-1}$:

$$\{t_{k(i-1)+1-\sum_{j=1}^{i-2}\nu_j}, \dots, t_{ki-\sum_{j=1}^{i-1}\nu_j}\} := \xi_i \quad \forall\ i = 2, \dots, l\ ,$$

(iii) the final $k$ knots are equal to $\xi_{l+1}$: $\{t_{n+j} : j = 1, \dots, k\} := \xi_{l+1}$.

Continuing the previous example suppose again that, $k = 4$, $\boldsymbol{\xi} = \{1, 2, 3, 4\}$, $\boldsymbol{\nu} = \{1, 3\}$, then the resultant knot sequence is:

$$\mathbf{t} = \{1, 1, 1, 1, 2, 2, 2, 3, 4, 4, 4, 4\}.$$

### 1.2.5 B-splines

Having established a procedure for constructing a general knot sequence, $\mathbf{t}$, we are in a position to begin the discussion about the general construction of B-splines. B-splines are themselves carefully constructed pp functions of specific order. In this sense, carefully constructed means that their design permits certain properties that are favourable for the numerical determinations of pp approximations by reducing round off errors and other phenomena that other structured polynomials can suffer, for example the truncated power basis[1, p84-85].

Firstly we define the basic B-splines, which are of order one:

$$B_{i,1}(x) = \begin{cases} 1, & \text{if } t_i \le x < t_{i+1} \\ 0, & \text{otherwise} \end{cases} .$$

For each knot sequence, $\mathbf{t}$, containing $n + k$ knots there are associated with it $n + k - 1$ basic B-splines. In regard to repeated knots the definition is such that the basic B-spline, $B_{i,1}$, takes the value zero when $t_i = x = t_{i+1}$. Additionally all of the basic B-splines defined above are observed to be right-continuous, and in order to be a well defined set on the basic interval, $I_{k,\mathbf{t}} = [t_k..t_{n+1}] = [a..b]$, then the rightmost basic B-spline must also be made to be left-continuous, so that we have the amended definition:

$$B_{n+k-1,1}(x) = \begin{cases} 1, & \text{if } t_{n+k-1} \le x \le t_{n+k} \\ 0, & \text{otherwise} \end{cases} .$$

Turning now to B-splines of order two, these are defined through the recurrence relation:

$$B_{i,2}(x) = \frac{x - t_i}{t_{i+1} - t_i} B_{i,1}(x) + \frac{t_{i+2} - x}{t_{i+2} - t_{i+1}} B_{i+1,1}(x)\ .$$

Again we run into some trouble where the denominator is zero for repeated knots. In this instance the corresponding B-spline of that coefficient will be zero for all $x$ and the over arching convention we adopt, and that is typical to adopt, is that "anything multiplied by zero is zero".

As an example for the knot sequence, as above, $\mathbf{t} = \{1, 1, 1, 1, 2, 2, 2, 3, 4, 4, 4, 4\}$ there are only three non-zero basic B-splines, $B_{4,1}$, $B_{7,1}$, and $B_{8,1}$. Using the above relation we determine the non-zero B-splines of order two, and illustrate them in figure 1.1:

$$B_{3,2}(x) = \frac{x-1}{1-1} B_{3,1}(x) + \frac{2-x}{2-1} B_{4,1}(x) = \begin{cases} 2-x, & \text{if } 1 \le x < 2 \\ 0, & \text{otherwise} \end{cases} ,$$

$$B_{4,2}(x) = \frac{x-1}{2-1} B_{4,1}(x) + \frac{2-x}{2-2} B_{5,1}(x) = \begin{cases} x-1, & \text{if } 1 \le x < 2 \\ 0, & \text{otherwise} \end{cases} ,$$

$$B_{6,2}(x) = \frac{x-2}{2-2} B_{6,1}(x) + \frac{3-x}{3-2} B_{7,1}(x) = \begin{cases} 3-x, & \text{if } 2 \le x < 3 \\ 0, & \text{otherwise} \end{cases} ,$$

$$B_{7,2}(x) = \frac{x-2}{3-2} B_{7,1}(x) + \frac{4-x}{4-3} B_{8,1}(x) = \begin{cases} x-2, & \text{if } 2 \le x < 3 \\ 4-x, & \text{if } 3 \le x < 4 \\ 0, & \text{otherwise} \end{cases} ,$$

$$B_{8,2}(x) = \frac{x-3}{4-3} B_{8,1}(x) + \frac{4-x}{4-4} B_{9,1}(x) = \begin{cases} x-3, & \text{if } 3 \le x < 4 \\ 0, & \text{otherwise} \end{cases} .$$

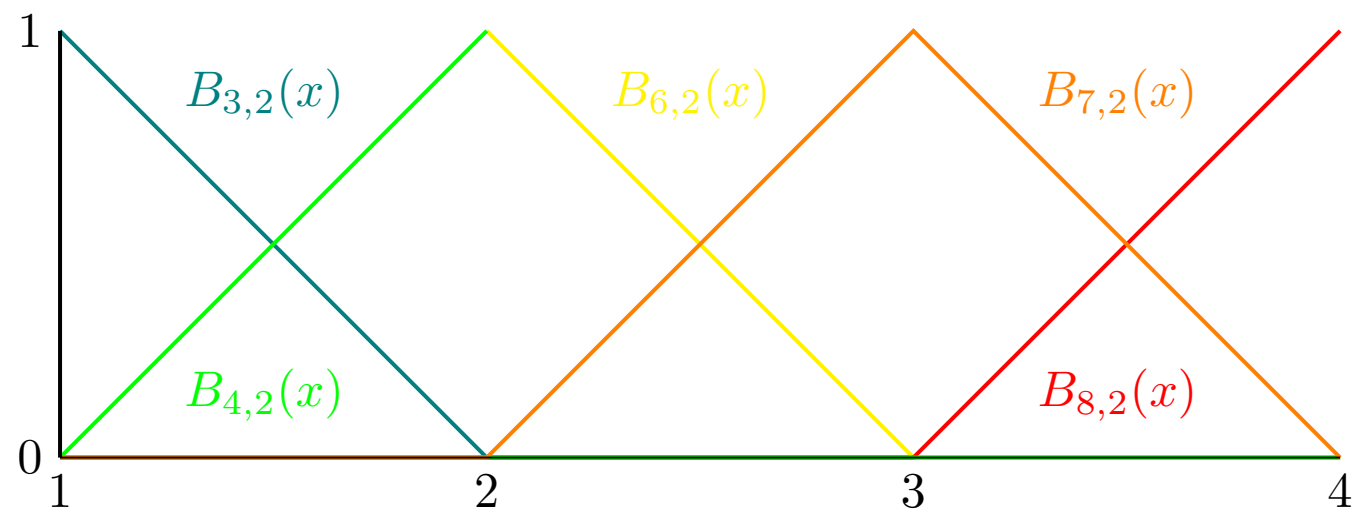

Figure 1.1: Graphical illustration of the five B-splines of order two for the example knot sequence, $\mathbf{t}$.

Generalising the recurrence relation, for any $\mathbf{t}$ with $n+k$ knots we can obtain upto $n$ B-splines of order $k$:

$$B_{i,k}(x) = \frac{x - t_i}{t_{i+k-1} - t_i} B_{i,k-1}(x) + \frac{t_{i+k} - x}{t_{i+k} - t_{i+1}} B_{i+1,k-1}(x) \ . \tag{1.1}$$

Continuing the above example, figure 1.2 depicts the eight, order four B-splines for the given knot sequence, $\mathbf{t}$, generated using the above recurrence relation. We highlight the distinction in the figure between B-splines about $x = 2$, and $x = 3$, which show a linear combination of these B-splines will always produce continuity in $f$ at $x = 2$ and $x = 3$, but will likely produce a discontinuity in first and second derivative of $f$ at $x = 2$ but not at $x = 3$, where it will be continuous. Of course this was expected from the configuration of $\boldsymbol{\nu} = \{1, 3\}$.

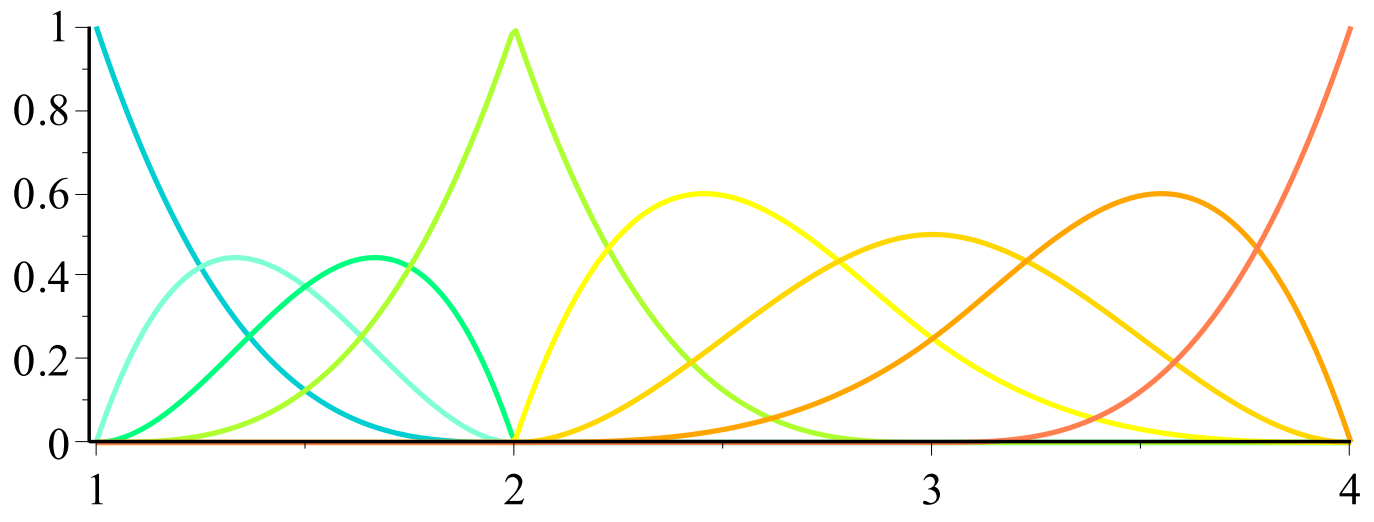

Figure 1.2: Graphical illustration of the eight B-splines of order four for the example knot sequence, $\mathbf{t}$.

Where it is perhaps not clear from context, $B_{i,k,\mathbf{t}}(x)$ denotes the $i$-th B-spline of order $k$ that has been generated from the knot sequence $\mathbf{t}$, evaluated at point $x$.

### 1.2.6 Properties of B-splines

The theory of solving differential equations using the collocation method relies on the linear combination of B-splines. Below we state here some of the important properties of B-splines that will be referred to in upcoming sections:

**Support and positivity[1, p91]**

$$B_{i,k,\mathbf{t}}(x) \begin{cases} > 0, & \text{for} \quad t_i < x < t_{i+k} \\ = 0, & \text{for} \quad x < t_i \ \text{ and } \ x > t_{i+k} \end{cases} . \tag{1.2}$$

**Positive and local, partition of unity[1, p96]**

$$\sum_i B_{i,k}(x) = 1, \quad x \in [a..b] \ . \tag{1.3}$$

**Derivatives**

From the definition it is clear that the derivative of basic B-splines is zero, and at the jump discontinuities we take the derivative from the appropriate left or right side so that:

$$DB_{i,1}(x) = 0 \ .$$

We assert that the generalised derivative of a B-spline is:

$$DB_{i,k}(x) = (k-1)\left( \frac{B_{i,k-1}(x)}{t_{i+k-1} - t_i} - \frac{B_{i+1,k-1}(x)}{t_{i+k} - t_{i+1}} \right) \ . \tag{1.4}$$

The proof of this is widely publicised, with the original[2] based upon the divided difference definition of B-splines (which is not covered in this study), another version[1] built upon the recurrence relation and dual functionals and a third[4] using an inductive argument based solely on the recurrence relation in equation 1.1.

**Endpoint support**

We also state, which is inferred from equations 1.1 and 1.4, the following:

$$B_{i,k}(a) = \begin{cases} 1, & \text{if } \ i = 1 \\ 0, & \text{otherwise} \end{cases} \ , \quad B_{i,k}(b) = \begin{cases} 0, & \text{otherwise} \\ 1, & \text{if } \ i = n \end{cases} \ , \tag{1.5}$$

$$\begin{aligned} D^m B_{i,k}(a) &\begin{cases} \neq 0, & \text{for } i = 1, \dots, m+1 \\ = 0, & \text{otherwise} \end{cases} \ , \\ D^m B_{i,k}(b) &\begin{cases} = 0, & \text{otherwise} \\ \neq 0, & \text{for } i = n-m, \dots, n \end{cases} \ . \end{aligned} \tag{1.6}$$

### 1.2.7 Spline space, $\$_{k,\mathbf{t}}$

We define the spline space of order $k$ with knot sequence $\mathbf{t}$, as described above, as any linear combination of B-splines of order $k$ generated from $\mathbf{t}$:

$$\$_{k,\mathbf{t}} := \left\{ \sum_{i=1}^{n} \alpha_i B_{i,k,\mathbf{t}}(x) \ : \text{for real } \ \alpha_i \right\} \ .$$

**Curry and Schoenberg theorem[6]**

This theorem states that for, $\boldsymbol{\xi}$ as defined above, $\boldsymbol{\nu}$ as defined above, and for dimension $n = \dim \mathbf{\Pi}_{<k,\boldsymbol{\xi},\boldsymbol{\nu}}$, as defined above, $\mathbf{t}$ as defined above, then the sequence of B-splines, $B_{1,k}, \dots, B_{n,k}$, each of order $k$ and generated from knot sequence $\mathbf{t}$ is a basis for $\mathbf{\Pi}_{<k,\boldsymbol{\xi},\boldsymbol{\nu}}$ considered as functions on $I_{k,\mathbf{t}} = [t_k..t_{n+1}] = [a..b]$. Specifically then a linear combination of B-splines span the space, i.e.:

$$\$_{k,\mathbf{t}} = \mathbf{\Pi}_{<k,\boldsymbol{\xi},\boldsymbol{\nu}} \quad \text{on} \quad I_{k,\mathbf{t}} \ .$$

### 1.2.8 Spline interpolation

As one might imagine the introduction previously of the spline space permits us now to discuss spline interpolation, which has been previously alluded to.

We seek to find our pp approximation $f := \sum_{i=1}^{n} \alpha_i B_{i,k,\mathbf{t}} \in \$_{k,\mathbf{t}}$, of $g$, such that $f$ agrees with $g$ at given datasites, $\boldsymbol{\tau} = \{\tau_1, \tau_2, \dots, \tau_m\}$, i.e.:

$$f(\tau_j) = g(\tau_j) \quad \forall \ j \ .$$

Given the dimension, $n$, of the basis of the space, it is natural to expect to accommodate $n$ interpolations, i.e $m = n$. Then,

$$\sum_{i=1}^{n} \alpha_i B_{i,k,\mathbf{t}}(\tau_j) = g(\tau_j) \quad \forall j = 1, \dots, n \ .$$

The matrix $\mathbf{B}(\boldsymbol{\tau}) = B_{i,k}(\tau_j)$ is termed the spline collocation matrix. Any deviation in the number of datasites results in an over-specified or under-specified system which is potentially unsolvable, unreliable or unstable. Additionally, one might expect some necessary relationship between the set of knots, $\mathbf{t}$, that defines the basis, and the set of datasites, $\boldsymbol{\tau}$, enforcing interpolation.

Indeed the following theorem clarifies this assertion:

**Schoenberg-Whitney theorem[5]**

Let $\boldsymbol{\tau}$ be strictly increasing and such that $a < t_i = \ldots = t_i + r = \tau_j < b$ implies $r < k-1$. Then the spline collocation matrix, $\mathbf{B}(\boldsymbol{\tau})$ of the above system is invertible, and the system uniquely solvable, if and only if,

$$B_{i,k}(\tau_i) \neq 0 \ ,$$

that is, if and only if, $t_i < \tau_i < t_{i+k} \ \forall \ i$ (except that $\tau_1 = t_1$ and $\tau_n = t_{n+k}$ are also permitted).

### 1.2.9 Approximation of linear, second order, ODE

Suppose $g(x)$ on $[a..b]$ is the solution to the differential equation:

$$D^2 g + p(x)Dg + q(x)g = s(x), \quad \beta_i g(a) = c_i, \ \ i = 1, 2 \ ,$$

where $\beta_i g = \beta_{i,1} g + \beta_{i,2} Dg$, for constants, $\beta_{i,1}, \beta_{i,2}$, and we assume further that these conditions are not linearly dependent. Also, for the avoidance of doubt, we highlight that these are initial conditions and that in general our discussion is specific to initial value problems except in the case of reference to other literature where the citation's mathematics are applicable to problems with boundary value conditions.

An existence and uniqueness theorem exists that states under the conditions of continuity for $p, q$, and $s$ a solution for $g$ exists and that that solution is unique.

Now suppose that we wish to approximate the solution $g$ by a pp $f \in \mathbf{\Pi}_{k,\boldsymbol{\xi},\boldsymbol{\nu}} = \$_{k,\mathbf{t}}$. The breakpoints, $\xi_i$ can be arbitrarily chosen but we impose two continuity conditions on the interior breakpoints (because the differential equation is of second order), so that $\boldsymbol{\nu} = \{2, \ldots, 2\}$, and then it follows that the order of the pp, $k \geq 3$. This results in dimension, $n = (k-2)l + 2$. Intuitively, this represents $(k-2)l$ interpolation site, and two initial, conditions that must be accommodated by our spline. The linear system then becomes:

$$\begin{aligned} \beta_j f(a) &= \sum_{i=1}^{n} \alpha_i \beta_j B_{i,k}(a) = c_j, \qquad j = 1, 2, \\ Lf(\tau_j) &= \sum_{i=1}^{n} \alpha_i L B_{i,k}(\tau_j) = s(\tau_j), \ \ j = 1, \ldots, (k-2)l, \end{aligned} \tag{1.7}$$

where $Lf(\tau_j) = D^2 f(\tau_j) + p(\tau_j)Df(\tau_j) + q(\tau_j)f(\tau_j)$.

We note firstly, using equations 1.5 and 1.6, that:

$$B_{i,k}(a) = \begin{cases} 1, & \text{if } \ i = 1 \\ 0, & \text{otherwise} \end{cases} ,$$

$$DB_{i,k}(a) = \begin{cases} \frac{1-k}{t_{k+1}-t_2}, & \text{if } \ i = 1 \\ \frac{k-1}{t_{k+1}-t_2}, & \text{if } \ i = 2 \\ 0, & \text{otherwise} \end{cases} .$$

The first two equations of the linear system representing the initial conditions reduce to the following in matrix form:

$$\begin{bmatrix} \beta_{1,1} - \beta_{1,2}\frac{k-1}{t_{k+1}-t_2}, & \beta_{1,2}\frac{k-1}{t_{k+1}-t_2}, & 0, & \ldots & 0 \\ \beta_{2,1} - \beta_{2,2}\frac{k-1}{t_{k+1}-t_2}, & \beta_{2,2}\frac{k-1}{t_{k+1}-t_2}, & 0, & \ldots & 0 \end{bmatrix} \boldsymbol{\alpha} = \begin{bmatrix} c_1 \\ c_2 \end{bmatrix} .$$

For the remaining $(k-2)l$ equations of the system we structure the datasites so that within each interval, $i \in \{1, \ldots, l\}$, we assign there to be specifically $(k-2)$ datasites, so that the following holds (and note this definition prohibits the coalescing of any such datasites);

$$\xi_i < \tau_{1+\gamma} < \ldots < \tau_{k-2+\gamma} < \xi_{i+1}, \quad \text{for} \quad \gamma = (i-1)(k-2) \ .$$

Within each interval the only B-splines that offer support are from the set $\{B_{1+\gamma,k}, \ldots, B_{k+\gamma,k}\}$, and hence in matrix form we can construct $(k-2)$ equations of the linear system for each interval, $i$:

$$\begin{bmatrix} \ldots & 0, & LB_{1+\gamma,k}(\tau_{1+\gamma}) & \ldots & LB_{k+\gamma,k}(\tau_{1+\gamma}) & 0, & \ldots \\ \vdots & & & \vdots & & & \vdots \\ \ldots & 0, & LB_{1+\gamma,k}(\tau_{k-2+\gamma}) & \ldots & LB_{k+\gamma,k}(\tau_{k-2+\gamma}) & 0, & \ldots \end{bmatrix} \boldsymbol{\alpha} = \begin{bmatrix} s(\tau_{1+\gamma}) \\ \vdots \\ s(\tau_{k-2+\gamma}) \end{bmatrix}$$

We refer to the complete, general linear system in matrix form as;

$$\mathbf{A}\boldsymbol{\alpha} = \mathbf{s}\ . \tag{1.8}$$

Provided that the datasites are independent of one another and don't coalesce, as per the above proposed structure, then the rows of the matrix, $\mathbf{A}$, are linearly independent and the it is invertible and the system uniquely solvable for $\boldsymbol{\alpha}$. Indeed by the choice of the positioning of the datasites within intervals and the nature of support of the B-splines and their derivatives the conditions of the Schoenberg-Whitney theorem can be observed to be satisfied in most cases. We make no reference here of the existence, or not, of a particular operator, $L$, which may be capable of producing a zero diagonal element for a particular datasite, $\tau_j$, say.

As an example let the number of intervals, $l = 3$, $k = 4$, then $n = 8$ and the linear system, where $X$ denotes a non-zero element detailed above, is:

$$\begin{bmatrix} X & X & 0 & 0 & 0 & 0 & 0 & 0 \\ X & X & 0 & 0 & 0 & 0 & 0 & 0 \\ \hline X & X & X & X & 0 & 0 & 0 & 0 \\ X & X & X & X & 0 & 0 & 0 & 0 \\ 0 & 0 & X & X & X & X & 0 & 0 \\ 0 & 0 & X & X & X & X & 0 & 0 \\ 0 & 0 & 0 & 0 & X & X & X & X \\ 0 & 0 & 0 & 0 & X & X & X & X \end{bmatrix} \begin{bmatrix} \alpha_1 \\ \alpha_2 \\ \alpha_3 \\ \alpha_4 \\ \alpha_5 \\ \alpha_6 \\ \alpha_7 \\ \alpha_8 \end{bmatrix} = \begin{bmatrix} c_1 \\ c_2 \\ \hline s(\tau_1) \\ s(\tau_2) \\ s(\tau_3) \\ s(\tau_4) \\ s(\tau_5) \\ s(\tau_6) \end{bmatrix} . \tag{1.9}$$

For this particular example, letting $\boldsymbol{\xi} = \{0, 1, 2, 3\}$ then the appropriate B-splines are graphically illustrated in figure 1.3, highlighting potential datasite positioning but ensuring the above structure is adhered to. It perhaps becomes clearer that, for this example, within each interval there are only ever four supporting B-splines, which factor into the rows of the matrix.

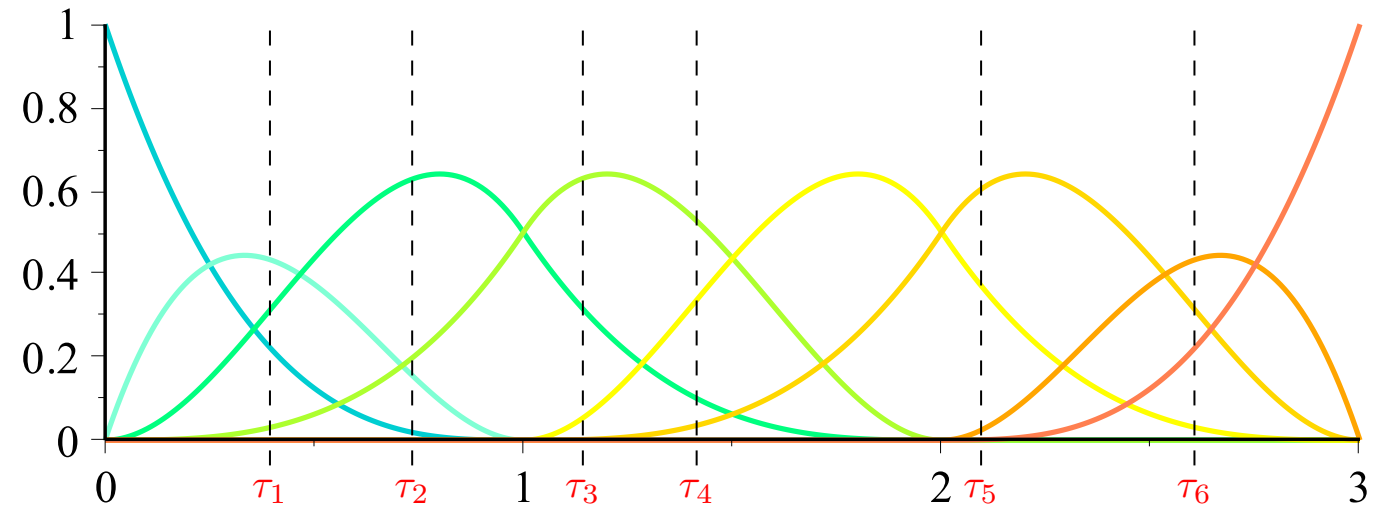

Figure 1.3: Graphical illustration of the eight B-splines lending support to each datasite, $\tau_j$, where here $\mathbf{t} = \{0, 0, 0, 0, 1, 1, 2, 2, 3, 3, 3, 3\}$.

### 1.2.10 Approximation of non-linear, second order, ODE

For $g(x)$ the solution, not necessarily unique, for the general second order differential equation with linear initial conditions:

$$D^2 g = F(x; g, Dg), \quad \beta_i g(a) = c_i, i = 1, 2, \tag{1.10}$$

then the linearised equation obtained through a Taylor expansion with parameter, $x$, about the point $\mathbf{g_0} = (g_0, Dg_0)$ is:

$$D^2 g \simeq F|_{\mathbf{g_0}} + \partial_{Dg} F|_{\mathbf{g_0}} (Dg - Dg_0) + \partial_g F|_{\mathbf{g_0}} (g - g_0),$$

$$D^2 g - \partial_{Dg} F|_{\mathbf{g_0}} Dg - \partial_g F|_{\mathbf{g_0}} g \simeq F|_{\mathbf{g_0}} - \partial_{Dg} F|_{\mathbf{g_0}} Dg_0 - \partial_g F|_{\mathbf{g_0}} g_0 .$$

This suggests an iterative scheme to repeatedly solve the associated linear collocation problem, above, in order to yield a pp solution for the ODE dependent upon an initial solution.

Indeed de Boor and Swartz [3, th. 3.1] state this explicitly as follows; there exists a unique collocation approximation $f_{\boldsymbol{\xi}}$, of equation 1.10 with equivalent side conditions, in the same neighbourhood of a solution, $g$, and Newton's method for approximating $f_{\boldsymbol{\xi}}$ converges quadratically in some neighbourhood of $f_{\boldsymbol{\xi}}$ for $|\boldsymbol{\xi}| \leq d$, where $|\boldsymbol{\xi}| = \max_i \{t_{i+1} - t_i\}$ and for some $d > 0$. Newton's method is given explicitly as:

$$D^2 f_{\boldsymbol{\xi},r+1} - \partial_{Dg} F|_{\mathbf{g_r}} D f_{\boldsymbol{\xi},r+1} - \partial_g F|_{\mathbf{g_r}} f_{\boldsymbol{\xi},r+1} = F|_{\mathbf{g_r}} - \partial_{Dg} F|_{\mathbf{g_r}} D f_{\boldsymbol{\xi},r} - \partial_g F|_{\mathbf{g_r}} f_{\boldsymbol{\xi},r}, \tag{1.11}$$

where $\mathbf{g_r} = (f_{\boldsymbol{\xi},r}, D f_{\boldsymbol{\xi},r})$, and $r = 0, 1, \dots$ .

Finally, and to be explicit, one must then observe that the collocation approximation, $f_{\boldsymbol{\xi},r+1}$, is obtained by employing the method outlined in section 1.2.9 for the linear system 1.7 /1.8, where we set:

$$\begin{aligned} p(x) &= -\partial_{Dg} F|_{\mathbf{g_r}}, \\ q(x) &= -\partial_g F|_{\mathbf{g_r}}, \\ s(x) &= F|_{\mathbf{g_r}} - \partial_{Dg} F|_{\mathbf{g_r}} D f_{\boldsymbol{\xi},r} - \partial_g F|_{\mathbf{g_r}} f_{\boldsymbol{\xi},r}, \end{aligned}$$

all of which are of course in terms of known quantities given the assumption of an initial supposed solution $f_{\boldsymbol{\xi},0}$.

### 1.2.11 Choice of collocation datasites, $\boldsymbol{\tau}$

The stipulation for there to be $(k-2)$ datasites within each interval has already been asserted and shown to produce the solvable linear system 1.7 /1.8 for a general nonlinear second order ODE in subsections 1.2.9 and 1.2.10. This far we have made no comment about the positioning of those datasites besides asserting that for the linear system to be solvable they must not coalesce. Continuing to suppose that $g(x)$ represents the solution of the general ODE then it is favourable to seek to position the datasites, $\boldsymbol{\tau}$, to minimise:

$$\|g - f_{\boldsymbol{\xi}}\|_\infty := \max_{a \leq t \leq b} |g(x) - f_{\boldsymbol{\xi}}(x)|.$$

Before presenting a theorem detailing the conditions under which this is achieved we first give an aside relevant to the application of said theorem.

#### Legendre polynomials and orthogonality

Two polynomials, $y(x)$ and $z(x)$, are said to be orthogonal to each other, with respect to weight function, $W(x)$, on an interval, $[a..b]$, if the following holds:

$$< y, z > := \int_a^b y(x) z(x) W(x) dx = 0$$

For $\mathbf{\Pi}_{<k}$ there exists a sequence of basis polynomials $(P_i)_{i=0}^{k-1}$, where deg $P_i = i$, that span the space, and for which the following holds:

$$< P_i, P_j > = K(i) \delta_{i,j},$$

The creation of such polynomials, with respect to the weight function, can be done iteratively starting with a basic polynomial, $P_0$ and then finding the coefficients of a polynomial of one order higher such that the definition holds, and repeating. The Gram-Schmidt process details this and is frequently adopted, commonly subject to $K(i) := 1$ in addition.

For our purposes we are interested in the property that $P_i$ is orthogonal to any polynomial $q(x) \in \mathbf{\Pi}_{<i}$. This is easily seen by noting $(P_j)_{j=0}^{i-1}$ is a basis and then $q(x) = \sum_{j=0}^{i-1} \alpha_j P_j$, for real $\alpha_j$, so that:

$$< q, P_i > = \sum_{j=0}^{i-1} \alpha_j < P_j, P_i > = 0.$$

The Legendre polynomials defined over the interval $[-1..1]$, and with respect to weight function $W(x) := 1$, are an important set of orthogonal polynomials. They may be expressed using Rodrigues' formula:

$$P_i(x) := \frac{1}{2^i i!} D^i \left[ (x^2 - 1)^i \right].$$

**Theorem on collocation datasite positioning[3, th. 4.1]**

Assume that the related linear problem 1.11 to equation 1.10 with linear side conditions is uniquely solvable, and assume that function $F$ in equation 1.10 is sufficiently smooth in a neighbourhood of the curve:

$$[a..b] \to \mathbb{R}^3 : x \longmapsto (x, g(x), Dg(x)).$$

Assume further that the collocation pattern, $\boldsymbol{\rho} = \{\rho_1, \rho_2, \ldots, \rho_{k-2}\}$ in the standard interval $[-1..1]$ is chosen so that:

$$\int_{-1}^{1} q(x) \prod_{j=1}^{k-2} (x - \rho_j) dx = 0, \tag{1.12}$$

for every $q(x) \in \mathbf{\Pi}_{<s}$. Then the collocation approximation, $f_{\boldsymbol{\xi}}$ (if it exists) near $g$ on $[a..b]$ of the associated problem 1.11 satisfies, globally:

$$\|D^m g - D^m f_{\boldsymbol{\xi}}\|_\infty \leq C|\boldsymbol{\xi}|^{(k-2)+\min(s,2-m)}, \quad m = 0, \ldots, 2, \tag{1.13}$$

where $C$ is a constant independent of $\boldsymbol{\xi}$.

Additionally at the breakpoints, $\boldsymbol{\xi}$, the error is of higher order, namely:

$$|D^m(g - f_{\boldsymbol{\xi}})(\xi_i)| \leq C|\boldsymbol{\xi}|^{(k-2)+s}, \quad m = 0, \ldots, 2, \;\; i = 1, \ldots, l+1. \tag{1.14}$$

Furthermore this is a superior error estimate to that which is specifed in [3, th. 3.1], relevant to section 1.2.10 in which a more general case of collocation datasite positioning is subsumed.

In the above theorem the collocation pattern, $\boldsymbol{\rho}$, defines the positioning of each of the $(k-2)$ datasites within each interval, $i$, in the following manner:

$$\tau_j := \frac{\xi_{i+1} + \xi_i}{2} + \rho_j \frac{\xi_{i+1} - \xi_i}{2}, \quad \in [\xi_i..\xi_{i+1}], \quad j = 1, \ldots, (k-2).$$

The entire collection of datasites can therefore be expressly written:

$$\tau_{j+\gamma} = \frac{\xi_{i+1} + \xi_i}{2} + \rho_j \frac{\xi_{i+1} - \xi_i}{2}, \quad i = 1, \ldots, l, \;\; j = 1, \ldots, (k-2).$$

By adopting, for $\boldsymbol{\rho}$, the zeroes of the Legendre polynomial $P_{k-2}$, we can be certain that for every $q(x) \in \mathbf{\Pi}_{<(k-2)}$ the equation 1.12 is satisfied, due to orthogonality, and thus we can set $s = (k-2)$ in inequalities 1.13 and 1.14. Due to the specificity of this theorem and the superiority of its resultant error estimate all succeeding sections of this study will adopt a collocation pattern, $\boldsymbol{\rho}$, equal to the zeroes of the Legendre polynomial, $P_{k-2}$, for the application of the collocation method to the Van der Pol equation with initial conditions.

## 1.3 The Van der Pol equation

### 1.3.1 Balthazar Van der Pol[11]

Balthazar Van der Pol, a Dutchman, died in 1959. He had been a pioneer in the field of radio and communications. His work and equations formed the basis of much of the modern theory of non-linear oscillations and his name is given to the most typical equation of the theory. While Van der Pol was working at Philips as an electrical engineer he found stable oscillations in electrical circuits employing vacuum tubes, which he termed relaxation-oscillations, and which are now known to be examples of limit cycles. The equation he subsequently derived has been studied in numerous different ways and is used to model the behaviour of a number of physical systems or effects. For instance as part of a model[12] for the action potential of neurons or to ascertain[13] displacement of tectonic plates across a geological fault.

The Van der Pol equation has been chosen specifically for this study, not least because of its physical applicability and therefore its practical use, which is something that appeals to me personally, but because it is hoped this study can add to the collective understanding of its solutions and create synergies for researchers analysing the cumulative works on the subject. Additionally the properties of the solutions discussed later in this section provide an excellent foundation on which to develop the ideas of collocation approximations proposed in this study.

### 1.3.2 The Van der Pol equation

The study of this dissertation is the numerical stability of the collocation method applied to the Van der Pol equation with varying parameters. Explicitly it is the investigation of the collocation approximation, $f_{\boldsymbol{\xi}}(x)$, to the solution, $g(x)$, of equation:

$$D^2g + \mu(g^2 - 1)Dg + g = 0, \qquad Dg(0) = 0, g(0) = 1, \qquad \mu > 0. \tag{1.15}$$

Going forward, any results or graphs depicted will adopt the above initial conditions unless otherwise specified.

Below we list a number of properties of the Van der Pol equation which direct and form part of the study.

### 1.3.3 Lack of analytic solutions

Using the substitution:

$$z_1 := g,$$
$$z_2 := Dg + \mu(\frac{1}{3}g^3 - g),$$

the Van der Pol equation can be expressed in terms of the Liénard system:

$$\begin{bmatrix} Dz_1 \\ Dz_2 \end{bmatrix} = \begin{bmatrix} z_2 - \mu(\frac{1}{3}z_1^3 - z_1) \\ -z_1 \end{bmatrix}.$$

Additionally since the original equation is an autonomous differential equation the substitution, $z = Dg$ leads the the first order differential equation:

$$z\frac{dz}{dg} + \mu(g^2 - 1)z + g = 0.$$

The above belongs to the class of Abel's equations of the second kind, and the lack of known (tabulated) functions of this equation infers there are no analytic solutions of the Van der Pol equation in terms of known (tabulated) functions[8].

### 1.3.4 Unique and stable limit cycle

The solution to the Van der Pol equation is a unique and stable limit cycle. The simplest way to demonstrate this is with the following theorem:

**Liénard's theorem[7]**

For;

$$D^2g + u(g)Dg + v(g) = 0,$$

where,

(i) $u, v$ are continuously differentiable functions on $\mathbb{R}$,

(ii) $u$ is an even function,

(iii) $v$ is an odd function, and $v > 0$ for all $x > 0$.

(iv) $\lim_{g\to\infty} U(g) := \lim_{g\to\infty} \int_0^g u(z)dz = \infty$.

(v) $U(g)$ has exactly one positive root at some value, $p$, where $U(g) < 0$ for $0 < g < p$ and $U(g) > 0$ and monotonic for $g > p$,

then the solution, $g(x)$, is a unique and stable limit cycle surrounding the origin.

Each of these conditions is satisfied by the Van der Pol equation where, $u(g) = \mu(g^2 - 1)$ and $v(g) = g$, so the result follows.

### 1.3.5 Nature of limit cycle and solutions

Figure 1.4 plots the unique limit cycle for varying parameter, $\mu$. For small $\mu$ the non-linear term is not dominant and the behaviour of the solution tends to simple harmonic oscillation. For larger values of $\mu$ the solution becomes more and more stiff, and for large values of $\mu$ the equation is sometimes referred to as a "relaxation oscillator" because the the stress built up over a longer time is relaxed in a sudden discharge. Equivalent to the phase plot, the solutions of the equation (with $Dg(0) = 0, g(0) = 2$) for select values, $\mu = \{0.05, 2, 10\}$, have been plotted in figure 1.5 to indicate the nature of the solution. Descriptively, as $\mu$ increases the period of the cycle increases and the graph becomes more and more 'saw-tooth'. For any $\mu$ the amplitude of oscillation always tends between $[-2..2]$, consistent with the phase plot.

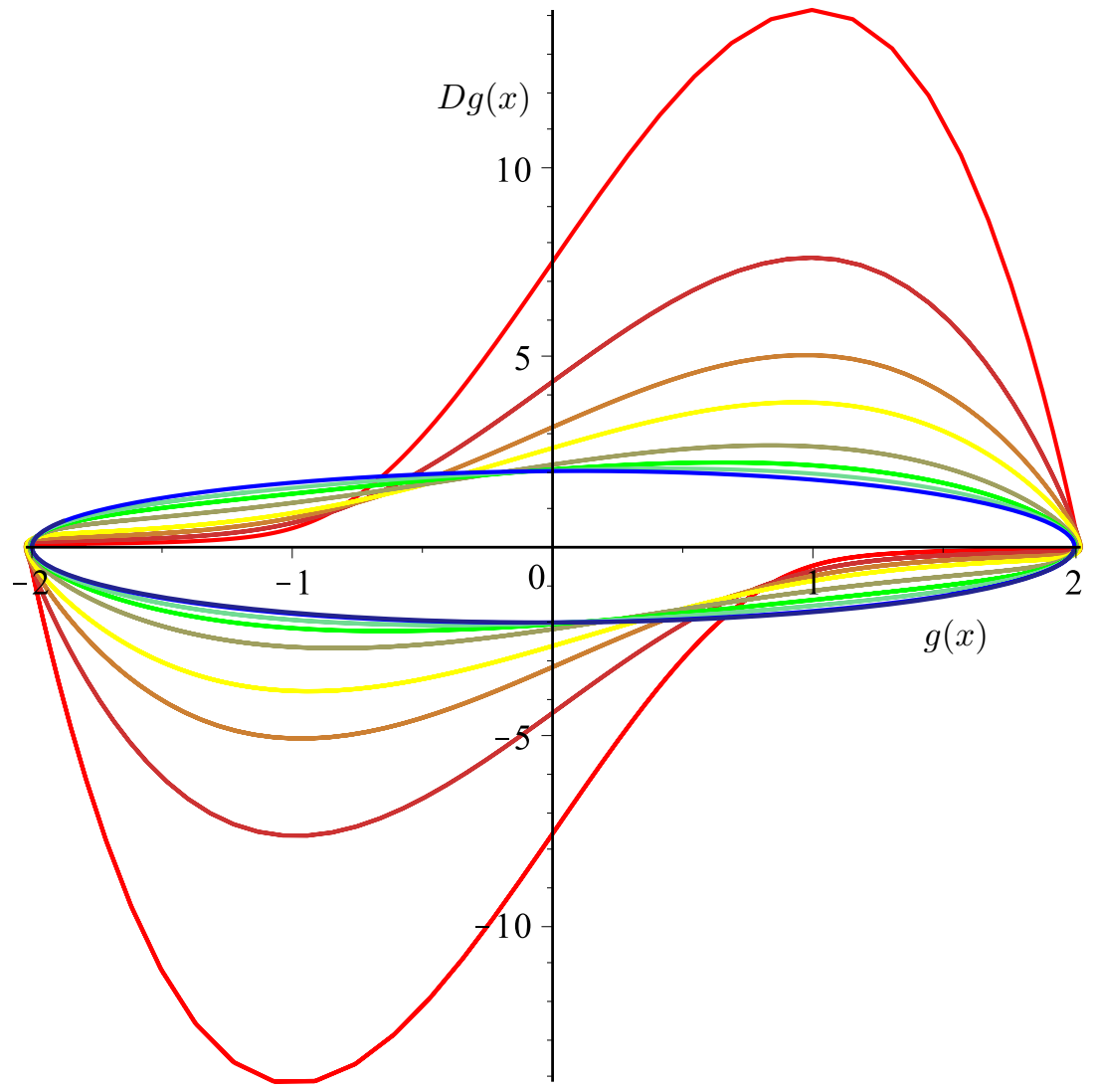

Figure 1.4: Traditional phase plot for $\mu = \{0.01, 0.05, 0.25, 0.5, 1, 2, 3, 5, 10\}$, where blue indicates the smallest $\mu$ and red the largest.

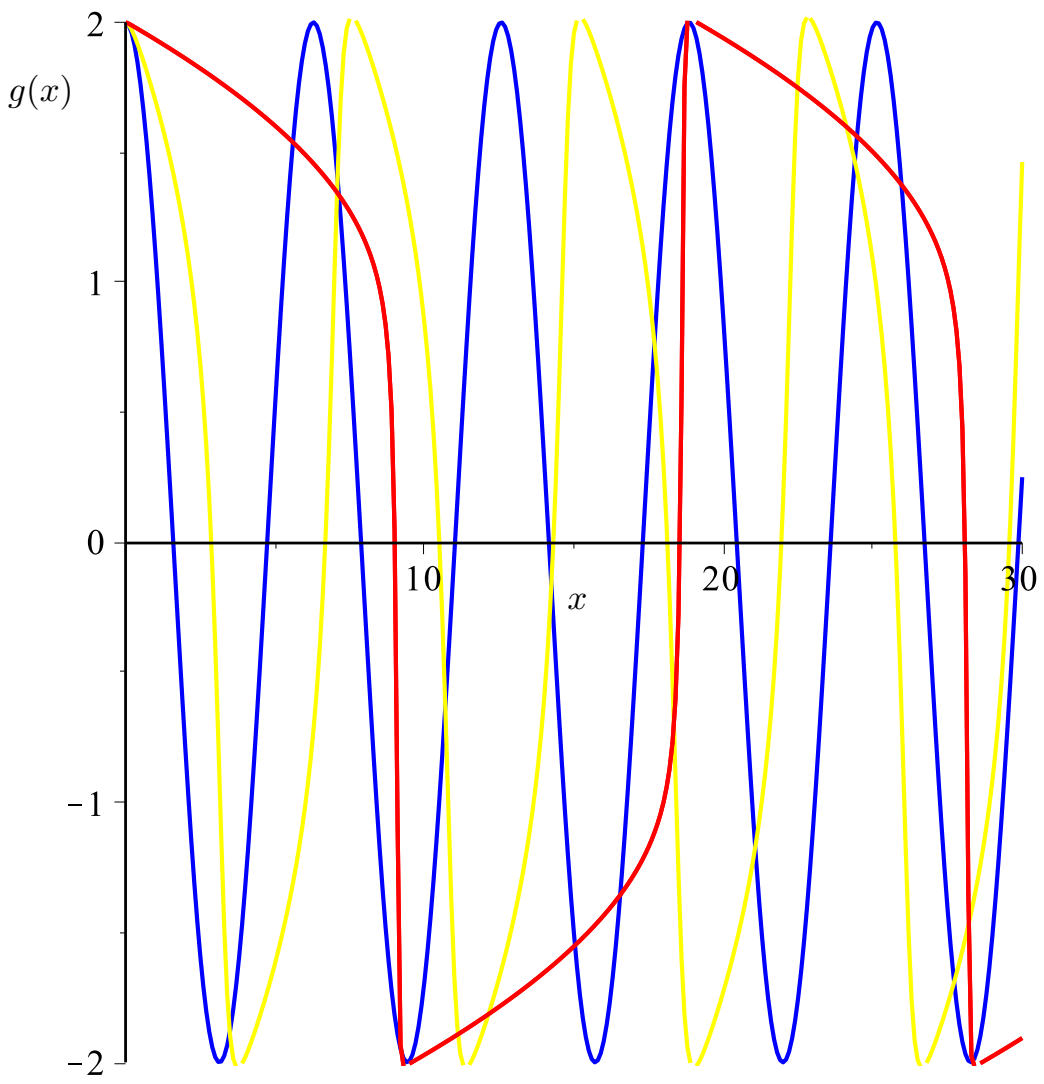

Figure 1.5: Solutions for $g(x)$ with $\mu = \{0.05, 2, 10\}$, where $Dg(0) = 0, g(0) = 2$.

### 1.3.6 Linearisation of the Van der Pol equation

Note that the initial conditions of the Van der Pol equation are already linear. Secondly, expressing the equation in the form:

$$D^2 g = F(g, Dg) = \mu(1 - g^2)Dg - g, \quad g(0) = 1, Dg(0) = 0, \tag{1.16}$$

then we observe that:

$$\begin{aligned} \partial_g F &= -2\mu g Dg - 1, \\ \partial_{Dg} F &= \mu(1 - g^2), \end{aligned}$$

so that by equation 1.11, the linearised form of the Van der Pol equation primed for iteration by Newton's method to the collocation approximation $f_{\boldsymbol{\xi}}$ to $g$ is:

$$D^2 f_{\boldsymbol{\xi},r+1} - \mu(1 - f_{\boldsymbol{\xi},r}^2)Df_{\boldsymbol{\xi},r+1} + (1 + 2\mu f_{\boldsymbol{\xi},r} Df_{\boldsymbol{\xi},r})f_{\boldsymbol{\xi},r+1} = \mu(1 - f_{\boldsymbol{\xi},r}^2)Df_{\boldsymbol{\xi},r} - f_{\boldsymbol{\xi},r}. \tag{1.17}$$

# Scope of study and methods

## 2.1 Practical considerations for the application of the collocation method

In this section we discuss; the application of the collocation method to generate approximations to the unique solution of the the Van der Pol equation 1.15, and the specific parameters and topics to analyse in the context of this study.

### 2.1.1 Properties of the equation and of the collocation method

#### Autonomous equations permit translations

The Van der Pol equation 1.15 is autonomous, such that the differential equation is invariant under the substitution $(x - t)$, and therefore permits translations in $x$ for unconditioned solutions.

#### Collocation approximations are unequal for different given conditions

A unique solution $g(x)$ exists to equation 1.15, and, specific to the initial conditions, there is a uniquely associated point of this solution, $g(b) = g_b$. Equivalently, given the boundary conditions $\{g(0) = 1, g(b) = g_b\}$, to the differential equation, the solution, $g(x)$ remains the same and has the associated derivative, $Dg(0) = 0$.

A collocation approximation $f_{\boldsymbol{\xi}}$ to $g$ does not have the same property. That is;

(i) for a collocation approximation, $f_{\boldsymbol{\xi}}^{ICs}$, based on initial conditions, we have;

$$f_{\boldsymbol{\xi}}^{ICs}(0) = g(0), \quad Df_{\boldsymbol{\xi}}^{ICs}(0) = Dg(0), \quad \text{but } \textit{not} \text{ necessarily,} \quad f_{\boldsymbol{\xi}}^{ICs}(b) = g(b), \quad Df_{\boldsymbol{\xi}}^{ICs}(b) = Dg(b),$$

(ii) for a collocation approximation, $f_{\boldsymbol{\xi}}^{BCs}$, based on boundary conditions, we have;

$$f_{\boldsymbol{\xi}}^{BCs}(0) = g(0), \quad f_{\boldsymbol{\xi}}^{BCs}(b) = g(b), \quad \text{but } \textit{not} \text{ necessarily,} \quad Df_{\boldsymbol{\xi}}^{BCs}(0) = Dg(0), \quad Df_{\boldsymbol{\xi}}^{BCs}(b) = Dg(b).$$

Although this property of potentially unequal collocation approximations is true in general for a change in supplied datasites, which leads to the theorem on collocation datasite positioning in section 1.2.11, we have highlighted it here in specific detail for the cases of initial conditions compared to boundary conditions supplied to the differential equation.

**Collocation approximations evolve in a left-to-right manner for initial value problems**

Superficially the solution, $g$, to equation 1.15 is better understood at $x = 0$, and becomes more uncertain towards $x = b$. The initial conditions are known and one can estimate the nature of the solution about the initial point but with increasing uncertainty the further one deviates from it. In an analogous manner this is how the collocation approximation evolves, becoming more accurate at rightmost values only after successful iterations.

This is made explicit by analysing B-spline coefficient dependence. Considering briefly the previous example equation 1.9 the inverse of matrix $\mathbf{A}$, which determines the collocation spline coefficients, in terms of non-zero elements can be shown to be:

$$\begin{bmatrix} X & X & 0 & 0 & 0 & 0 & 0 & 0 \\ X & X & 0 & 0 & 0 & 0 & 0 & 0 \\ X & X & X & X & 0 & 0 & 0 & 0 \\ X & X & X & X & 0 & 0 & 0 & 0 \\ 0 & 0 & X & X & X & X & 0 & 0 \\ 0 & 0 & X & X & X & X & 0 & 0 \\ 0 & 0 & 0 & 0 & X & X & X & X \\ 0 & 0 & 0 & 0 & X & X & X & X \end{bmatrix}^{\mathbf{-1}} = \begin{bmatrix} X & X & 0 & 0 & 0 & 0 & 0 & 0 \\ X & X & 0 & 0 & 0 & 0 & 0 & 0 \\ X & X & X & X & 0 & 0 & 0 & 0 \\ X & X & X & X & 0 & 0 & 0 & 0 \\ X & X & X & X & X & X & 0 & 0 \\ X & X & X & X & X & X & 0 & 0 \\ X & X & X & X & X & X & X & X \\ X & X & X & X & X & X & X & X \end{bmatrix} .$$

The important property, in generality, being the upper triangular pattern of zeroes can be shown to be maintained in the inverse for $\mathbf{A}$ for any permissible $k$, $l$. The specific dependence of the elements in the inverse matrix are defined in greater detail in section 7.1. In turn this allows for the statement of dependence of the spline coefficients, $\alpha_i$, in this example case, where we denote by $A_{x,y,z}$ and $s_{x,y,z}$ the rows $x, y, z$ of matrices $\mathbf{A}$ and $\mathbf{s}$ respectively:

$$\begin{bmatrix} \alpha_1 \\ \alpha_2 \\ \alpha_3 \\ \alpha_4 \\ \alpha_5 \\ \alpha_6 \\ \alpha_7 \\ \alpha_8 \end{bmatrix} = \begin{bmatrix} \alpha_1(A_{1,2}, s_{1,2}) \\ \alpha_2(A_{1,2}, s_{1,2}) \\ \alpha_3(A_{1,..,4}, s_{1,..,4}) \\ \alpha_4(A_{1,..,4}, s_{1,..,4}) \\ \alpha_5(A_{1,..,6}, s_{1,..,6}) \\ \alpha_6(A_{1,..,6}, s_{1,..,6}) \\ \alpha_7(A_{1,..,8}, s_{1,..,8}) \\ \alpha_8(A_{1,..,8}, s_{1,..,8}) \end{bmatrix} ,$$

and in general, taken again from section 7.1, the dependence for any B-spline coefficient is observed to be as follows:

$$\alpha_i = \alpha_i\Big(A_{1,..,m+2}, s_{1,..,m+2}\Big), \quad m := \left(\text{floor}\left(\frac{i-3}{k-2}\right) + 1\right)(k-2).$$

From the support and positivity property of B-splines in equation 1.2, rightmost values of the collocation approximation are determined by coefficients with higher $i$, and are therefore dependent upon the accurate convergence of approximation values to the left, proving the approximation can only evolve in this manner. As an example, figure 2.1 demonstrates this rather nicely.

Furthermore, because of this dependence, the collocation approximation actually evolves in a manner that is akin to applying the collocation method interval by interval. Once an interval, $i$ say, has converged on $[\xi_i..\xi_{i+1}]$ then all the B-spline coefficients for those B-splines lending support in that interval are determined. This allows the determination of $f_{\boldsymbol{\xi}}(\xi_{i+1})$ and $Df_{\boldsymbol{\xi}}(\xi_{i+1})$ from precisely two coefficients, since at any interior breakpoints there are only ever two B-splines which lend support to function values or derivatives;

$$f_{\boldsymbol{\xi}}(\xi_{i+1}) = \alpha_{k-1+\gamma} B_{k-1+\gamma,k}(\xi_{i+1}) + \alpha_{k+\gamma} B_{k+\gamma,k}(\xi_{i+1}) \ ,$$
$$Df_{\boldsymbol{\xi}}(\xi_{i+1}) = \alpha_{k-1+\gamma} DB_{k-1+\gamma,k}(\xi_{i+1}) + \alpha_{k+\gamma} DB_{k+\gamma,k}(\xi_{i+1}) \ .$$

The collocation approximation continues to compute, essentially inheriting these values as the new initial conditions (albeit strictly speaking these are imposed continuity conditions) for solving the remaining $(k-2)$ coefficients relevant to the interval $[\xi_{i+1}..\xi_{i+2}]$.

By explanation in a practical sense for the continued example equation 1.9;

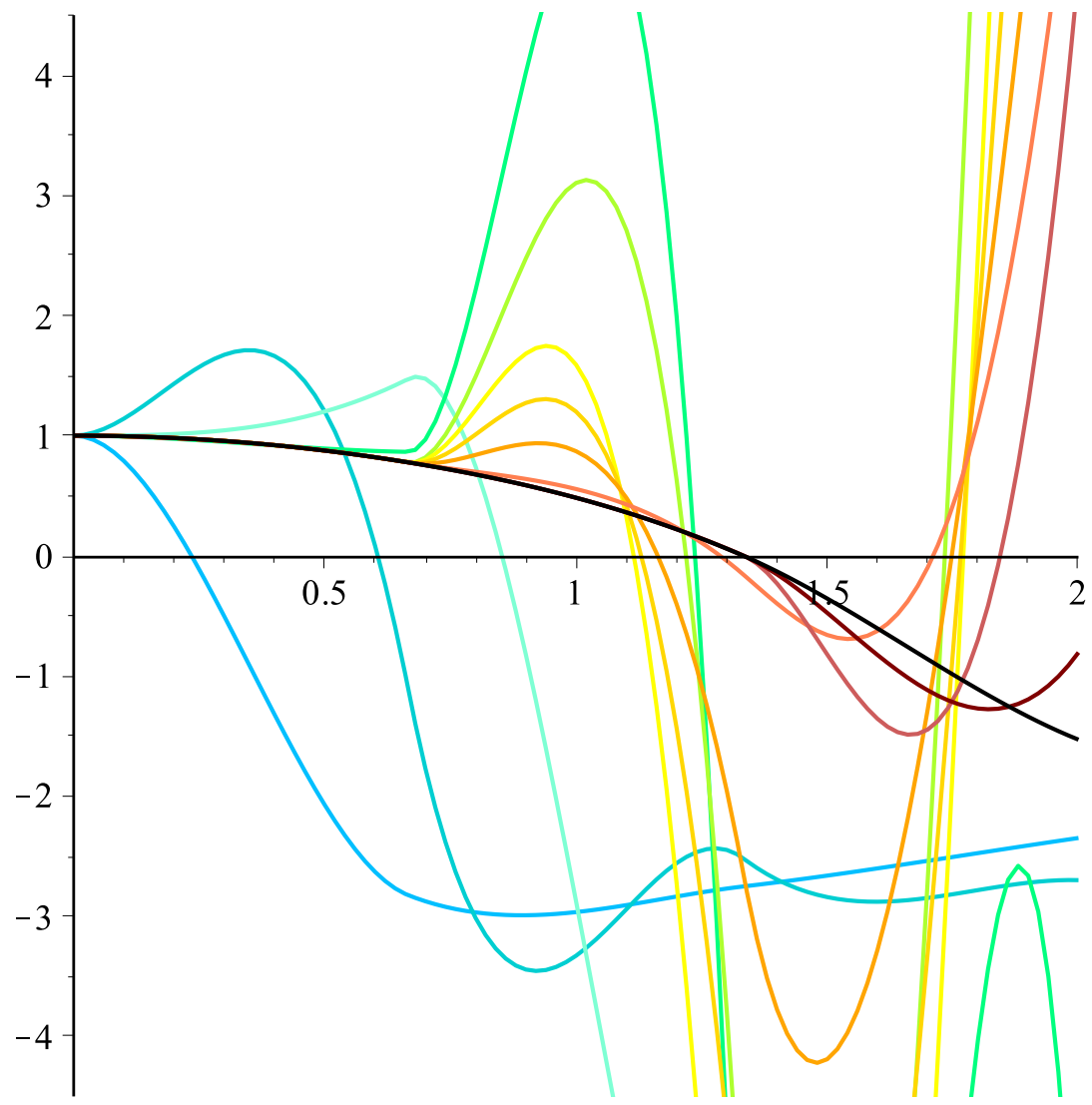

Figure 2.1: Plot for $\{f_{\boldsymbol{\xi},r} : r = 1, \ldots, 12\}$, where $\boldsymbol{\xi} = \{0, \frac{2}{3}, \frac{4}{3}, 2\}$, $\mu = 1.33$, showing convergence from blue, through green, yellow, red, to black, in a left-to-right manner.

(i) the coefficients, $\alpha_1, \alpha_2$ are determined by the initial conditions at $x = 0$,

(ii) the coefficients, $\alpha_3, \alpha_4$ are then determined through collocation approximation converging, such that the differential equation is satisfied at datasites, $\tau_1, \tau_2$. The pp approximation on $[0..\xi_2]$ is then established,

(iii) the coefficients, $\alpha_5, \alpha_6$ are determined in the same manner for datasites, $\tau_3, \tau_4$, subject to the new initial conditions from (ii) of $f_{\boldsymbol{\xi}}(\xi_2)$ and $Df_{\boldsymbol{\xi}}(\xi_2)$. This establishes the pp on $[\xi_2..\xi_3]$.

(iv) $\alpha_7, \alpha_8$ are determined in the same manner as (iii) appropriate to an increase by one in $i$, to produce the complete $f_{\boldsymbol{\xi}}$.

### 2.1.2 Measure of collocation approximation error

#### Considerations to error

The global error estimate provided by the theorem on collocation datasite positioning in section 1.2.11 outlines constant, $C$, scaling the upper bound, to be independent of $|\boldsymbol{\xi}|$. However we assert that for initial value problems $C$ is indirectly dependent upon choice of $\boldsymbol{\xi}$. We note in passing that in their article[3] the authors comment they "cannot recommend [initial value problems] unconditionally without further analysis" which they claim is due to the similarity of the scheme to that of Loscalzo and Talbot[9], which can be numerically unstable.

We observe that collocation approximations obtained through boundary value conditions supplied to equation 1.15 on $[0..b]$ will always be constrained with equality at each side, and there is no direct or indirect dependence of $C$ on $|\boldsymbol{\xi}|$. Though, for initial value problems that are not rightmost constrained one would naturally expect for an arbitrary range $[0..b]$ that the global error estimate of an approximation to increase with $b$.

The combination of all three above properties of section 2.1.1 lead to the following assertions for collocation approximations with supplied initial conditions;

(i) after convergence of the first or any subsequent interval in a left-to-right manner, it is *not* necessarily true that;

$$f_{\boldsymbol{\xi}}(\xi_i) = g(\xi_i), \quad Df_{\boldsymbol{\xi}}(\xi_i) = Dg(\xi_i) \ ,$$

(ii) this then inherited error from previous intervals can compound to produce greater uncertainty through calculation of subsequent intervals and the further one deviates from the point $x = 0$,

(iii) the collocation approximation can become translated in $x$ against the real solution $g$, but satisfy the differential equation fully in a particular interval.

(iv) each interval can then accumulate error. The global error becomes indirectly dependent upon the total number of intervals of $\xi$, and also the error assumed in each interval.

An equivalent mathematical argument of this concept is shown in section 7.2.

**Establishing an error measure for numerical results**

We acknowledge there are at least two properties of an approximation one might analyse;

(i) the approximation's, $f_{\boldsymbol{\xi}}$, error when compared to the real solution, $g$,

(ii) the error in the differential equation evaluated at $f_{\boldsymbol{\xi}}(x)$.

We define the error in the differential equation;

$$err(x) := D^2 f_{\boldsymbol{\xi}}(x) + \mu(f_{\boldsymbol{\xi}}(x)^2 - 1)Df_{\boldsymbol{\xi}}(x) + f_{\boldsymbol{\xi}}(x) \ .$$

For $g$ and $f_{\boldsymbol{\xi}}$ continuous we observe that:

$$\begin{aligned} \|g - f_{\boldsymbol{\xi}}\|_\infty \to 0, \quad &\implies \quad f_{\boldsymbol{\xi}}(x) \to g(x) \ , \\ \|err(x)\|_\infty \to 0, \quad &\implies \quad f_{\boldsymbol{\xi}}(x) \to g(x) \ . \end{aligned}$$

From section 1.3.3 the real solution is unknown so (i) is numerically indeterminable, so we focus in this study on (ii) as our measure, which is numerically determinable. We remain vigilant of the mentioned instances where this can be misleading, and look to sample many thousands of $x$-values in the solution range in an attempt to ensure much larger errors, existing precisely at unsampled $x$-values, are not misreported.

Adopting an assumed traditional boundary value collocation approach we derive, in section 7.3, the bound for $err(x)$ to characterise an inherent relationship between $\|err(x)\|_\infty$ and $\|g - f_{\boldsymbol{\xi}}\|_\infty$ in that sense:

$$\|err(x)\|_\infty \le C\big(|\boldsymbol{\xi}|^{(k-2)} + |\boldsymbol{\xi}|^k + \mu(3|\boldsymbol{\xi}|^{(k-1)} + \frac{1}{4}|\boldsymbol{\xi}|^{3k-1})\big) + 4\mu^2|\boldsymbol{\xi}|^k \ . \tag{2.1}$$

### 2.1.3 Measure of collocation approximation reliability

As well as perceived approximation error it is useful to have a concept of the cost of an approximation in terms of the amount of work expended in achieving the approximation. We have developed this estimate theoretically in terms of a count of floating point operations (flop count), but is equivalent to a CPU time estimate. The appendix details the derivation but we state here that a cost function associated with an optimally efficient algorithm for finding a collocation approximation is:

$$cost(f_{\boldsymbol{\xi}}) = O((k^3 + k^2 + k)Nl), \tag{2.2}$$

where $N$ represents the number of iterations required for the approximation to converge to the specified tolerance.

### 2.1.4 Initial supposed solution

The choice of initial supposed solution affects the number of iterations it will take before a collocation approximation spanning the range converges to a given tolerance. We prefer to promote the assumption that accuracy is unlikely to be ascertained at the outset, so we adopt the following approach to satisfy the known quantities and properties;

(i) for any initial point, $a$ say, the initial supposed solution and its derivative, $f_{\boldsymbol{\xi},0}(a)$ and $Df_{\boldsymbol{\xi},0}(a)$, should be with a small neighbourhood of $g(a)$ and $Dg(a)$,

(ii) the initial supposed solution begins as a straight line and collapses to zero over an appropriate number of units, to limit uncontrolled initial behaviour.

The precise definition of the initial supposed solution adopted in this study is the following, let:

$$f_0(x) := (g(a) + Dg(a)(x-a))\left(\frac{1 - tanh(x - a - 3)}{2}\right) \ ,$$

then the initial supposed solution, $f_{\boldsymbol{\xi},0}$ is the spline interpolation of this function:

$$f_{\boldsymbol{\xi},0}(x) \in \$_{k,\mathbf{t}} : f_{\boldsymbol{\xi},0}(a) = f_0(a), f_{\boldsymbol{\xi},0}(b) = f_0(b), f_{\boldsymbol{\xi},0}(\tau_i) = f_0(\tau_i) \forall i = 1, \dots, (k-2)l.$$

## 2.2 Parameters of study

In this section we outline the aspects of study, in terms of which quantities we seek to vary and which to hold constant, and to what extent.

### $\mu$ and $b$

The periodicity and stiffness of the solution is determined by parameter, $\mu$. A key question posed is whether large values of $\mu$ lead to unstable, or unreliable solutions, in the sense their associated cost functions are too great.

As well, the solution range must display enough information about the solution to be thorough. The parameter $b$ is typically chosen so that two periodic cycles are broadly attained within the range, or if convergence to a limit cycle is slow (very small $\mu$) or fast (large $\mu$) then we include more or less cycles respectively. Based on the values in table 2.1, we adopt, to allow direct comparison, a constant range [0..40] except occasionally where $\mu = 40$, and we sometimes extend the range.

| $\mu$ | $b$ |
|---|---|
| 0.01 | 40 |
| 0.05 | 40 |
| 0.25 | 30 |
| 0.5 | 25 |
| 1.0 | 20 |
| 2.0 | 20 |
| 3.0 | 20 |
| 5.0 | 25 |
| 10.0 | 40 |
| 20.0 | 40 |
| 40.0 | 70 |

Table 2.1: Appropiate solution range $[0..b]$ for varying $\mu$.

### $\boldsymbol{\xi}$

In this study only a uniform breakpoint sequence is used but varied in terms of the number of intervals. $l$.

### $\boldsymbol{k}$

Most results are produced for typical, lower values of $k$ such as 4,5 or 6, but some numerical results are recorded with values of $k$ upto 11.

### $\boldsymbol{\rho}$

Section 1.2.11 outlines the choice of this study to keep $\boldsymbol{\rho}$ constant.

### $\boldsymbol{f_{\xi,0}}$

As outlined in section 2.1.4 a constant form of initial solution is maintained in all results.

### *tol*

As outlined in section 2.3 a tolerance is required to ascertain the assessed convergence of the collocation approximation. We adopt $tol = 10^{-4}$. Given the final three iterated approximations are conditioned and that convergence is quadratic then we assume this to be sufficient and expect the convergence to be superior to that specified.

## 2.3 Numerical procedures

Maple 16 is the software used to produce the numerical results. The procedures created and used can be seen in section 7.7, and are available in electronic format as a Maple worksheet. As an overview, the main procedure is designed to complete the following tasks;

(i) allow inputs to control and parametrise the collocation approximation,

(ii) create the vectors, $\boldsymbol{\xi}, \boldsymbol{\tau}, \mathbf{t}, \boldsymbol{\rho}$ according to inputs,

(iii) for each, $w$ [1], create $f_{\boldsymbol{\xi},0}$ from the specified function with a spline interpolation procedure,

(iv) before Newton iteration, evaluate, to an array, all B-spline values and their derivatives at each required collocation datasite,

(v) perform Newton iteration, firstly creating matrices $\mathbf{A}$ and $\mathbf{s}$ in a block structure,

(vi) then using Maple LinearSolve(method=SparseDirect) to calculate the B-spline coefficients, $\boldsymbol{\alpha}$,

(vii) repeat for a maximal number of iterations or until a tolerance has been met that satisfies the following condition:

$$|f_{\boldsymbol{\xi},r}(b*) - f_{\boldsymbol{\xi},r-1}(b*)| < tol \quad \text{and} \quad |f_{\boldsymbol{\xi},r}(b*) - f_{\boldsymbol{\xi},r-2}(b*)| < tol,$$

where $b*$ is the rightmost point for the specific $w$ in scope,

(viii) collate data for use in plots and other output,

(ix) repeat for each $w$, if $w > 1$,

(x) finally display results.

[1] $w$ is introduced in section 4.1

# Preliminary results

## 3.1 Verification of $cost(f_{\boldsymbol{\xi}})$ estimate

In section 2.1.3 we developed a theoretical estimate of the cost of generating collocation approximations. Table 3.1 displays numerical results in terms of the number of seconds taken for the CPU to complete the given task, which we equivocate to a flop count. We conclude that equation 2.2 is likely an accurate reflection on an efficient algorithm, however, our results are in fact consistent with a method of slightly less efficiency, that of flop count of $O((k^3 + k^2 + k)Nl^{\frac{4}{3}})$, see figure 3.1. The suspected reason for this is the inclusion of Maple's LinearSolve subroutine using a sparsedirect solve-method. Indeed a note on sparsedirect solvers[10] suggests this order is typical across some of those systems. This certainly demonstrates that in conjunction with $\boldsymbol{\xi}$; $k, l, N$ must be considered in the context of generating a collocation approximation with an appropriate accuracy for a feasible cost.

| $l/N$ | 10 | 100 | 1000 |
|---|---|---|---|
| 20 | 0.09 | 0.67 | 7.92 |
| 40 | 0.16 | 1.99 | 15.4 |
| 80 | 0.32 | 3.32 | 32.0 |
| 160 | 0.70 | 6.82 | 70.3 |
| 320 | 1.64 | 17.6 | 173 |

(a) $k = 4$

| $l/N$ | 10 | 100 | 1000 |
|---|---|---|---|
| 20 | 0.15 | 1.31 | 15.4 |
| 40 | 0.30 | 3.27 | 30.5 |
| 80 | 0.63 | 6.75 | 64.3 |
| 160 | 2.02 | 14.5 | 124 |
| 320 | 4.04 | 36.7 | 359 |

(b) $k = 5$

| $l/N$ | 10 | 100 | 1000 |
|---|---|---|---|
| 20 | 0.26 | 2.28 | 25.8 |
| 40 | 0.52 | 5.61 | 51.4 |
| 80 | 1.10 | 10.7 | 112 |
| 160 | 3.10 | 25.3 | 253 |
| 320 | 7.35 | 62.2 | 629 |

(c) $k = 6$

| $l/N$ | 10 | 100 | 1000 |
|---|---|---|---|
| 20 | 0.44 | 3.57 | 40.2 |
| 40 | 0.82 | 8.28 | 82.0 |
| 80 | 1.75 | 17.7 | 175 |
| 160 | 4.91 | 39.3 | 397 |
| 320 | 13.3 | 101 | 999 |

(d) $k = 7$

Table 3.1: $cost(f_{\boldsymbol{\xi}})$ in terms of CPU seconds to complete $N$ iterations for a given, $l$ and $k$, with $\mu = 1.0$ on [0..40].

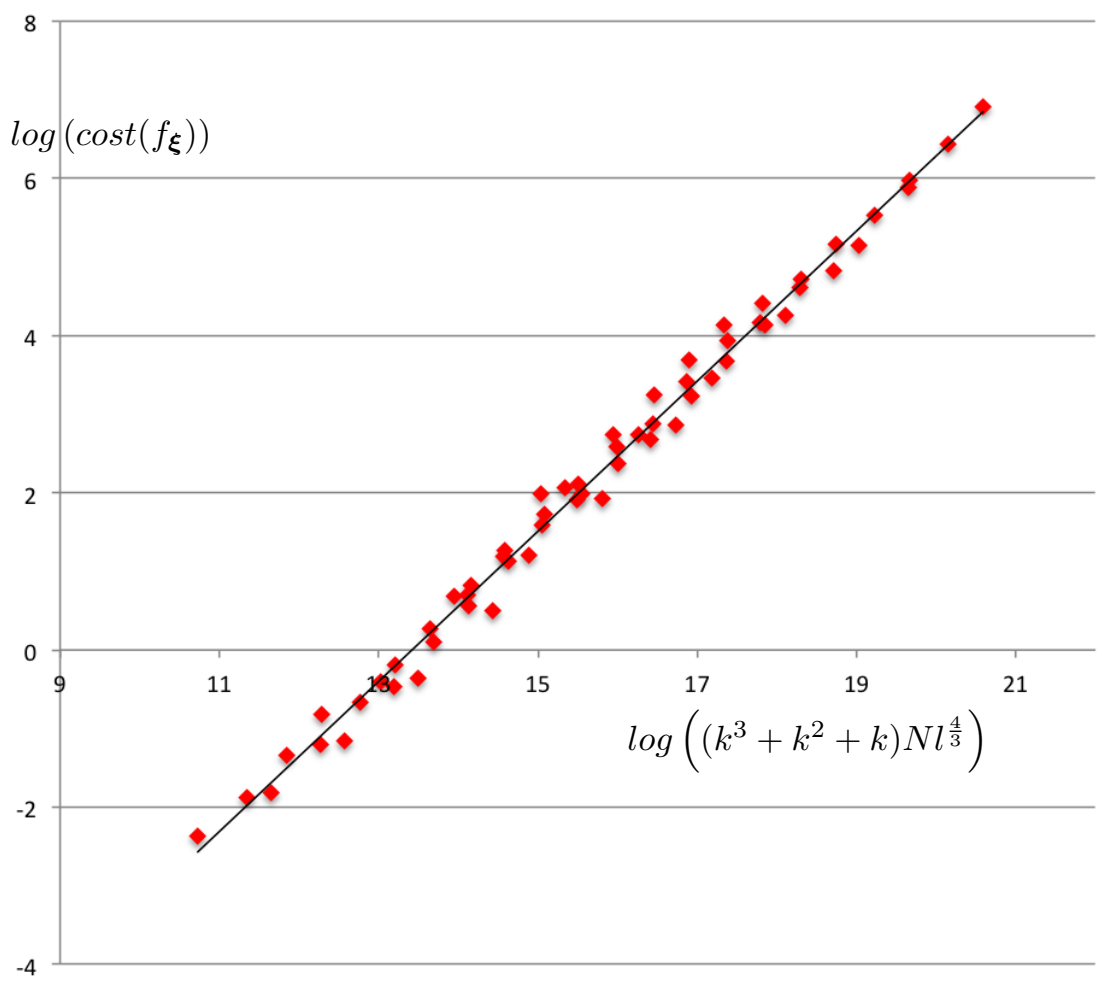

Figure 3.1: $cost(f_{\boldsymbol{\xi}})$ plotted against $(k^3+k^2+k)Nl^{\frac{4}{3}}$ using data from table 3.1.

## 3.2 Small values of $\mu$, with $k = 5$

### 3.2.1 Number of iterations to converge

See table 3.2.

| $l$ / $\mu$ | 0.01 | 0.05 | 0.25 | 0.5 | 1.0 | 2.0 |
|---|---|---|---|---|---|---|
| 5 | ** | ** | ** | | | |
| 10 | 4 | 6 | 76* | ** | ** | |
| 20 | 4 | 6 | 82 | 149* | 285* | >1000* |
| 40 | 4 | 6 | 15 | 164 | 288 | 821 |
| 80 | 4 | 6 | 15 | 152 | 379 | 791 |
| 140 | 4 | 6 | 15 | 28 | 641 | 109 |
| 160 | 4 | 6 | 15 | 82 | 744 | 109 |
| 200 | | 6 | 15 | 28 | >1000 | 109 |
| 300 | | | | 28 | 819 | 109 |
| 400 | | | | | >1000 | 109 |

Table 3.2: Number of iterations for $f_{\boldsymbol{\xi},r}(b)$ to converge to a tolerance of $10^{-4}$ for varying parameters, $\mu$, and $l$. (*poor approximation, **unacceptable approximation)

### 3.2.2 Error estimates for converged approximations

See figures 3.2, 3.3, 3.4.

### 3.2.3 Preliminary conclusions

We infer from the results in table 3.2 that for small values of $\mu$ the collocation method and resulting approximation is viable in terms of cost, particularly for $\mu \leq 0.25$. It is difficult to discern any pattern, particularly a useful or reliable one, for predicting the number of iterations, in the case where $\mu \geq 0.5$, for the approximation to converge. However, for small $\mu$ the approximations do converge for a typically acceptable cost. It was observed in passing that varying the initial supposed solution also changes the iteration path so could be adopted in an attempt to reduce the number of required iterations, but we do not comment further on this.

Figures 3.2 to 3.4 depicting $err(x)$ serve to demonstrate, along with table 3.2, that as $\mu$ increases not only does $cost(f_{\boldsymbol{\xi}})$ increase due to the requirement to perform more iterations, but also that the approximation suffers from increased error in the differential equation. The error can be somewhat diminished by increasing $l$, again as shown in the graphs, but will of course increase $cost(f_{\boldsymbol{\xi}})$ too. The number of iterations can decrease as $l$ increases, as a mitigating factor, but this is not true in general and cannot be relied upon. Additionally, increasing $k$ is another approach that might be adopted to reduce error, but again this has an associated cost and indeed may to lead to more or less iterations required.

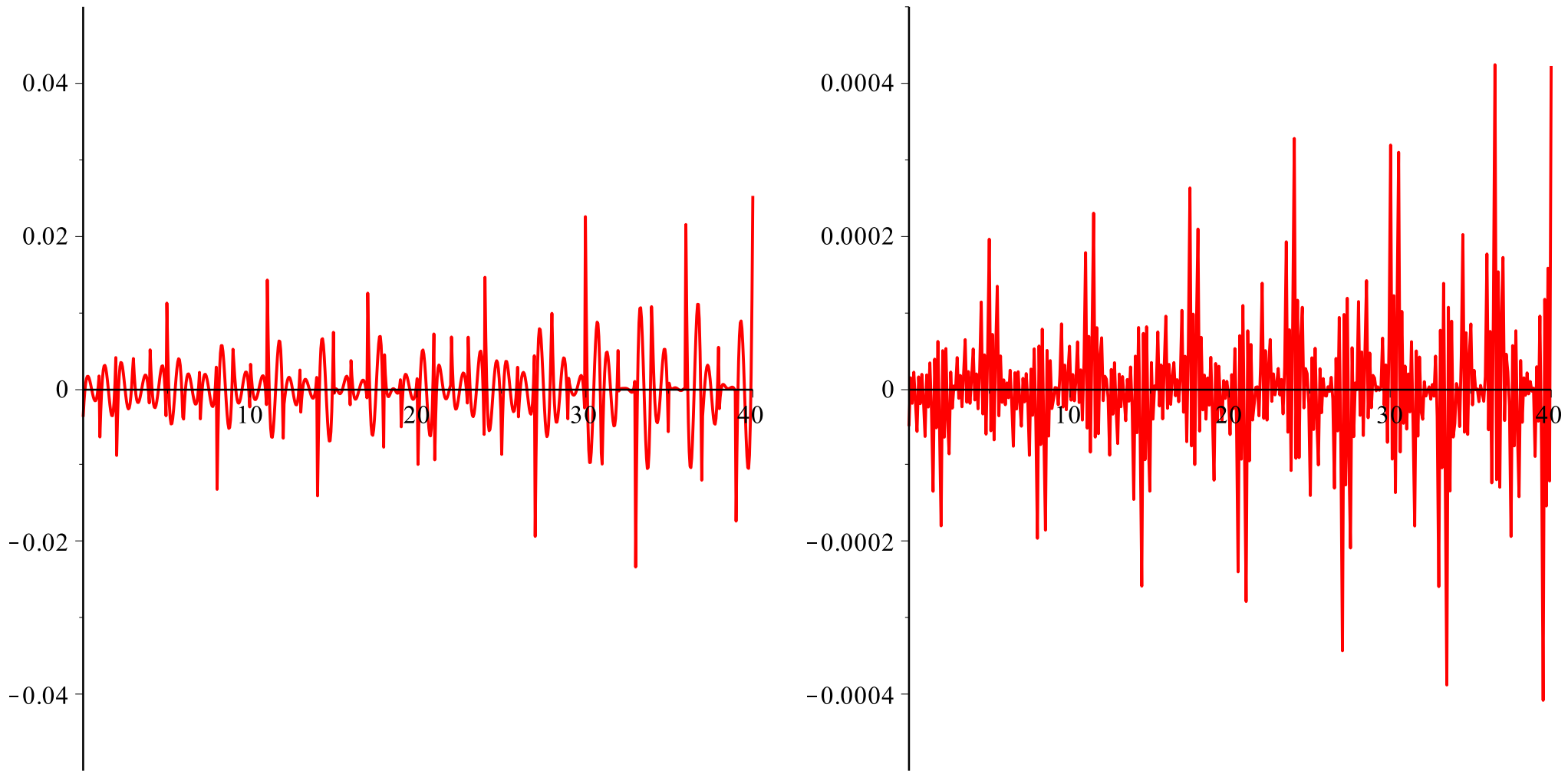

(a) $l = 40$ (b) $l = 160$

Figure 3.2: $err(x)$ for $\mu = 0.05$ after the given number of iterations in table 3.2

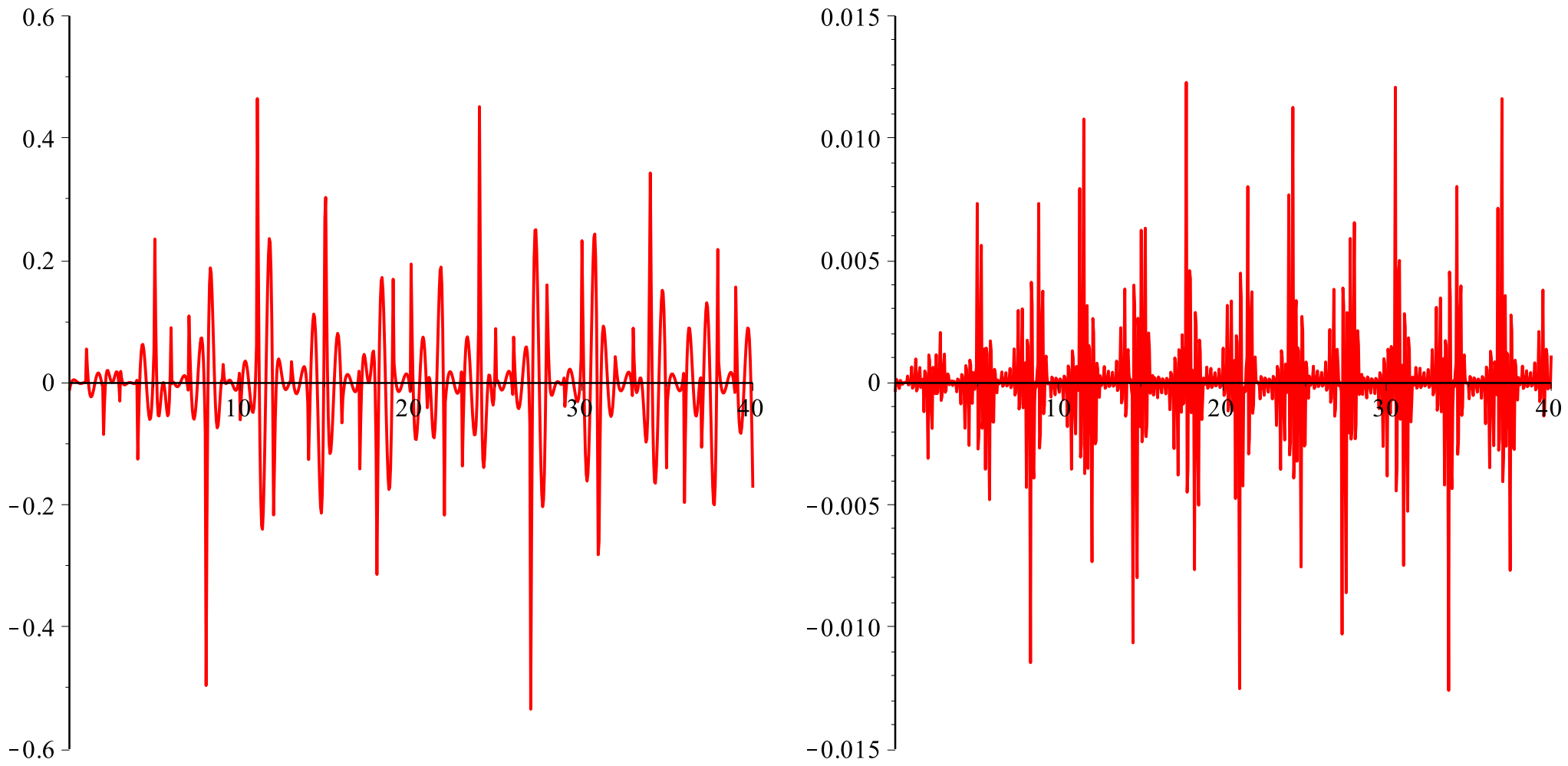

(a) $l = 40$ (b) $l = 160$

Figure 3.3: $err(x)$ for $\mu = 0.5$ after the given number of iterations in table 3.2

## 3.3 Large values of $\mu$, with $k = 5$

### 3.3.1 Collocation approximations and error estimates

See figures 3.5, 3.6, 3.7.

### 3.3.2 Preliminary conclusions

Figures 3.5, 3.6 and 3.7 show the viable collocation approximations for $\mu = \{3, 5, 10\}$. The costs of generating these solutions are quite high with convergence after 1138, 1969 and 5606 iterations respectively, where the number of intervals, $l = 160$. The number of iterations and therefore cost required to produce

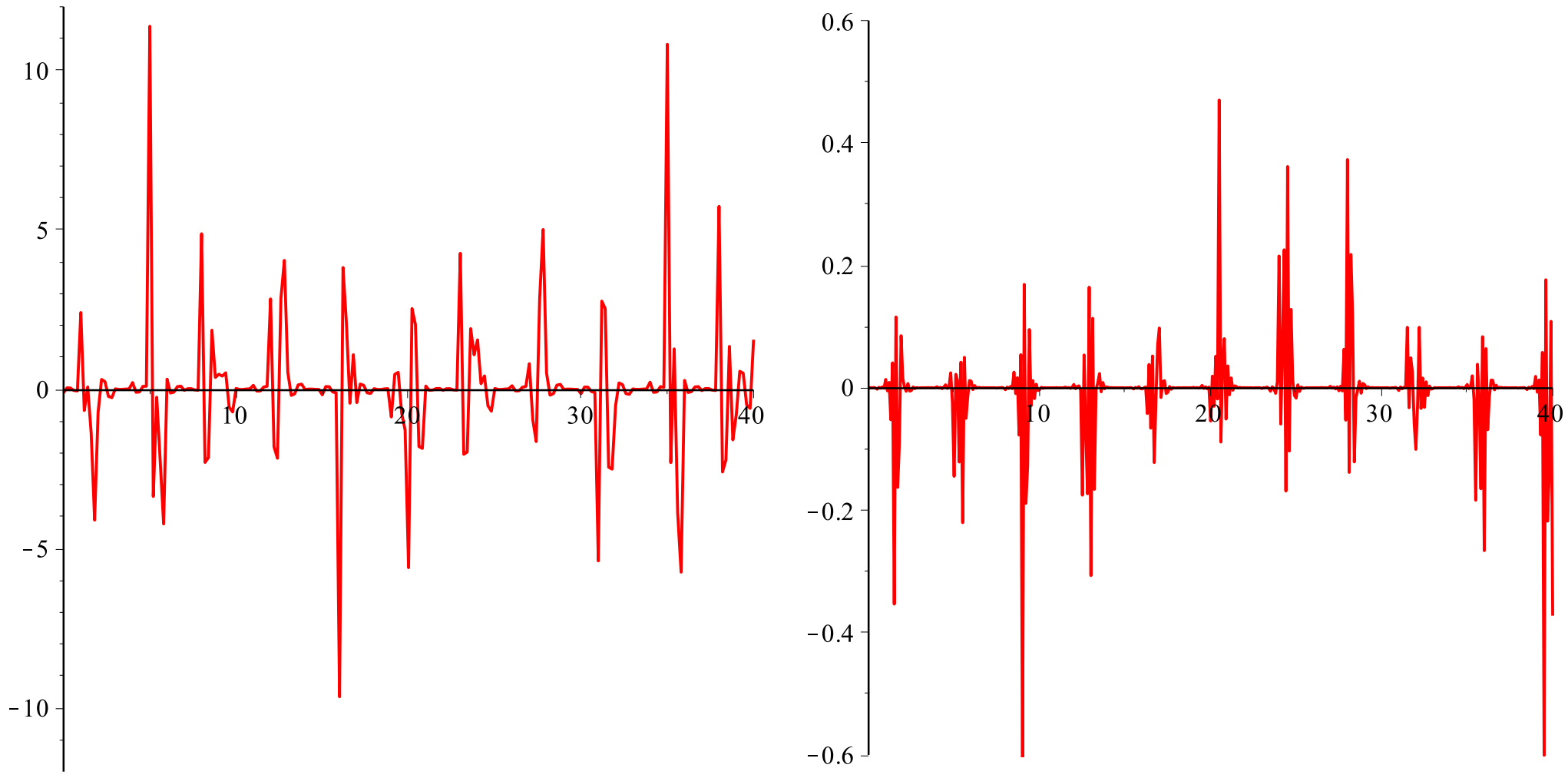

(a) $l = 40$ (b) $l = 160$

Figure 3.4: $err(x)$ for $\mu = 2.0$ after the given number of iterations in table 3.2

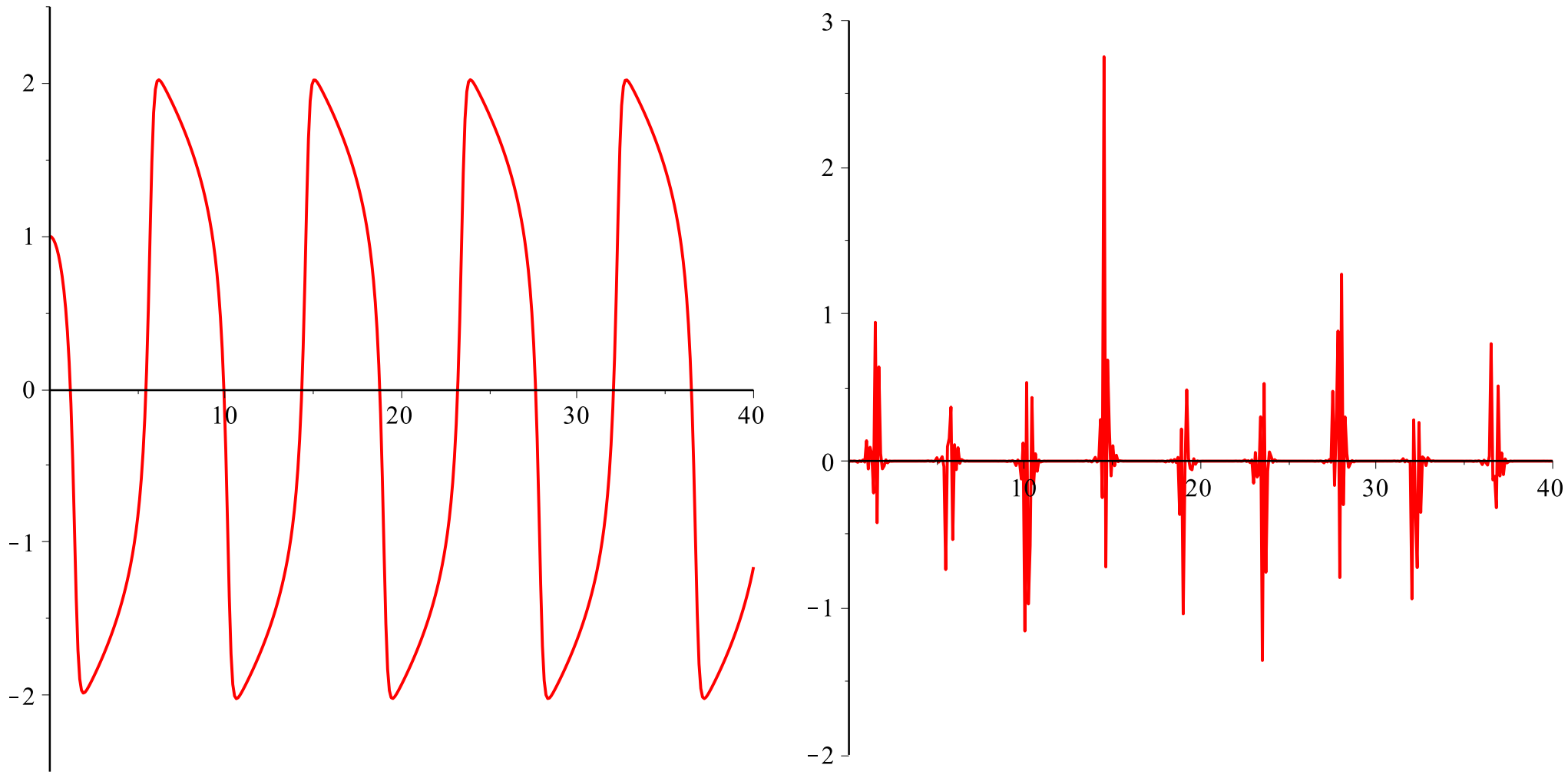

(a) Collocation approximation, $f_{\boldsymbol{\xi}}(x)$ (b) $err(x)$

Figure 3.5: Numerical results for $\mu = 3.0$, $l = 160$, after 1138 iterations for $f_{\boldsymbol{\xi},r}(b)$ to converge to a tolerance of $10^{-4}$.

$f_{\boldsymbol{\xi}}$ for the test cases $\mu = \{20, 40\}$ was so large, that we are unable to present any results. From the results that are available, however, it is possible to recognise the associated pattern of diminished accuracy for a given number of intervals as $\mu$ increases.

Additionally, it is now very apparent that there is a distinction between the $err(x)$ for the part of the solution which exhibits stiff compared against non-stiff behaviour. For the rapid changes in $f_{\boldsymbol{\xi}}(x)$ and, therefore, large absolute derivative values, the uniform mesh size is clearly suboptimal, and there is correspondingly large $err(x)$.

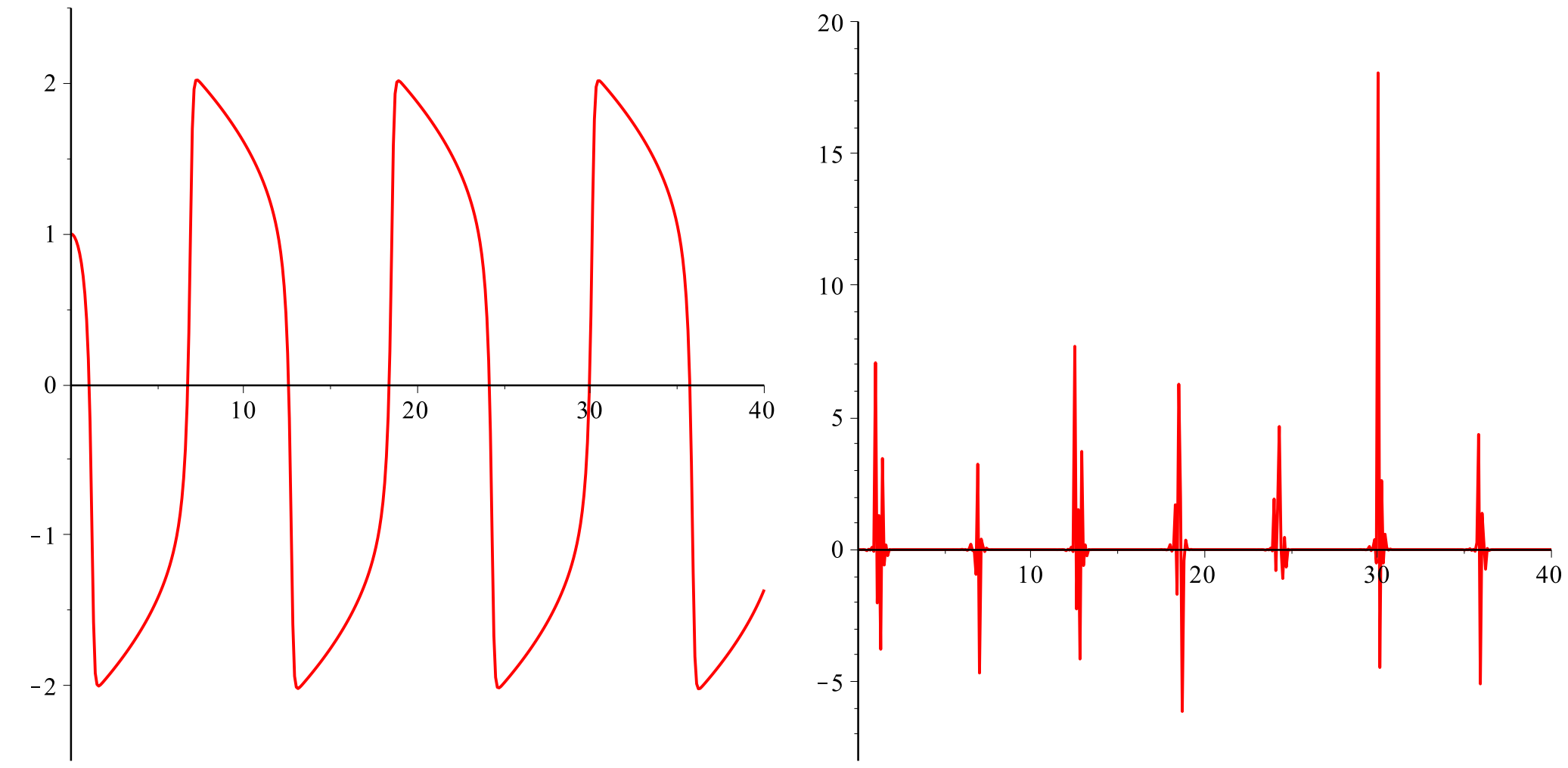

(a) Collocation approximation, $f_{\boldsymbol{\xi}}(x)$ (b) $err(x)$

Figure 3.6: Numerical results for $\mu = 5.0$, $l = 160$, after 1969 iterations for $f_{\boldsymbol{\xi},r}(b)$ to converge to a tolerance of $10^{-4}$.

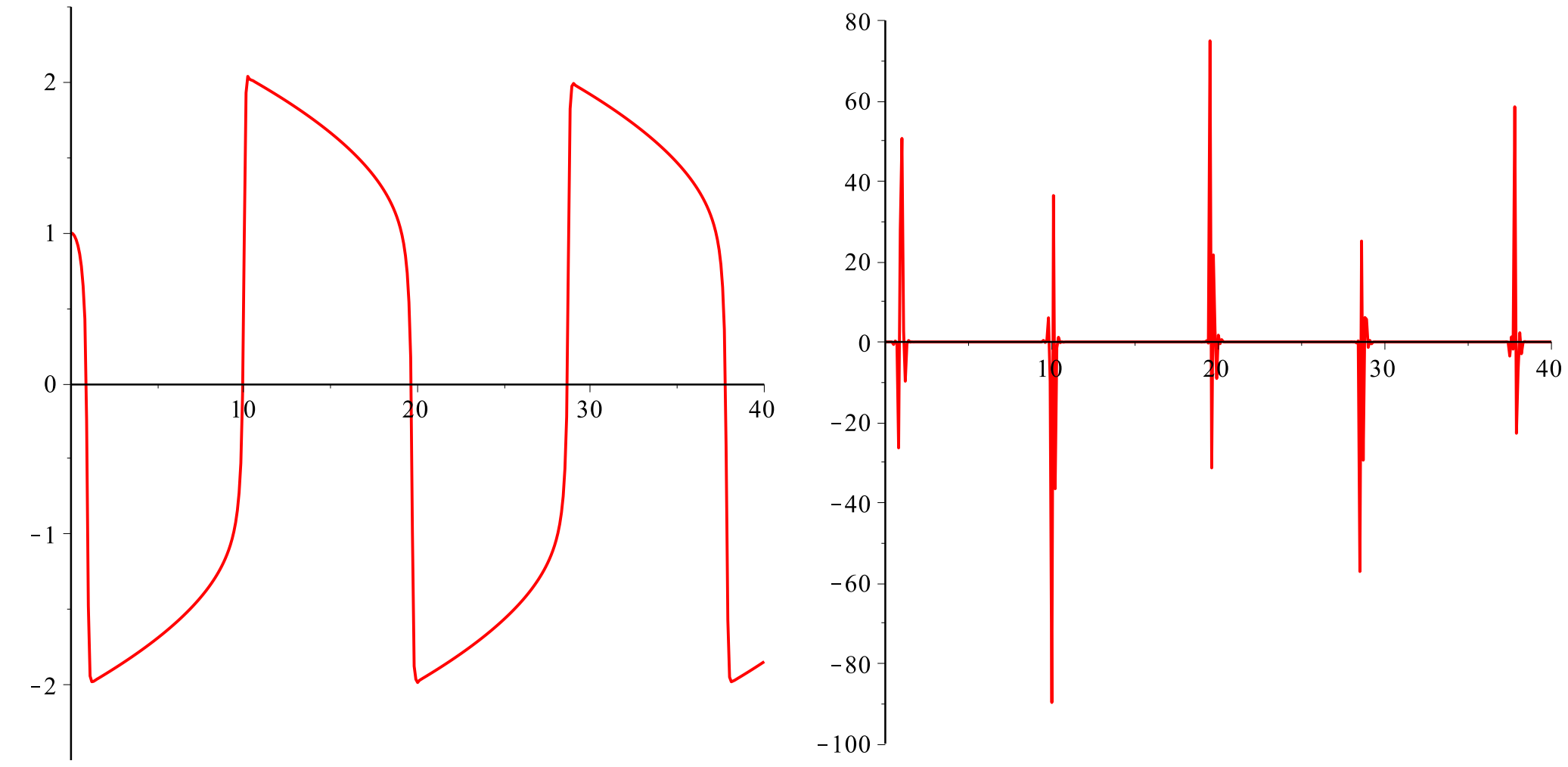

(a) Collocation approximation, $f_{\boldsymbol{\xi}}(x)$ (b) $err(x)$

Figure 3.7: Numerical results for $\mu = 10.0$, $l = 160$, after 5606 iterations for $f_{\boldsymbol{\xi},r}(b)$ to converge to a tolerance of $10^{-4}$.

# Method adaptations

## 4.1 Adaptation: piecewise collocation approximation methods

In this section we bring together the preliminary numerical results and relevant discussed theory to propose an adaptation to the collocation method as it stands. The previous sections have served to illustrate two major problems with the collocation method, even though it has already generated some viable approximations;

(i) where approximations are viable, $cost(f_{\boldsymbol{\xi}})$ becomes prohibitive as $l$, $k$ or $\mu$ is increased further, even for values one might deem reasonably small.

(ii) for $\mu = \{20, 40\}$ we are yet to produce any collocation approximations because $cost(f_{\boldsymbol{\xi}})$ is prohibitive from the outset.

Analogous, perhaps, to the design of pps, we propose that for initial value problems the collocation method benefits from construction in a piecewise manner, as such we term an adaptation of this form a piecewise collocation method.

### 4.1.1 Expanding range collocation

From section 2.1.1 we have seen the approximation is necessarily constructed in a left-to-right manner. This means that for early stages of the iteration process the rightmost part of the solution range simply adds to the cost function, without the approximation benefitting from its inclusion. This first adaptation to the original collocation method we suggest is defined as follows;
Given a solution range $[0..b]$ and uniform $\boldsymbol{\xi}$ choose a second parameter, $w$, where $w$ divides $l$, which defines the number of piecewise collocation approximations to find. So that:

$$\boldsymbol{\xi}_i = \{\xi_1, \xi_2, \dots, \xi_{1+i\frac{l}{w}}\}, \quad i = 1, \dots, w.$$

Then obtain the collocation approximation, $f_{\boldsymbol{\xi}_i}(x)$ for $x \in \left[0..\xi_{1+i\frac{l}{w}}\right]$, to the specified tolerance for $f_{\boldsymbol{\xi}_i}(\xi_{1+i\frac{l}{w}})$, which will necessarily have, $n = (k-2)i\frac{l}{w} + 2$, B-spline coefficients. Once found, seek the next collocation approximation $f_{\boldsymbol{\xi}_{i+1}}(x)$ with an initial supposed solution where the first $\left((k-2)i\frac{l}{w} + 2\right)$ B-spline coefficients are carried over, and the remaining coefficients are set to, for example, zero. The B-splines close to $\xi_{1+i\frac{l}{w}}$ do change in the course of the range extension so their coefficients are no longer equivalent but the iteration process will continue to converge in any case and a minor amount of misaligned B-splines if of no overall concern.

The final approximation is then:

$$f_{\boldsymbol{\xi}}(x) = f_{\boldsymbol{\xi}_w}(x)\ .$$

As an assessment of the likely cost saving, we first make the assumption that the overall number of iterations, $N$, for a piecewise collocation approximation to be found remains the same as when applying the original collocation method. There is little reason at this stage, other than the insertion of a new initial supposed solution for each $i$, to suggest it differs greatly. Secondly we purport an averaging assumption that the number of iterations to converge each $f_{\boldsymbol{\xi}_i}$ are broadly the same. This is unlikely to be true given the discrepancy between stiff and non-stiff solution ranges, but taken over the entirety of the solution range it may be reasonable. Then we have that:

$$\begin{aligned} cost(f_{\boldsymbol{\xi}}) &= \sum_{i=1}^{w} cost(f_{\boldsymbol{\xi}_i}) = \sum_{i=1}^{w} O\left((k^3+k^2+k)\frac{N}{w} i \frac{l}{w}\right) \\ &= O\left((k^3+k^2+k)Nl\frac{w(w+1)}{2w^2}\right) \end{aligned} \tag{4.1}$$

So this adaptation could potentially halve the cost if there exists the scope for $w$ to be chosen high enough.

### 4.1.2 Segmented collocation

This second adaptation makes use of the fact that once the leftmost solution has converged then it also adds to the cost function without the approximation benefitting from its inclusion. We define this method as follows;
With the definition of $\boldsymbol{\xi}_i$ above let $f_{\boldsymbol{\xi}_i}(x)$ be the collocation approximation for $x \in \left[\xi_{1+(i-1)\frac{l}{w}} .. \xi_{1+i\frac{l}{w}}\right]$, where the initial values required to determine $f_{\boldsymbol{\xi}_i}(x)$ are taken as $f_{\boldsymbol{\xi}_{i-1}}(\xi_{1+(i-1)\frac{l}{w}})$ and $Df_{\boldsymbol{\xi}_{i-1}}(\xi_{1+(i-1)\frac{l}{w}})$. The supposed initial solution for each segment is determined, for example, akin to section 2.1.4.

The final approximation is then:

$$f_{\boldsymbol{\xi}}(x) = \left\{ f_{\boldsymbol{\xi}_i}(x) : x \in \left[\xi_{1+(i-1)\frac{l}{w}} .. \xi_{1+i\frac{l}{w}}\right]\right\}$$

As an assessment of cost, again we make the same assumptions as before to obtain:

$$\begin{aligned} cost(f_{\boldsymbol{\xi}}) &= \sum_{i=1}^{w} cost(f_{\boldsymbol{\xi}_i}) = \sum_{i=1}^{w} O\left((k^3+k^2+k)\frac{N}{w}\frac{l}{w}\right) \\ &= O\left((k^3+k^2+k)N\frac{l}{w}\right) \end{aligned} \tag{4.2}$$

This suggests if $l = w$ we can potentially achieve a cost which is dependent only on the number of iterations that the original method takes to converge, which is a considerable improvement in the case for large $\mu$ where correspondingly large $l$ are required.

# Peremptory results

## 5.1 Segmented collocation method applied to large $\mu$

### 5.1.1 Revisited collocation approximations

Figures 5.1 and 5.2 show the results that we were unable to produce using the original collocation method in section 3.3, by applying the segmented collocation method for $\mu = \{20, 40\}$. As a further demonstration of the improvement, figure 5.3 presents a comparison in the case where $\mu = 10$ to figure 3.7 where we note the that $\|err(x)\|_\infty$ is an $O(10^{-1})$ lower and its cost, of 6.1 CPU seconds, is of $O(10^{-2})$ that of the previous method.

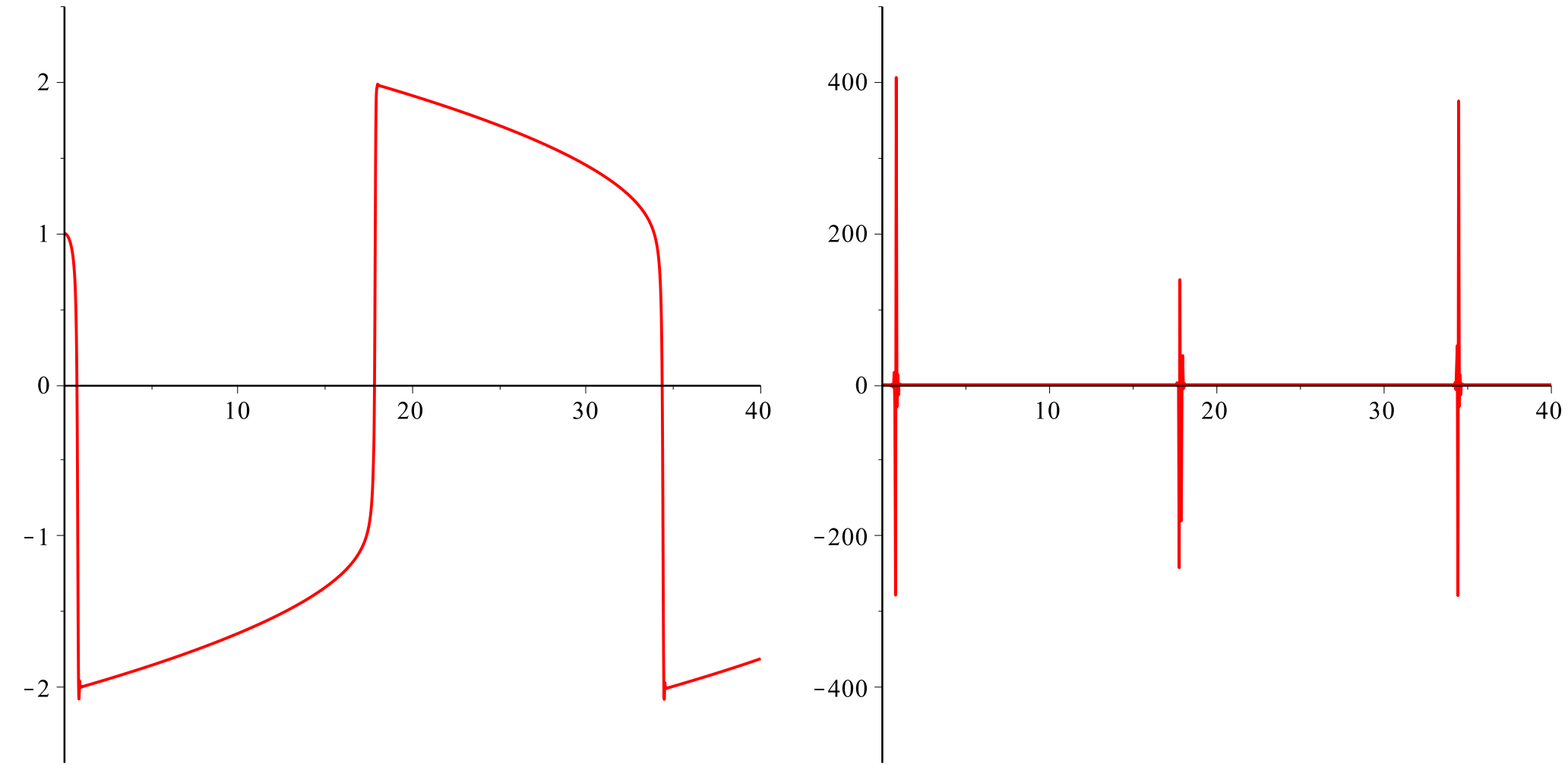

(a) Segmented collocation approximation, $f_{\boldsymbol{\xi}}(x)$ (b) $err(x)$

Figure 5.1: Numerical results on [0..40] for $\mu = 20.0$, $k = 5$, $l = 320$, $w = 40$ after 199 iterations for $f_{\boldsymbol{\xi}_i,r}(\xi_{1+(i-1)\frac{l}{w}})$ to converge to a tolerance of $10^{-4}$.

#### Approximation phenomena

Figure 5.4 highlights a phenomena observed in the results from figures 5.1 and 5.2.

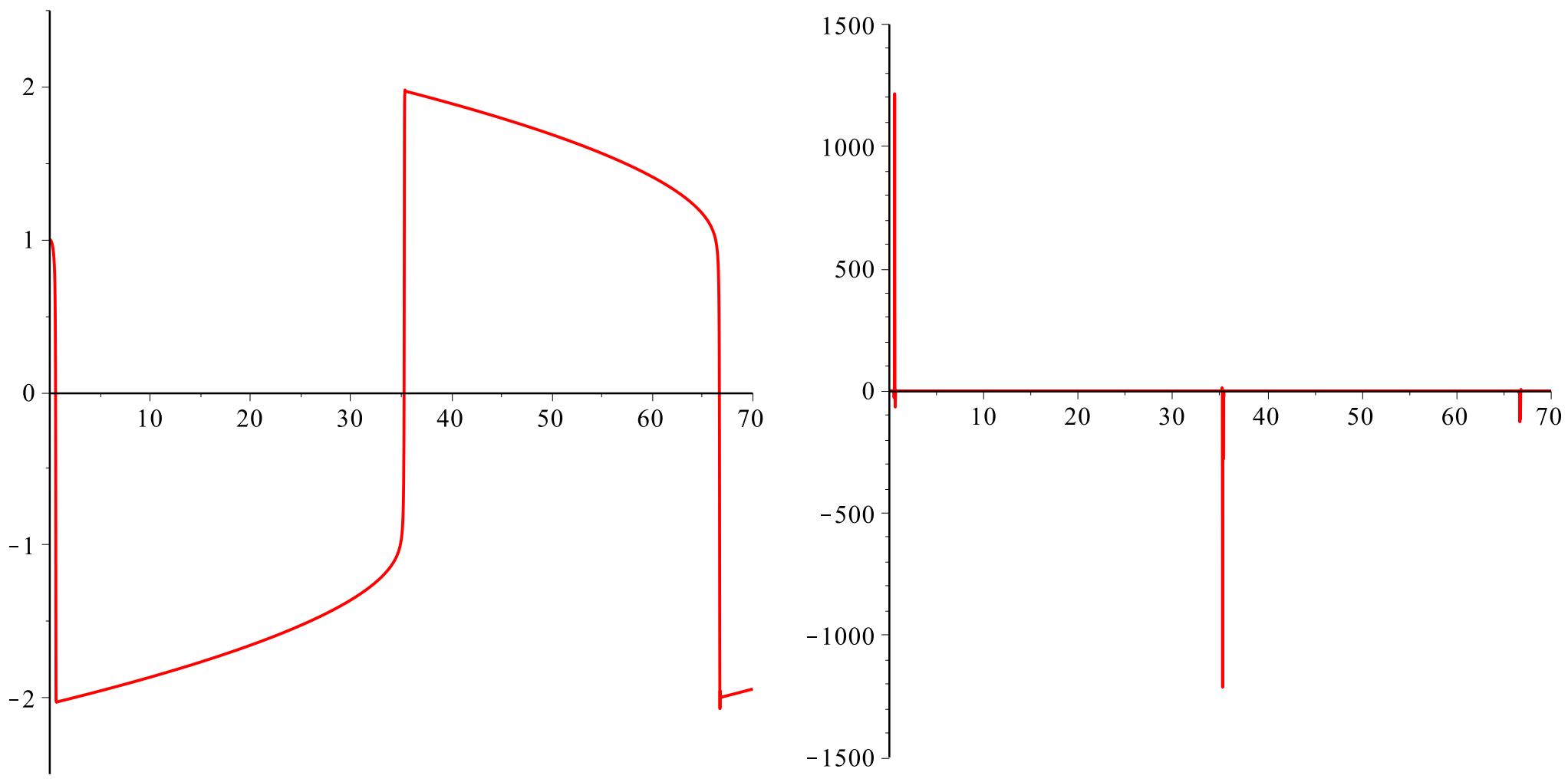

(a) Segmented collocation approximation, $f_{\boldsymbol{\xi}}(x)$ (b) $err(x)$

Figure 5.2: Numerical results on [0..70] for $\mu = 40.0$, $k = 5$, $l = 1120$, $w = 70$ after 815 iterations for $f_{\boldsymbol{\xi}_i,r}(\xi_{1+(i-1)\frac{l}{w}})$ to converge to a tolerance of $10^{-4}$.

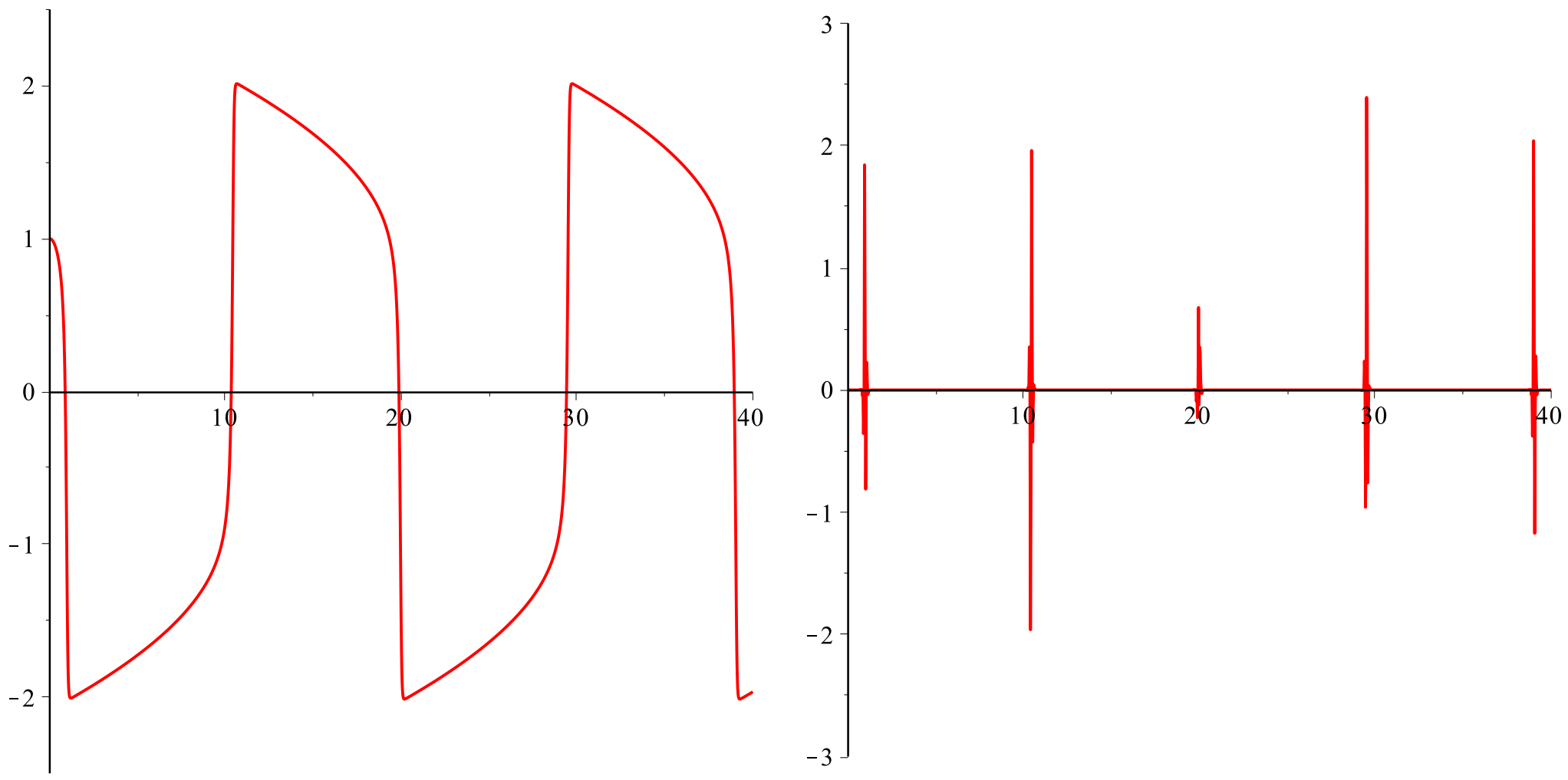

(a) Segmented collocation approximation, $f_{\boldsymbol{\xi}}(x)$ (b) $err(x)$

Figure 5.3: Numerical results for $\mu = 10.0$, $k = 6$, $l = 640$, $w = 40$ after 207 iterations for $f_{\boldsymbol{\xi}_i,r}(\xi_{1+(i-1)\frac{l}{w}})$ to converge to a tolerance of $10^{-4}$.

### 5.1.2 Numerical data collected on segmented collocation method

Table 5.1 documents a collection of results for the study $\mu = 40$, with figure 5.5 highlighting the behaviour of each of the eight generated approximations about a key stiff behavioural point, $x \approx 67$.

Table 5.2 documents a collection of results ascertained through the application of the segmented collocation method for various parameters.

Table 5.3 documents specific choices of $k, l$ to ascertain approximations of a similar degree of accuracy for each $\mu$. Specific $w$ are also shown to highlight the concept of an optimal choice.

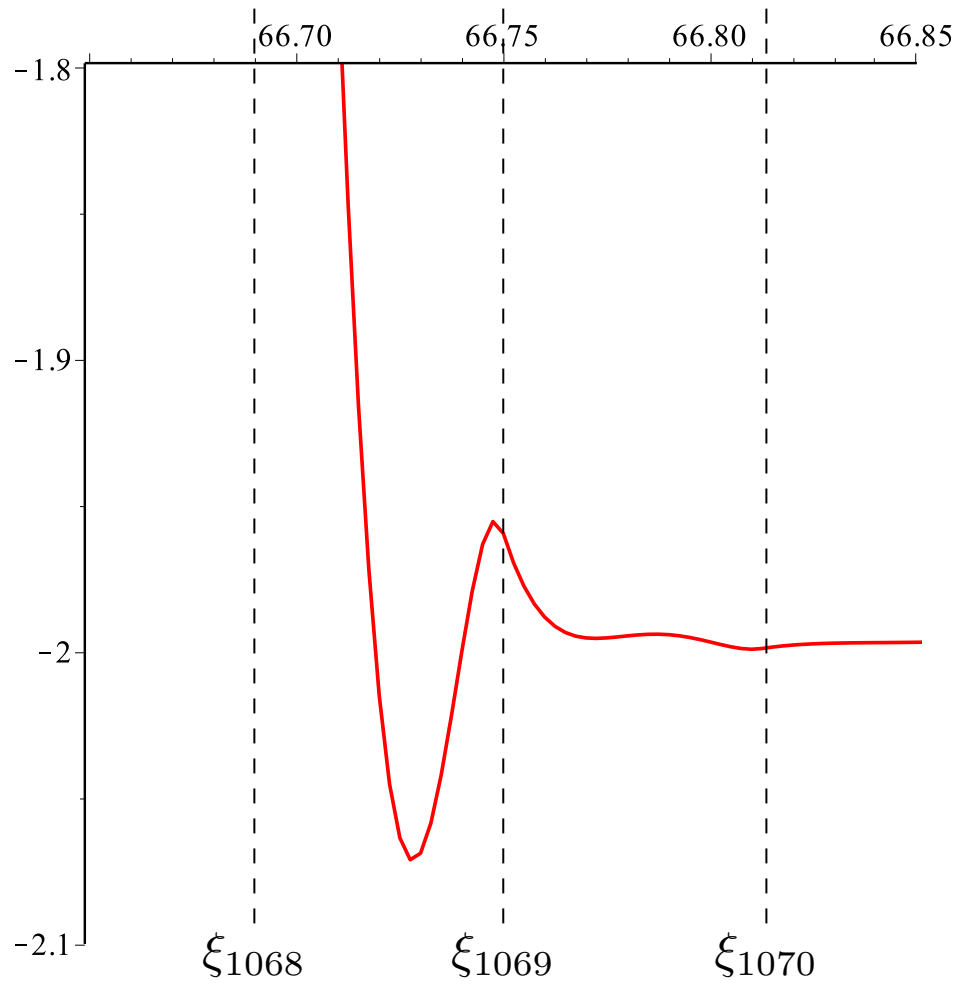

Figure 5.4: Highlighted phenomena in the presented segmented collocation approximation for $\mu = 40.0$.

| Method | $\mu$ | $l$ | $w$ | $k$ | $N$ | $cost(f_{\boldsymbol{\xi}})$ | $\|err(x)\|_\infty$ |
|---|---|---|---|---|---|---|---|
| Segmented | 40.0 | 700 | 140 | 6 | 815 | 14.9 s | 3.2E3 |
| Segmented | 40.0 | 875 | 175 | 6 | 673 | 13.9 s | 3.9E3 |
| Segmented | 40.0 | 1000 | 200 | 6 | 630 | 14.7 s | 2.0E3 |
| Segmented | 40.0 | 1400 | 350 | 6 | 1078 | 16.1 s | 1.0E3 |
| Segmented | 40.0 | 1750 | 350 | 6 | 1077 | 17.1 s | 6.0E2 |
| Segmented | 40.0 | 2000 | 400 | 6 | 1228 | 18.6 s | 4.1E2 |
| Segmented | 40.0 | 3500 | 700 | 6 | 2127 | 21.1 s | 8.0E1 |
| Segmented | 40.0 | 7000 | 1400 | 6 | 4223 | 26.9 s | 5.0E0 |

Table 5.1: $N$, $cost(f_{\boldsymbol{\xi}})$, in terms of CPU seconds, and $\|err(x)\|_\infty$ for the study $\mu = 40$ using different parameter values $l, w$ but constant $k$.

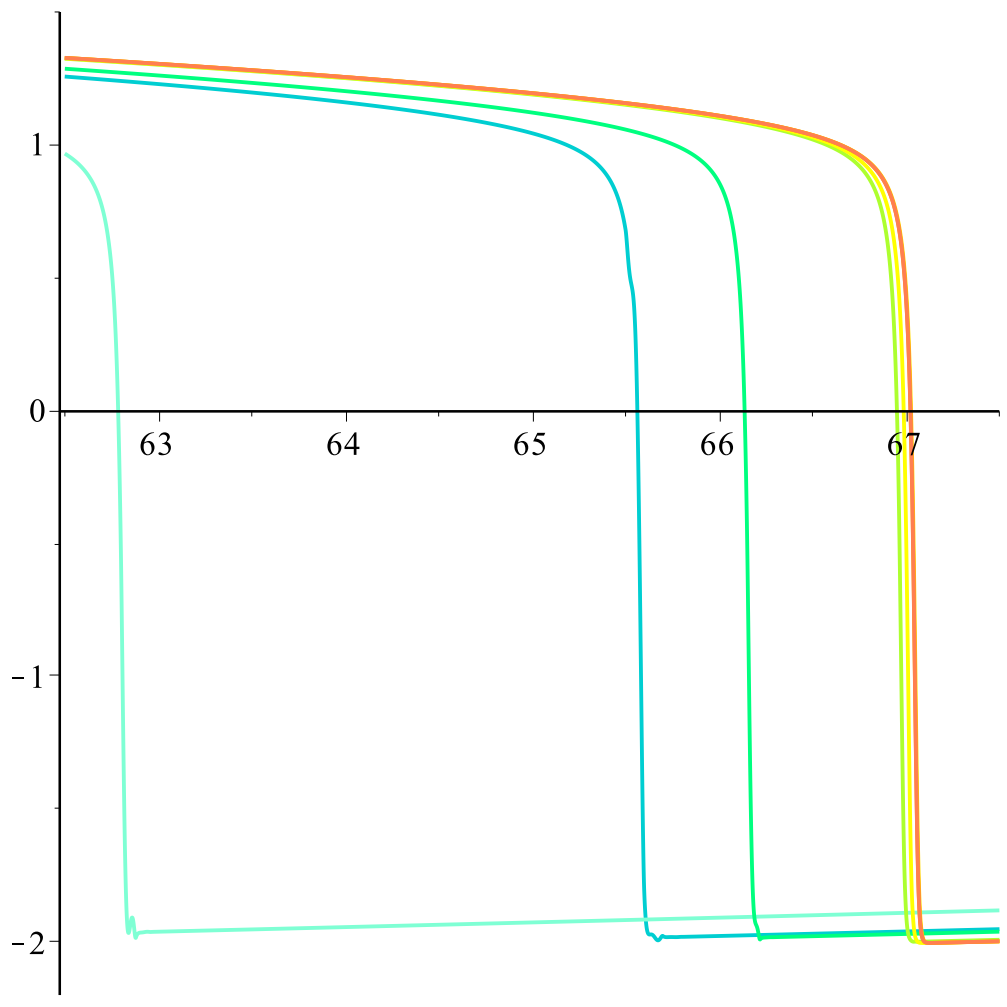

(a) Range displayed [62.5 .. 67.5]

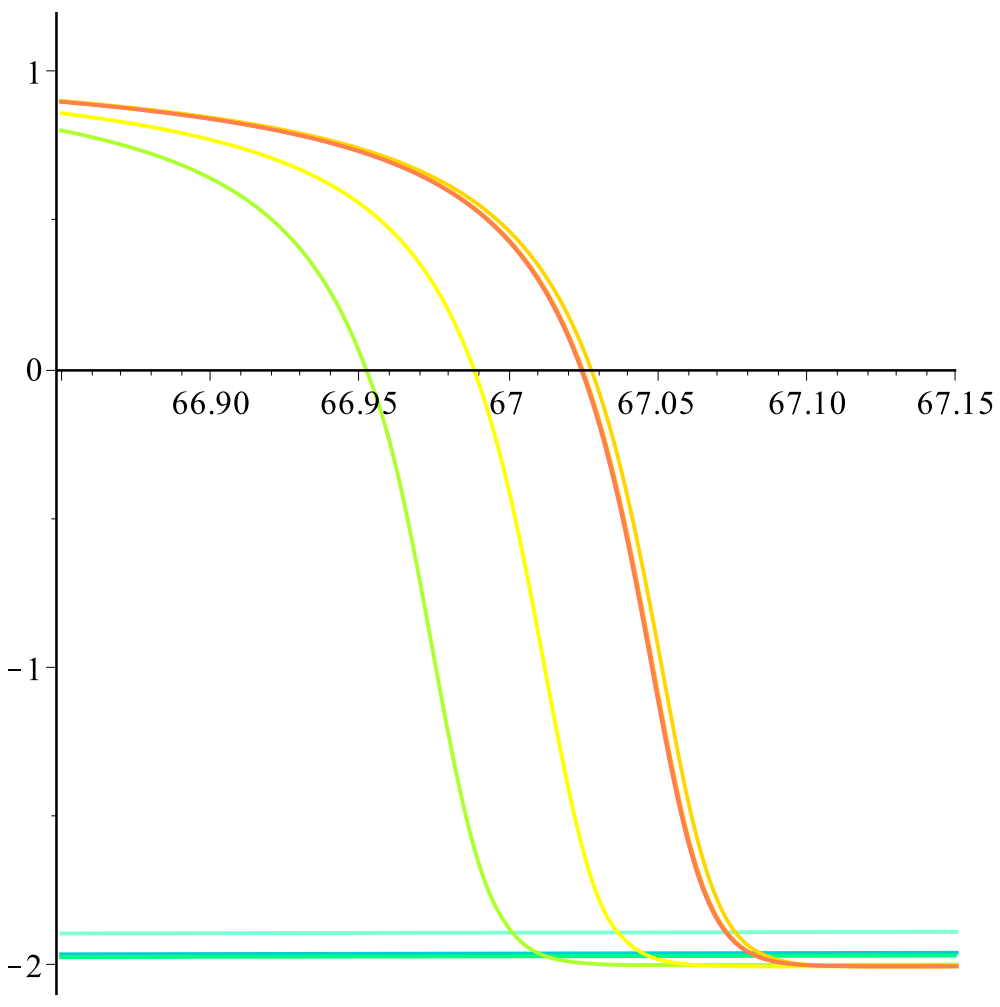

(b) Range displayed [66.85 .. 67.15]

Figure 5.5: Displaying the eight collocation approximations generated in table 5.1 for specific range of $x$, colour wise for $l = \{700, \ldots, 7000\}$ equivalent to {blue,..,green,..,yellow,..,red}, highlighting translation error.

| Method | $\mu$ | $l$ | $w$ | $k$ | $N$ | $cost(f_{\boldsymbol{\xi}})$ | $\|err(x)\|_\infty$ |
|---|---|---|---|---|---|---|---|
| Original | 10.0 | 160 | 1 | 5 | 5606 | 695 s | 137 |
| Segmented | 10.0 | 160 | 2 | 5 | 4432 | 268 s | 137 |
| Segmented | 10.0 | 160 | 4 | 5 | 2562 | 75.0 s | 137 |
| Segmented | 10.0 | 160 | 8 | 5 | 1258 | 18.4 s | 137 |
| Segmented | 10.0 | 160 | 16 | 5 | 369 | 3.71 s | 137 |
| Segmented | 10.0 | 160 | 20 | 5 | 216 | 1.95 s | 137 |
| Segmented | 10.0 | 160 | 40 | 5 | 194 | 1.88 s | 137 |
| Segmented | 10.0 | 160 | 80 | 5 | 292 | 2.57 s | 137 |
| Segmented | 10.0 | 160 | 160 | 5 | 524 | 3.94 s | 137 |
| Original | 3.0 | 160 | 1 | 5 | 1138 | 162 s | 2.73 |
| Segmented | 3.0 | 160 | 2 | 5 | 697 | 41.8 s | 2.73 |
| Segmented | 3.0 | 160 | 4 | 5 | 279 | 8.03 s | 2.73 |
| Segmented | 3.0 | 160 | 8 | 5 | 146 | 2.21 s | 2.73 |
| Segmented | 3.0 | 160 | 16 | 5 | 122 | 1.07 s | 2.73 |
| Segmented | 3.0 | 160 | 20 | 5 | 135 | 1.01 s | 2.73 |
| Segmented | 3.0 | 160 | 40 | 5 | 190 | 0.85 s | 2.73 |
| Segmented | 3.0 | 160 | 80 | 5 | 304 | 0.87 s | 2.73 |
| Segmented | 3.0 | 160 | 160 | 5 | 539 | 1.05 s | 2.73 |
| Original | 3.0 | 80 | 1 | 5 | 1202 | 76.8 s | 12.2 |
| Segmented | 3.0 | 80 | 2 | 5 | 859 | 23.7 s | 12.2 |
| Segmented | 3.0 | 80 | 4 | 5 | 473 | 6.66 s | 12.2 |
| Segmented | 3.0 | 80 | 8 | 5 | 221 | 1.65 s | 12.2 |
| Segmented | 3.0 | 80 | 16 | 5 | 121 | 0.61 s | 12.2 |
| Segmented | 3.0 | 80 | 20 | 5 | 136 | 0.62 s | 12.2 |
| Segmented | 3.0 | 80 | 40 | 5 | 187 | 0.61 s | 12.2 |
| Segmented | 3.0 | 80 | 80 | 5 | 304 | 2.10 s | 12.2 |

Table 5.2: $N$, $cost(f_{\boldsymbol{\xi}})$, in terms of CPU seconds, and $\|err(x)\|_\infty$ for various parameter values $\mu, l, w$, using different collocation methods for constant $k$. Highlighting; the relationship between $w$ and $N$, the relationship between $w$ and $cost(f_{\boldsymbol{\xi}})$, and the consistent $\|err(x)\|_\infty$ for variable $w$.

| Method | $\mu$ | $l$ | $w$ | $k$ | $N$ | $cost(f_{\boldsymbol{\xi}})$ | $\|err(x)\|_\infty$ |
|---|---|---|---|---|---|---|---|
| Segmented | 0.5 | 40 | 2 | 4 | 48 | 0.34 s | 1.10 |
| Segmented | 1.0 | 80 | 4 | 4 | 68 | 0.64 s | 1.19 |
| Segmented | 2.0 | 120 | 12 | 5 | 61 | 0.69 s | 1.31 |
| Segmented | 3.0 | 210 | 21 | 5 | 138 | 1.35 s | 1.37 |
| Segmented | 3.0 | 210 | 42 | 5 | 193 | 1.04 s | 1.37 |
| Segmented | 5.0 | 250 | 25 | 6 | 154 | 2.57 s | 1.91 |
| Segmented | 5.0 | 250 | 50 | 6 | 218 | 2.05 s | 1.91 |
| Segmented | 5.0 | 250 | 125 | 6 | 425 | 2.20 s | 1.91 |
| Segmented | 10.0 | 400 | 10 | 8 | 777 | 128 s | 1.18 |
| Segmented | 10.0 | 400 | 16 | 8 | 253 | 19.5 s | 1.18 |
| Segmented | 10.0 | 400 | 20 | 8 | 174 | 11.3 s | 1.18 |
| Segmented | 10.0 | 400 | 25 | 8 | 174 | 9.56 s | 1.18 |
| Segmented | 10.0 | 400 | 40 | 8 | 207 | 7.45 s | 1.18 |
| Segmented | 10.0 | 400 | 50 | 8 | 224 | 6.79 s | 1.18 |
| Segmented | 10.0 | 400 | 80 | 8 | 284 | 5.98 s | 1.18 |
| Segmented | 10.0 | 400 | 100 | 8 | 344 | 6.11 s | 1.18 |
| Segmented | 10.0 | 400 | 200 | 8 | 636 | 6.43 s | 1.18 |
| Segmented | 20.0 | 800 | 16 | 9 | 608 | 160 s | 1.47 |
| Segmented | 20.0 | 800 | 20 | 9 | 162 | 38.5 s | 1.47 |
| Segmented | 20.0 | 800 | 25 | 9 | 201 | 30.6 s | 1.47 |
| Segmented | 20.0 | 800 | 40 | 9 | 197 | 20.3 s | 1.47 |
| Segmented | 20.0 | 800 | 50 | 9 | 240 | 19.9 s | 1.47 |
| Segmented | 20.0 | 800 | 100 | 9 | 333 | 15.5 s | 1.47 |
| Segmented | 20.0 | 800 | 160 | 9 | 507 | 15.2 s | 1.47 |
| Segmented | 20.0 | 800 | 200 | 9 | 624 | 15.2 s | 1.47 |
| Segmented | 20.0 | 800 | 400 | 9 | 1223 | 15.6 s | 1.47 |
| Segmented | 40.0 | 1600 | 40 | 10 | 197 | 66.4 s | 1.19 |
| Segmented | 40.0 | 1600 | 50 | 10 | 254 | 62.3 s | 1.19 |
| Segmented | 40.0 | 1600 | 80 | 10 | 287 | 44.6 s | 1.19 |
| Segmented | 40.0 | 1600 | 100 | 10 | 328 | 41.5 s | 1.19 |
| Segmented | 40.0 | 1600 | 160 | 10 | 507 | 40.7 s | 1.19 |
| Segmented | 40.0 | 1600 | 200 | 10 | 621 | 40.1 s | 1.19 |
| Segmented | 40.0 | 1600 | 320 | 10 | 979 | 40.2 s | 1.19 |
| Segmented | 40.0 | 1600 | 800 | 10 | 2415 | 43.3 s | 1.19 |
| Segmented | 80.0 | 3125 | 625 | 11 | 1886 | 109 s | 1.63 |

Table 5.3: $N$, $cost(f_{\boldsymbol{\xi}})$, in terms of CPU seconds, and $\|err(x)\|_\infty$ for various parameter values $\mu, l, w, k$. Highlighting; the relationship between $w$ and $N$, the relationship between $w$ and $cost(f_{\boldsymbol{\xi}})$, and the achievement of consistent accuracy in terms of $\|err(x)\|_\infty$ for variable $\mu$ with careful choice of $k, l$.

# Conclusions and proposed further study

## 6.1 Conclusions

The original collocation method outlined in sections 1.2.9 and 1.2.10, used in conjunction with an initial supposed solution proposed in section 2.1.4, proves incapable of reliably producing approximations for the Van der Pol equation 1.15 for large $\mu$. This is due to the associated cost of generating the approximations which increases with $\mu$ due to;

(i) an expressed dependence on $\mu$ of the total number of iterations, $N$, required for the approximation to converge to a specified tolerance,

(ii) the requirement of an increased number of intervals, $l$, for a given solution range, $[0..b]$, to obtain sufficient accuracy about exhibited stiff behavioural points, to produce stable and permissible approximations,

(iii) the requirement of an increase to the order of pp, $k$, to obtain improved accuracy, again specifically about exhibited stiff behavioural points, to allow stable and permissible approximations.

In some cases even for small $\mu$ the costs associated with increasing $k, l$ can become prohibitive. As a result, adaptations are proposed to the original collocation method, seeking to benefit, in cost terms, from the structural dependence of the method, outlined in section 2.1.1. One such adaptation, which is termed here the ”segmented collocation method”, is developed in section 4.1.2 and shown numerically to produce approximations of equivalent error to the original collocation method for a fraction of the cost. This permits a greater scope of parameters to be feasibly tested and examined in this study. The resultant approximations generated in this way are reliably attainable and flexible in terms of their parametrisation for a given level of accuracy, as demonstrated in the figures and tables of chapter 5.

### 6.1.1 Approximation error

**Phenomena**

Figure 5.4 shows a type of phenomena which materialises in the case of large $\mu$. It bears similarity to that of Gibbs phenomenon in which the Fourier series of a piecewise continuously differentiable periodic function behaves at a jump discontinuity, or that of signal overshoot when a signal or function exceeds its target. We can only assert that for larger $\mu$ this phenomena is more difficult to overcome, and can yield unstable or unreliable approximations. This phenomena shows form because the real solution typically changes too rapidly for the pp approximation to represent it accurately.

**Translation error**

Section 2.1.1 already outlines the theoretical expectation that the approximation can become translated in $x$ and any subsequent evaluation of an accurate approximation to the differential equation suffers from this kind of error. Specifically, about stiff behavioural points the $|err(x)|$ can become very large, for insufficiently chosen $k, l$, as in figure 5.2(b), and equivalently this can result in greater overall approximation error about those points, which is frequently carried forward as translation error, see figure 5.5.

To counteract both of the above, we require larger $k, l$, or a better designed break sequence $\boldsymbol{\xi}$. Again, figure 5.5 shows the eradication of phenomena for larger $l$ and the convergence of the approximations toward one with specific $x$-axis intercept points. We assert a relationship between $\|err(x)\|_\infty$ and $\|g - f_{\boldsymbol{\xi}}\|_\infty$, where for the case of collocation with boundary conditions we state an explicit relationship in equation 2.1. Convergence in one must imply convergence in the other, and therefore conclude that convergence is in fact toward the real solution, for sufficiently increased $k, l$.

### 6.1.2 Cost of approximations generated with segmented and original collocation methods

Equations 2.2 and 4.2 provide a theoretical basis for cost analysis. They are, however, based on an optimally efficient algorithm which we commented in section 3.1 is not quite achieved in the numerical method used in this study, although it is broadly similar. Nonetheless, we restate our previously derived theoretical cost functions:

$$\begin{aligned} cost(f_{\boldsymbol{\xi}}^{ori.}) &= O\left((k^3+k^2+k)N^{ori.}l\right), \\ cost(f_{\boldsymbol{\xi}}^{seg.}) &= O\left((k^3+k^2+k)N^{seg.}\frac{l}{w}\right), \end{aligned} \tag{6.1}$$

where the superscripts indicate a collocation solution by the original or segmented method.

In the derivation, in section 4.1.2, of the second cost function we made the assumption that $N^{seg.} = N^{ori.}$ without any information to the contrary. However, from table 5.2, it is now possible to derive the following relationship:

$$N^{seg.} = \frac{N^{ori.}}{w^\lambda} + N^{min}w, \tag{6.2}$$

where $N^{min}$ represents the minimum number of iterations it generally takes to converge to a specified tolerance for a single small interval.

We interpret this result in the following way; for large $w$ each segment represents a single interval which typically takes $N^{min}$ iterations to converge and therefore registers a tally of $N^{min}w$ total iterations, and for the $\frac{N^{ori.}}{w^\lambda}$ term we conjecture that it is related to the original total number of iterations divided by the "resetting" of the initial supposed solution for each $w$, which in turn were designed at the outset with the purpose of assisting convergence.

If $N^{ori.}$ and $N^{min}$ are known then minimum $N^{seg.}$ is obtained if:

$$w := \left(\frac{\lambda N^{ori.}}{N^{min}}\right)^{\frac{1}{1+\lambda}}, \tag{6.3}$$

and substituting this into equation 6.1 produces the significant and greatly reduced cost function, independent of $N^{ori.}$:

$$cost(f_{\boldsymbol{\xi}}^{seg.}) = O\left((k^3+k^2+k)\left(1+\frac{1}{\lambda}\right)N^{min}l\right), \quad \text{for } w = \left(\frac{\lambda N^{ori.}}{N^{min}}\right)^{\frac{1}{1+\lambda}}. \tag{6.4}$$

Section 7.5 outlines the determination of $\lambda = 1.0$ from numerical results and asserts that $N^{min}$ typically averages just above 3.0.

However, by directly substituting 6.2 into 6.1 it becomes clear that as $w \to \infty$ we achieve slightly better, indicating smaller intervals save more costs than those obtained by specifying a $w$ to minimise the total number of iterations:

$$cost(f_{\boldsymbol{\xi}}^{seg.}) = O\left((k^3 + k^2 + k)N^{min}l\right), \quad \text{for } w \to \infty. \tag{6.5}$$

Tables 5.2 and 5.3 serve to verify this theoretical analysis. They show that the minimum number of CPU seconds to generate the collocation approximation is typically achieved for $w$ greater than that which results in minimum $N^{seg.}$ for specified parameters, $\mu, k, l$. More than this, the tables indicate that minimum $cost(f_{\boldsymbol{\xi}})$ is about half that of the $cost(f_{\boldsymbol{\xi}})$ when minimum $N^{seg.}$ is achieved predicted by equations 6.4 and 6.5. At odds with the theoretical analysis is the slight increase in $cost(f_{\boldsymbol{\xi}})$ measured in CPU seconds when $w$ tends to $l$. This remains an item of investigation but it is suspected that in section 7.4 where the approximate cost function is derived there are some terms ignored, relevant to the initialisation of the problem and calculation of B-spline values and derivative values, that become more relevant as $w$ tends to $l$, because these processes are performed more frequently.

### 6.1.3 Numerical limits and measure of equation stiffness

The general theme of this study has shown that as $\mu$ increases then so too does the stiffness of the equation, that is the character of the solution about specific $x$-values rapidly changes resulting in very high first and second derivative values. With regard to approximations this can mean more instability, increased time to converge (when approximations do converge), and increased inaccuracy.

There are a number of ways one might choose to characterise stiffness and consider the limits of a numerical process. This study has opted to consider the costs of generating approximations. Data from table 5.3, analysed in section 7.5, indicates an estimated relationship that for an approximation, $f_{\boldsymbol{\xi}}^{seg.}$, on [0..40] to have $\|err(x)\|_\infty \leq 2.0$, for $\mu \geq 2$:

$$cost(f_{\boldsymbol{\xi}}^{seg.}) \approx \frac{1}{4}\mu^{\frac{7}{5}} \text{ CPU seconds.}$$

We conclude that utilising the segmented collocation method presented in this study, one should expect approximations to be permissible, for appropriate accuracy and in reasonable time, for $\mu$ in the low hundreds. Of course this is subject to processing speed and one's level of patience, and the expectation of extrapolated relationships and no new phenomena appearing.

## 6.2 Further study

### 6.2.1 $N^{ori.}$ and $N^{min}$ as functions of $k$

Given the dependence of $cost(f_{\boldsymbol{\xi}}^{ori.})$ on $N^{ori.}$ it did seem that an analysis of $N^{ori.}$ as a function of $k$ should be potentially be included. It was observed in passing that a higher $k$ typically resulted in lower $N^{ori.}$ but not always and so this was an unreliable conclusion. The developed theory on segmented collocation ultimately made this rather moot, because that more efficient method is dependent on $N^{min}$ which is typically small enough to not exhibit dependence to $k$, although the further development of this method, for example in the areas proposed below, might open up such an item for exploration.

### 6.2.2 Chaos

An interesting result observed in passing was that the value, $N^{ori.}$, was dependent upon different initial supposed solutions, $f_{\boldsymbol{\xi},0}$. The path of the converging approximation seemed to demonstrate no pattern associated with a change in the initial supposed solution, indicative of a chaotic process, which may prove to be interesting to explore in addition.

### 6.2.3 Non-uniform break sequence, $\boldsymbol{\xi}$

The highlighted phenomena discussed in section 6.1.1 serves to demonstrate that the break sequence is an integral part of the applicability of a collocation method. It stands to reason that the development of an optimum choice for $\boldsymbol{\xi}$ would considerably reduce costs and perhaps eliminate the phenomena observed. By considerably reducing the number of intervals a study might also show that the collocation method in its original form can be used to generate approximations for large $\mu$, by effectively reducing the overall cost function.

A study of this nature might well set out with the purpose of tackling the following;

(i) developing theory to establish an optimum break sequence (with respect to minimising $cost(f_{\boldsymbol{\xi}})$, as opposed to minimising error or eliminating phenomena) for a known limit cycle period for given $\mu$,

(ii) developing or researching a numerical or analytical method of establishing the period of the limit cycle and predicting key $x$-values around which the dominant share of intervals might be expected and required in (i).

I envisage that if the above were to be accomplished solely in a numerical sense that the researcher would have to be careful to not lose any gains generated through using optimum breakpoints, by iteratively solving for the period of the limit cycle and for determining the key $x$-values used to calculate such optimum breakpoints.

### 6.2.4 Comparison to other methods

In my personal opinion the most interesting area of further study would be the comparison of using the segmented collocation method to generate approximations versus other known numerical procedures, such as the Runge-Kutta family. An investigation of this nature could take a variety of forms, it could for example;

(i) consider the cost of using different methods to produce approximations of assumed accuracy,

(ii) consider the suitability of different methods for their ability to handle stiff equations, and any arising phenomena,

(iii) as well as (i) compare flexibility of methods when generating approximations, such as utilisation of non-uniform breakpoint sequences in the case of the collocation method, or utilising adaptive Runge-Kutta methods as opposed to standard explicit Runge-Kutta methods,

(iv) consider the numerical limits of each method and test which are more suited toward dealing with super-stiff equations.

(v) consider any similarities in methods in terms of the forms of approximations they yield.

### 6.2.5 Applicability to other equations

The Van der Pol equation as shown in section 1.3 has a unique and stable limit cycle and solution given initial conditions. It is asserted that an area of development and testing could be into equations which do not necessarily have unique solutions or predictable solutions. It may result that spurious behaviour is attained, of a type that is impossible to reproduce when studying equations with similar properties to those of the Van der Pol equation, but we neglect to speculate further here on the sorts of behaviour that may materialise.

### 6.2.6 Preferential outcome

After exploring all areas of further study one may determine;

(i) whether the segmented collocation method is equally applicable to other non-linear, second order, ODEs with initial conditions and with an initial supposed solution in some neighbourhood of the target solution,

(ii) whether there exists a general numerical method, applicable to (i), to generate an optimum or certifiably good non-uniform break sequence, for a comparatively small cost,

(iii) whether the segmented collocation method is competitive, in terms of cost and accuracy, to the arguably more widely used adaptive Runge-Kutta family of methods, again for the problems proposed in (i).

As a potential candidate as a general use, common tool in the efficient and accurate approximation of non-linear, second order, ODEs with given initial conditions, the areas of further study proposed above have been structured in such a way to make possible such an assessment, one way or the other, and earnest determination whether segmented collocation stands to be a worthwhile method in this context.

# Appendix

## 7.1 B-spline coefficient dependence referenced in section 2.1.1

The generic matrix $\mathbf{A}$ has the following form, as described in section 1.2.9, but in more explicit detail here, where we have omitted the order subscript in the specification of B-splines, and the $x$-value of the specific interpolation datasite on which the operator $LB_i$ typically acts, to ease notation:

$$
\begin{bmatrix}
\beta_1 B_1 & \beta_1 B_2(0) & 0 & & & & & & & & & & & \\
\beta_2 B_1 & \beta_2 B_2(0) & 0 & 0 & 0 & 0 & & & & & & & & \\
LB_1 & LB_2 & \cdots & LB_{k-1} & LB_k(\tau_1) & 0 & & & & & & & & \\
\vdots & & & & \vdots & & & & & & & & & \\
 & & & & & 0 & & & & & & & & \\
LB_1 & LB_2 & \cdots & LB_{k-1} & LB_k(\tau_{k-2}) & 0 & 0 & 0 & 0 & & & & & \\
0 & 0 & 0 & LB_{k-1} & LB_k & \cdots & LB_{2k-3} & LB_{2k-2}(\tau_{k-1}) & 0 & & & & & \\
 & & & \vdots & & & & \vdots & & & & & & \\
 & & 0 & & & & & & 0 & & & & & \\
 & & 0 & LB_{k-1} & LB_k & \cdots & LB_{2k-3} & LB_{2k-2}(\tau_{2k-4}) & 0 & 0 & 0 & 0 & & \\
 & & 0 & 0 & 0 & 0 & \bullet & \bullet & \ddots & \bullet & \bullet & 0 & 0 & 0 \\
 & & & & & 0 & 0 & 0 & 0 & LB_{1+n-k} & LB_{2+n-k} & \cdots & LB_{n-1} & LB_n(\tau_{1+n-k}) \\
 & & & & & & & & & \vdots & & & & \vdots \\
 & & & & & & & & 0 & & & & & \\
 & & & & & & & & 0 & LB_{1+n-k} & LB_{2+n-k} & \cdots & LB_{n-1} & LB_n(\tau_{n-2})
\end{bmatrix}
$$

The inverse of matrix $\mathbf{A}$ has the structure below, where we used the notation $R_{x,y,z}$ to define a non-zero element's dependence upon the elements of rows $x, y, z$ of $\mathbf{A}$.

$$
\begin{bmatrix}
R_{1,2} & R_{1,2} & 0 & & & & & & & & & \\
R_{1,2} & R_{1,2} & 0 & 0 & 0 & 0 & & & & & & \\
R_{1,..,k} & R_{1,..,k} & R_{3,..,k} & \cdots & R_{3,..,k} & 0 & & & & & & \\
\vdots & & & & \vdots & & & & & & & \\
 & & & & & 0 & & & & & & \\
R_{1,..,k} & R_{1,..,k} & R_{3,..,k} & \cdots & R_{3,..,k} & 0 & 0 & 0 & 0 & & & \\
R_{1,..,2k-2} & R_{1,..,2k-2} & R_{3,..,2k-2} & \cdots & R_{3,..,2k-2} & R_{k+1,..,2k-2} & \cdots & R_{k+1,..,2k-2} & 0 & & & \\
\vdots & & & & & & & \vdots & & & & \\
 & & & & & & & & 0 & & & \\
R_{1,..,2k-2} & R_{1,..,2k-2} & R_{3,..,2k-2} & \cdots & R_{3,..,2k-2} & R_{k+1,..,2k-2} & \cdots & R_{k+1,..,2k-2} & 0 & 0 & 0 & 0 \\
R_{1,..,3k-4} & R_{1,..,3k-4} & R_{3,..,3k-4} & \cdots & R_{3,..,3k-4} & R_{k+1,..,3k-4} & \cdots & R_{k+1,..,3k-4} & R_{2k-1,..,3k-4} & \cdots & R_{2k-1,..,3k-4} & 0 \\
\vdots & & & & & & & & & & \vdots & \\
 & & & & & & & & & & & 0 \\
R_{1,..,3k-4} & R_{1,..,3k-4} & R_{3,..,3k-4} & \cdots & R_{3,..,3k-4} & R_{k+1,..,3k-4} & \cdots & R_{k+1,..,3k-4} & R_{2k-1,..,3k-4} & \cdots & R_{2k-1,..,3k-4} & 0 \\
\vdots & & & & & & & & & & & \ddots
\end{bmatrix}
$$

## 7.2 Global error consideration referenced in section 2.1.2

**Non-chaotic nature of $g(x)$**

Assume that for;

$$D^2g + \mu(g^2-1)Dg + g = 0, \quad g(a) = g_A, \;\; Dg(a) = g_B,$$
$$D^2h + \mu(h^2-1)Dh + h = 0, \quad h(a) = g_A + \epsilon_A, \;\; Dh(a) = g_B + \epsilon_B,$$

there exists $\phi > 0$ and $\theta > 0$ such that for $\epsilon_A, \epsilon_B < \phi$ then;

$$|h(x) - g(x)| < \theta\phi \;\text{ and }\; |Dh(x) - Dg(x)| < \theta\phi \quad \text{on any interval}[a..a + |\boldsymbol{\xi}|] \;,$$

such that these tend to zero with $\phi$.

**Error estimate for the first interior breakpoint**

After convergence of the collocation approximation in the first interval $[\xi_1 = 0..\xi_2]$ then we have by 1.14:

$$|(f_{\boldsymbol{\xi}} - g)(\xi_2)| \leq C|\boldsymbol{\xi}|^{2k-4}, \quad |D(f_{\boldsymbol{\xi}} - g)(\xi_2)| \leq C|\boldsymbol{\xi}|^{2k-4} \;.$$

**Error estimate at subsequent breakpoints**

From section 2.1.1 and reference to the above we assert the collocation approximation error at $\xi_3$ after convergence of the first 2 intervals is;

$$|f_{\boldsymbol{\xi}}(\xi_3) - g(\xi_3)| = |f_{\boldsymbol{\xi}}(\xi_3) - h(\xi_3) + h(\xi_3) - g(\xi_3)| \leq |f_{\boldsymbol{\xi}}(\xi_3) - h(\xi_3)| + |h(\xi_3) - g(\xi_3)| \;,$$

with $h$ being the target solution of the collocation approximation having inherited error from the first interval, such that;

$$|f_{\boldsymbol{\xi}}(\xi_3) - g(\xi_3)| \leq C_2|\boldsymbol{\xi}|^{2k-4} + \theta_1 C_1 |\boldsymbol{\xi}|^{2k-4}| \;,$$
$$|f_{\boldsymbol{\xi}}(\xi_3) - g(\xi_3)| \leq C\left(|\boldsymbol{\xi}|^{2k-4} + \theta_1|\boldsymbol{\xi}|^{2k-4}\right)| \;, \text{for } C = \max\{C_1, C_2\}$$

By induction this leads to general inequalities for $f_{\boldsymbol{\xi}}$ and $Df_{\boldsymbol{\xi}}$, the collocation approximation having converged for $m$ intervals;

$$\begin{aligned} |(f_{\boldsymbol{\xi}} - g)((\xi_{1+m})| &\leq C\left(|\boldsymbol{\xi}|^{2k-4} + \theta_{m-1}|\boldsymbol{\xi}|^{2k-4} + \ldots + \theta_{m-1}\ldots\theta_1|\boldsymbol{\xi}|^{2k-4}\right) \;, \\ &\leq C\left(|\boldsymbol{\xi}|^{2k-4} + \theta|\boldsymbol{\xi}|^{2k-4} + \ldots + \theta^{m-1}|\boldsymbol{\xi}|^{2k-4}\right) \;, \\ &= C\left(\frac{\theta^m - 1}{\theta - 1}\right)|\boldsymbol{\xi}|^{2k-4} \;, \;\; C = \max\{C_1, .., C_m\}, \;\; \theta = \max_{i=1..(m-1)}\left\{\sum_{j=i}^{m-1}\frac{\theta_j}{m-i}\right\} \;, \end{aligned}$$

and by the same reasoning,

$$|D(f_{\boldsymbol{\xi}} - g)((\xi_{1+m})| \leq C\left(\frac{\theta^m - 1}{\theta - 1}\right)|\boldsymbol{\xi}|^{2k-4} \;.$$

We note due to the specification of $\theta$ as the maximum of the arithmetic average of combinations of $\theta_j$ that it presents a rather conservative estimate, but it serves to give a mathematical description of that offered descriptively in the text that;

(i) $\theta$ being indirectly dependent on $|\boldsymbol{\xi}|$ impacts the errors at subsequent breakpoints,

(ii) errors are compounded for each successive interval to which the collocation method is applied.

## 7.3 Bounds for $err(x)$ under assumptions referenced in section 2.1.2

We define $err(x)$ to be the error in the differential equation calculated for the collocation approximation $f_{\boldsymbol{\xi}}$, which is numerical determinable as;

$$err(x) := D^2 f_{\boldsymbol{\xi}}(x) + \mu(f_{\boldsymbol{\xi}}(x)^2 - 1)Df_{\boldsymbol{\xi}}(x) + f_{\boldsymbol{\xi}}(x) \;,$$

Sections 2.1.1 and 2.1.2 discuss why $err(x)$ for an arbitrary interval $[A..B] \in [a..b]$ is not, pointwise or locally, necessarily dependent on $|g(x) - f_{\boldsymbol{\xi}}(x)|$. If we were to assume a traditional collocation problem with side conditions we can obtain an upper bound for this error in the following manner:

$$\begin{aligned}\|err(x)\|_\infty &= \|err(x) - D^2 g(x) - \mu(g(x)^2 - 1)Dg(x) - g(x)\|_\infty \ ,\\ &\leq \|D^2 f_{\boldsymbol{\xi}} - D^2 g\|_\infty + \|f_{\boldsymbol{\xi}} - g\|_\infty + \mu\|(f_{\boldsymbol{\xi}}^2 - 1)Df_{\boldsymbol{\xi}} - (g^2-1)Dg\|_\infty \ ,\\ &\leq C\big(|\boldsymbol{\xi}|^{(k-2)} + |\boldsymbol{\xi}|^k\big) + \mu\|(f_{\boldsymbol{\xi}}^2 - 1)Df_{\boldsymbol{\xi}} - (g^2-1)Dg\|_\infty \ ,\end{aligned}$$

where in the last step we have used equation 1.13 (assuming $k \geq 4$). To bound $\mu\|(f_{\boldsymbol{\xi}}^2 - 1)Df_{\boldsymbol{\xi}} - (g^2 - 1)Dg\|_\infty$, consider;

$$h = \frac{g + f_{\boldsymbol{\xi}}}{2}, \quad g = f_{\boldsymbol{\xi}} + \epsilon, \quad \text{where } \epsilon \leq |\boldsymbol{\xi}|^k,$$
$$\implies h = \begin{cases} g - \frac{1}{2}\epsilon \\ f_{\boldsymbol{\xi}} + \frac{1}{2}\epsilon \end{cases},$$

such that,

$$\begin{aligned}\mu\|(f_{\boldsymbol{\xi}}^2 - 1)Df_{\boldsymbol{\xi}} - (g^2-1)Dg\|_\infty &= \mu\|(h^2 - h\epsilon + \frac{1}{4}\epsilon^2 - 1)Df_{\boldsymbol{\xi}} - (h^2 + h\epsilon + \frac{1}{4}\epsilon^2 - 1)Dg\|_\infty \\ &\leq \mu\big(\|(h^2-1)(Df_{\boldsymbol{\xi}} - Dg)\|_\infty + \|\frac{1}{4}\epsilon^2(Df_{\boldsymbol{\xi}} - Dg)\|_\infty + \|h\epsilon Df_{\boldsymbol{\xi}}\|_\infty + \|h\epsilon Dg\|_\infty\big) \ .\end{aligned}$$

Ignoring phenomena such as overshoot we make the approximation that $|h| \leq 2$, $h \approx g$ and $Df_{\boldsymbol{\xi}} \approx Dg$, and again with equation 1.13,

$$\begin{aligned}\mu\|(f_{\boldsymbol{\xi}}^2 - 1)Df_{\boldsymbol{\xi}} - (g^2-1)Dg\|_\infty &\leq C\mu\left(3|\boldsymbol{\xi}|^{(k-1)} + \frac{|\boldsymbol{\xi}|^{2k}}{4}|\boldsymbol{\xi}|^{(k-1)}\right) + \mu|\boldsymbol{\xi}|^k\left(\|hDf_{\boldsymbol{\xi}}\|_\infty + \|hDg\|_\infty\right) \\ &\approx C\mu\left(3|\boldsymbol{\xi}|^{(k-1)} + \frac{|\boldsymbol{\xi}|^{2k}}{4}|\boldsymbol{\xi}|^{(k-1)}\right) + 2\mu|\boldsymbol{\xi}|^k\|gDg\|_\infty \ .\end{aligned}$$

Finally we make reference to $\|gDg\|_\infty$ as being obtainable from the phase plot as the rectangle of greatest area with diagonal vertices at $(0,0)$ and $(g, Dg)$. Figure 7.1 shows the relationship which leads to the supposition that,

$$\|gDg\|_\infty \leq 2\mu,$$

and hence to the overall estimate,

$$\|err(x)\|_\infty \leq C\big(|\boldsymbol{\xi}|^{(k-2)} + |\boldsymbol{\xi}|^k + \mu(3|\boldsymbol{\xi}|^{(k-1)} + \frac{1}{4}|\boldsymbol{\xi}|^{3k-1})\big) + 4\mu^2|\boldsymbol{\xi}|^k \ .$$

For given $\|err(x)\|_\infty$ we can rearrange the above inequality to assess a minimum value $C$, which can be substituted into equation 1.13 to give:

$$\|g - f_{\boldsymbol{\xi}}\|_\infty \leq \left(\frac{\|err(x)\|_\infty|\boldsymbol{\xi}|^k - 4\mu^2|\boldsymbol{\xi}|^{2k}}{|\boldsymbol{\xi}|^{(k-2)} + |\boldsymbol{\xi}|^k + \mu(3|\boldsymbol{\xi}|^{(k-1)} + \frac{1}{4}|\boldsymbol{\xi}|^{3k-1})}\right), \quad k \geq 4,$$

which permits for a numerical determination of the accuracy of the approximation to the unknown real solution, in the case of collocation with linear side conditions.

## 7.4 Derivation of $cost(f_{\boldsymbol{\xi}})$ referenced in section 2.1.3

Here we seek to define reliability by developing a measure associated with the cost of generating a collocation approximation under different parameters. This requires analysis of various stages of the numerical procedure in obtaining such an approximation, and we seek to create the measure in terms of an estimate of a floating point operation (flop) count in obtaining the approximation, which can be proxied and compared to the number of CPU seconds to generate the numerical approximations.

**Initialisation of problem,** $cost \approx O((k-2)l)$

The configuring of the solution, in terms of the initialisation of; $\boldsymbol{\xi}$ and $\mathbf{t}, \boldsymbol{\rho}$ and $\boldsymbol{\tau}$, $f_{\boldsymbol{\xi},0} = \sum_{i=1}^n \alpha_{i,0}B_{i,k,\mathbf{t}}$, we view as having a cost equivalent to the direct population of $\alpha_{i,0}$, and assigned the flop count above, and is essentially considered a priori to the problem.

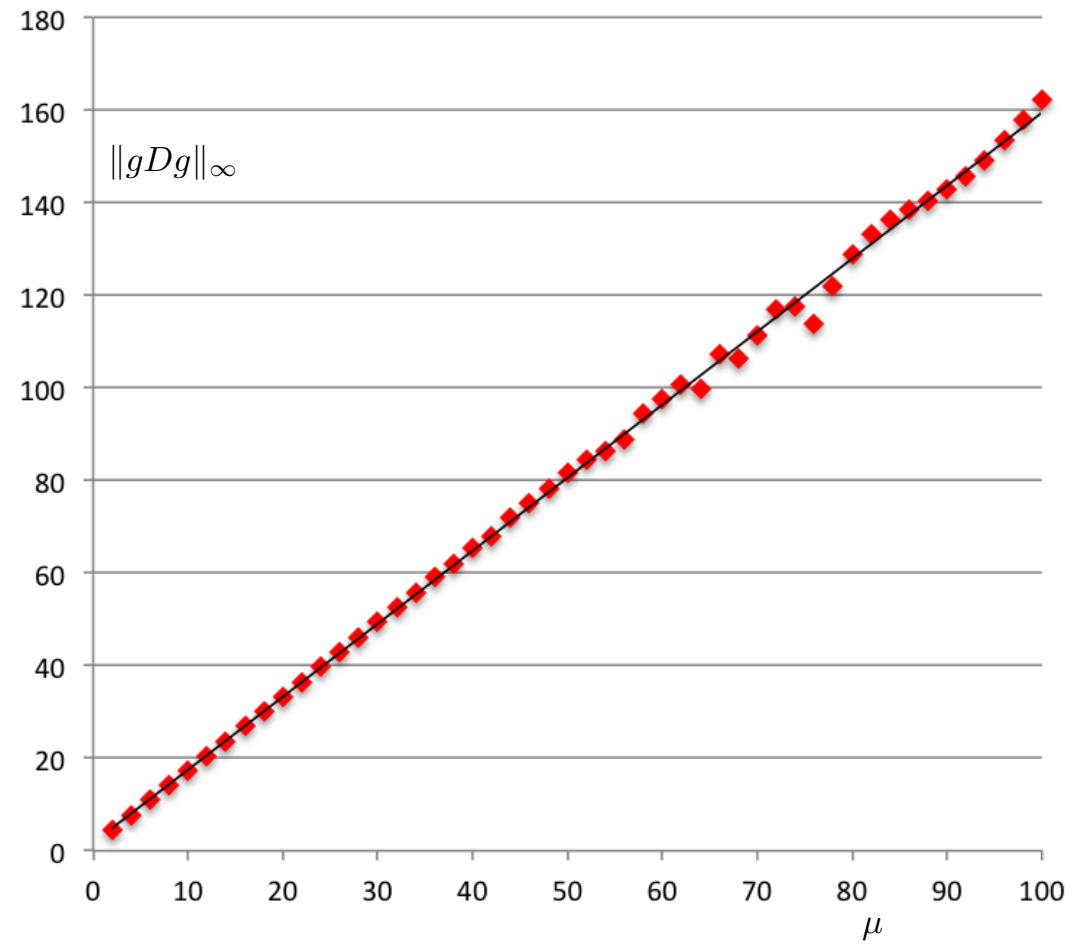

Figure 7.1: Scatter plot of $\|gDg\|_\infty$ against $\mu$ using data generated numerically from multiple sampling of phase plot data.

### B-spline value and derivatives calculation, *cost* $\approx O((6k! + 14k)(1 + l(k-2)))$

For each initialised knot sequence, $\mathbf{t}$, B-spline values and their first and second derivatives are required at $x = 0$ and $x = \boldsymbol{\tau}$. This can be calculated outside of Newton iteration. For any given value, $x$ the recurrence relation of equation 1.1 permits the calculation of B-spline values in $O(6k!)$ flops, simply counting the number of additions, divisions and multiplications, in the recurrence relation and the number of times it is used, recursively.

With these values equation 1.4 can be used to generate derivatives in $O(7k)$ more flops, and in turn with these derivatives the second derivative values can be generated in another $O(7k)$ flops. This is done for each $x$.

### Construction of matrix, A, and vector, s, in each iteration, *cost* $\approx O\left((k-2)(9k+12)lN + 16N\right)$

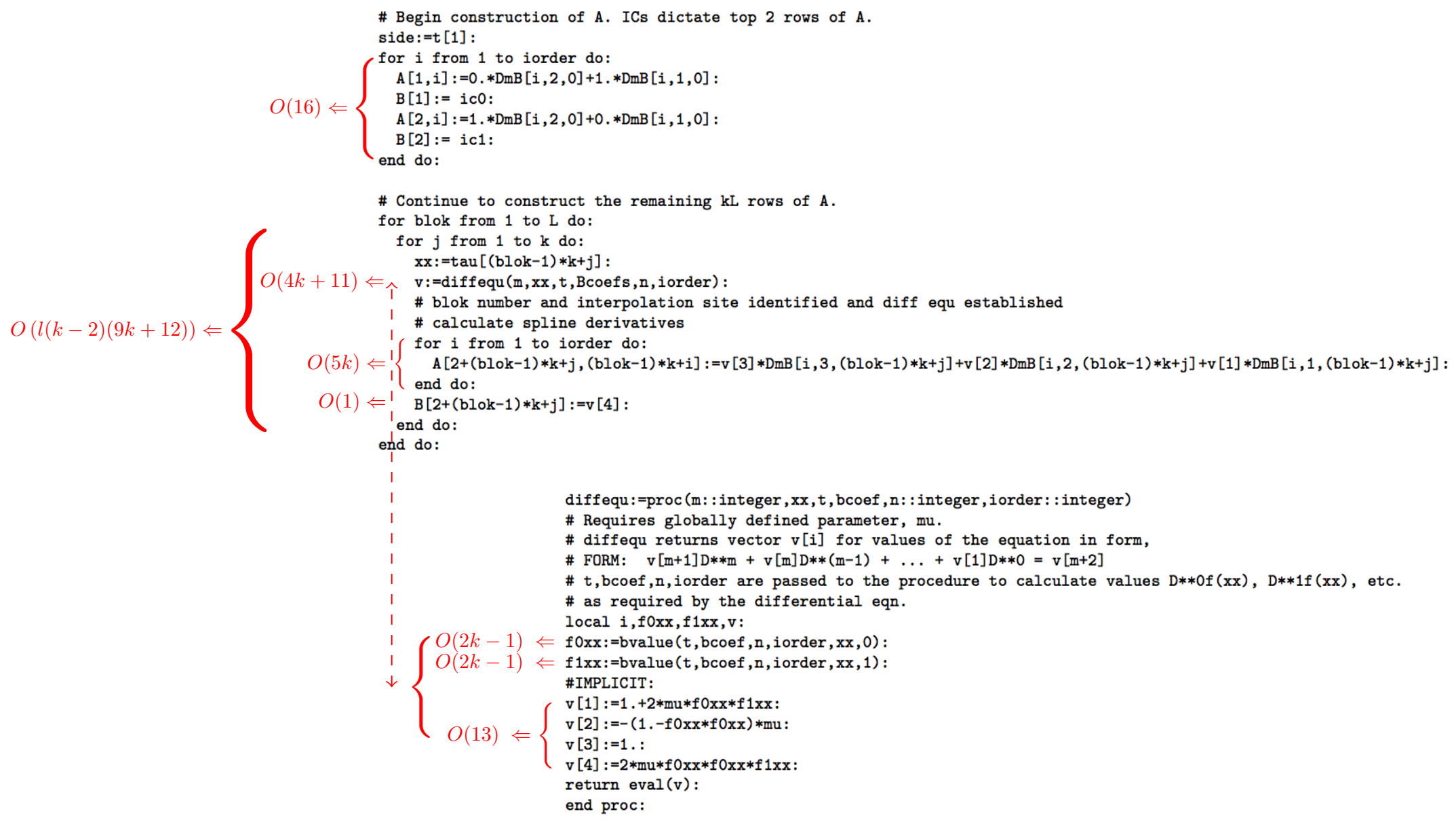

Figure 7.2: Outlining the derivation of an optimal flop count estimate for the construction of $\mathbf{A}$ and $\mathbf{s}$ in each iteration.

The construction of $\mathbf{A}$ requires the population of $k(k-2)l+4$ non-zero elements. Figure 7.2 attempts to outline the flop count estimates derived for an optimally efficient algorithm, which we stress is not what this algorithm here purports to be;

(i) establish the top left four elements related to initial conditions, $\approx O(16)$,

(ii) obtain differential equation coefficients, evaluating function and derivative values by utilising known B-spline values and known B-spline coefficients, $\alpha_{i,r}$, $\approx O(4k+11)$,

(iii) cycle through each relevant element of a specific row of $\mathbf{A}$ populating with combinations of known quantities, $\approx O(5k+1)$,

(iv) repeat (ii) and (iii) for $(k-2)l$ rows, $\approx O\left((k-2)(9k+12)lN\right)$.

**LinearSolve in each iteration,** *cost* $\approx O((\frac{2}{3}(k-1)^3 - (k-1)^2 + \frac{13}{3}(k-2) + \frac{1}{3})lN)$

To solve for $\boldsymbol{\alpha}$ requires the inversion of the banded matrix $\mathbf{A}$. Considering solving the system using row operations then the following broadly outlines a flop count;

(i) set $A_{1,1} = 1, A_{1,2} = 0$: 1d+2m+2a+1d $\approx O(6)$,

(ii) set $A_{2,1} = 0, A_{2,2} = 1$: 1m+1a+1d $\approx O(3)$,

(iii) set $\{A_{3,1}, \dots, A_{k,1}\} = 0$: (k-2)(1m+1a) $\approx O(2(k-2))$,

(iv) set $\{A_{3,2}, \dots, A_{k,2}\} = 0$: (k-2)(1m+1a) $\approx O(2(k-2))$,

(v) set $A_{3,3} = 1$, $\{A_{3,4}, \dots, A_{3,2+k}\} = 0$: [1d+(k-2)m+(k-2)a+1m](k-3) $\approx O(2(k-2)(k-3))$,

(vi) set $\{A_{4,3}, \dots, A_{k,3}\} = 0$: (1m+1a)(k-3) $\approx O(2(k-3))$,

(vii) repeat (v) and (vi) for a smaller $(k-3)\times(k-3)$ sub matrix, where previous consideration was for the top row and left column of a $(k-2)\times(k-2)$ sub matrix, and continue to repeat until all rows and all columns of the the original sub matrix are complete. $\sum_{j=1}^{k-2} O(2(j+1)(j-1)) \approx \sum_{j=1}^{k-2} O(2j^2) = O(\frac{2}{3}(k-1)^3 - (k-1)^2 + \frac{1}{3}(k-2) + \frac{1}{3})$

(viii) repeat steps (iii) to (vii) $l$ times for the consideration of each interval block used in the contsruction of $\mathbf{A}$.

**Flop Tally**

Combining all operations we find an approximate flop count after $N$ iterations, for $N, l \gg k$:

$$\begin{aligned}
cost(f_{\boldsymbol{\xi}}) &\approx O\Bigg( (k-2)l + (1+(k-2)l)(6k!+14k) + (k-2)(9k+12)lN \\
&\qquad + (\frac{2}{3}(k-1)^3 - (k-1)^2 + \frac{13}{3}(k-2) + \frac{1}{3})lN + 16N \Bigg), \\
&= O\Bigg( 6k! + 14k + (1+6k!+14k)(k-2)l + 16N \\
&\qquad + ((k-2)(9k+12) + \frac{2}{3}(k-1)^3 - (k-1)^2 + \frac{13}{3}(k-2) + \frac{1}{3})lN \Bigg), \\
&\approx O\Bigg( (\frac{2}{3}(k-1)^3 + 9k(k-2) - (k-1)^2 + \frac{49}{3}(k-2) + \frac{1}{3})lN \Bigg), \\
&\approx O\Bigg( (k^3 + k^2 + k)lN \Bigg),
\end{aligned}$$

where we highlight figure 7.3 shows the last step in the approximation is valid for the values of $k$ typically used in this study and we also observe that the dominant contributing factor to $cost(f_{\boldsymbol{\xi}})$ has been assumed to be the construction of matrix $\mathbf{A}$ and the row operations involved in solving the linear system, and the initialisation of the problem and the B-spline value and derivatives calculation has been neglected in deriving this order.

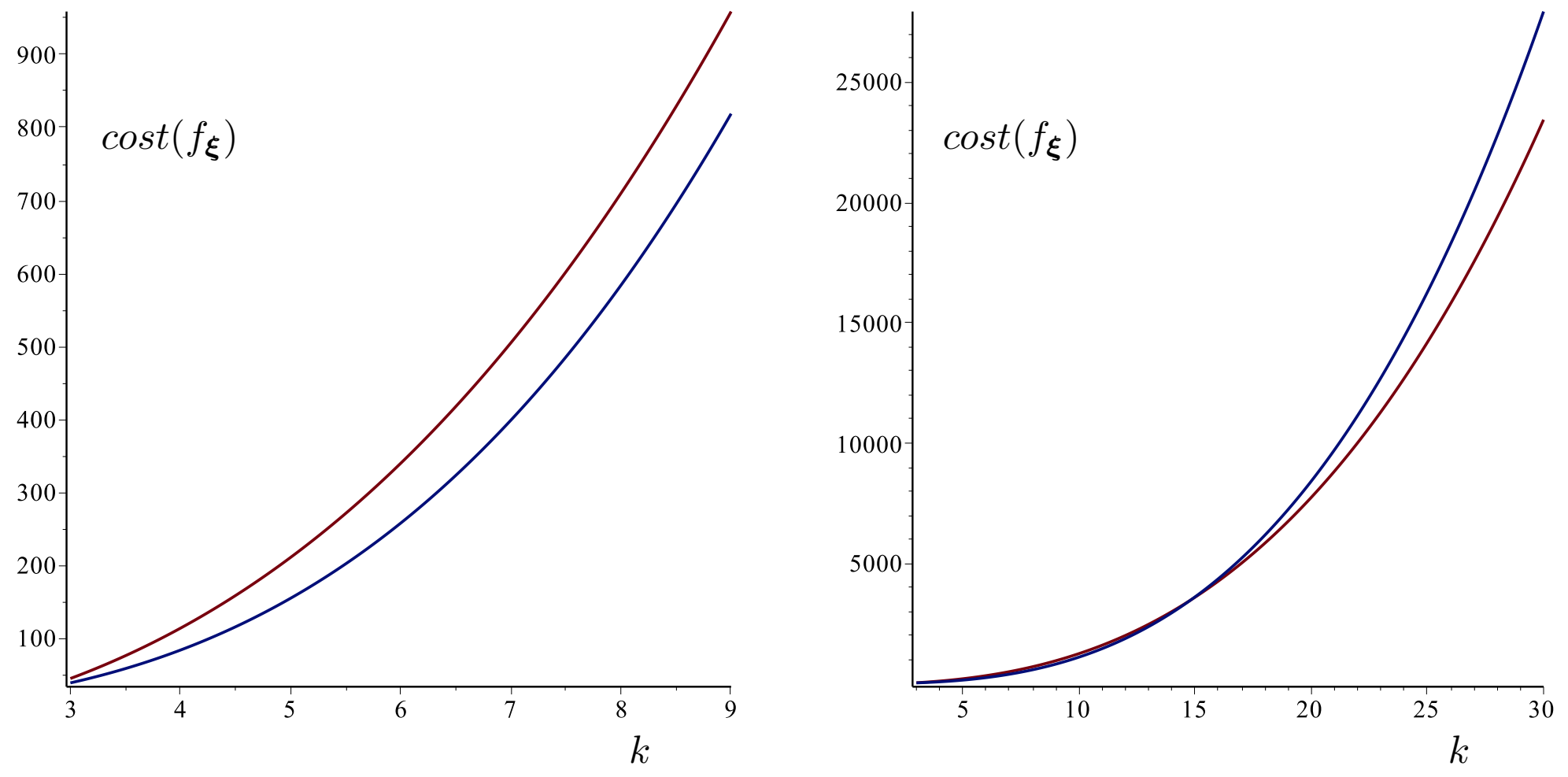

(a) Range $k = \{3, \ldots, 9\}$ (b) Range $k = \{3, \ldots, 30\}$

Figure 7.3: Plot of $\{cost(f_{\boldsymbol{\xi}}) = \frac{2}{3}(k-1)^3 + 9k(k-2) - (k-1)^2 + \frac{49}{3}(k-2) + \frac{1}{3},\ cost(f_{\boldsymbol{\xi}}) = k^3 + k^2 + k\}$ in the respective colours {red, blue}.

## 7.5 Determination of $\lambda$ referenced in section 6.1.2

We rearrange equation 6.2 to:

$$log(N^{seg.} - wN^{min.}) - log(N^{min}) = -\lambda log(w),$$

and produce relevant scatter plots, in figures 7.4 and 7.5, from results tables 5.2 and 5.3, to ascertain $\lambda$ as the gradient of the line minimising the error of least squares regression.

We conclude by comparing the trend line gradients of -0.876 and -1.078, that the appropriate value to adopt is $\lambda = 1.0$, which also pertains to be an aesthetic and qualitatively resolvable determination.

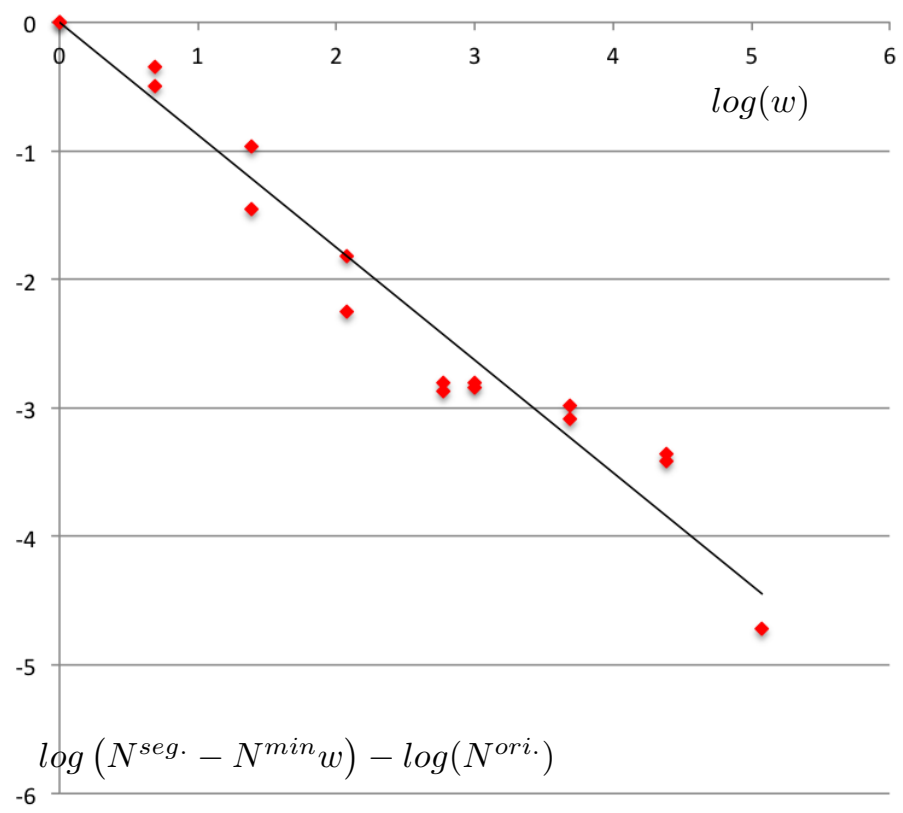

Figure 7.4: For the case $\mu = 3.0$, using the data from table 5.2, where the minimum least squares regression error has been obtained by setting $N^{min} = 3.305$, calculated iteratively, the gradient of the trend line is observed to be -0.876.

## 7.6 Determination of relationship between $\mu$ and $cost(f_{\boldsymbol{\xi}})$

We assert the relationship:

$$cost(f_{\boldsymbol{\xi}}) = A\mu^m,$$

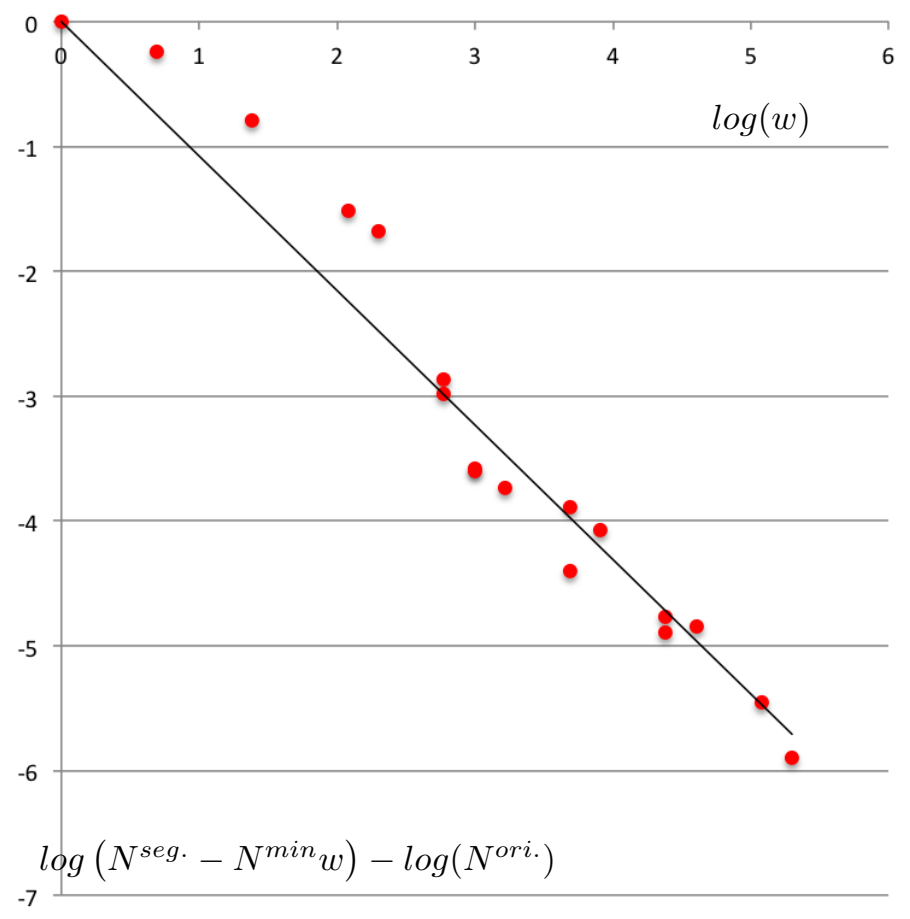

Figure 7.5: For the case $\mu = 10.0$, using the data from tables 5.2 and 5.3, where the minimum least squares regression error has been obtained by setting $N^{min} = 3.125$, calculated iteratively, and estimating for the case $l = 400$ that $N^{ori.} = 4000$, the gradient of the trend line is observed to be -1.078.

and rearrange to obtain the familiar:

$$log(cost(f_{\boldsymbol{\xi}})) = m\ log(\mu) + logA,$$

Figure 7.6 shows the scatter plot using the data in table 5.3 corresponding to the minimal cost achieved for producing an approximation of specific accuracy for a specific value $\mu \geq 2$. The gradient of the line, approximately $\frac{7}{5}$ gives:

$$cost(f_{\boldsymbol{\xi}}) \approx \frac{1}{4}\mu^{\frac{7}{5}} \text{ CPU seconds.}$$

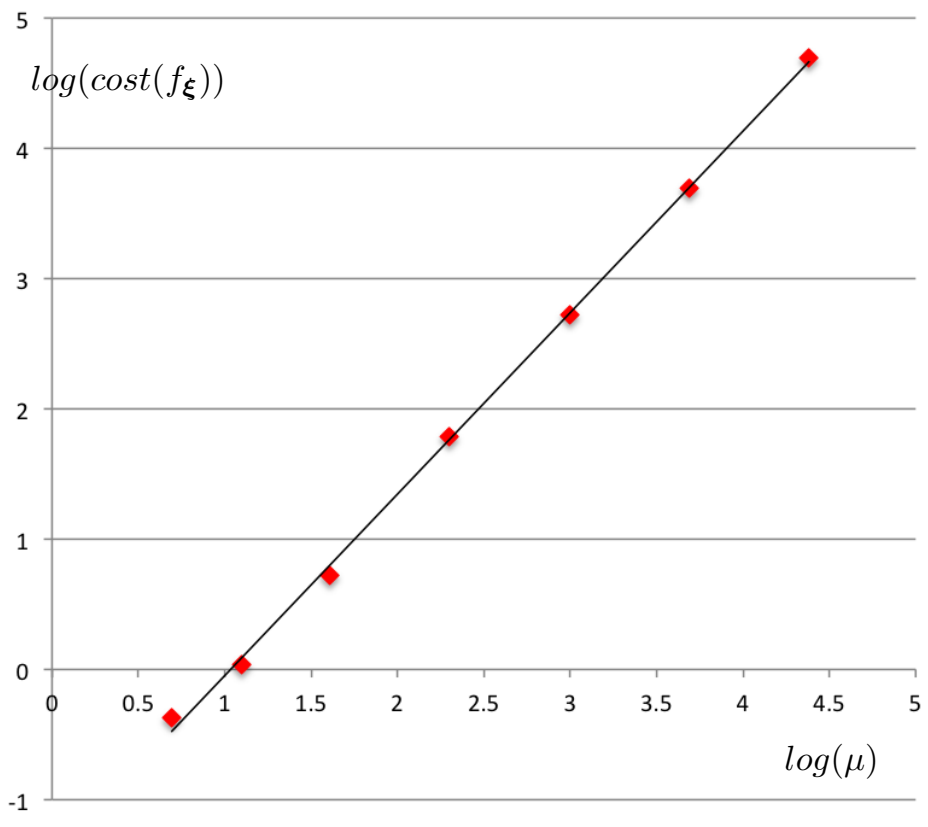

Figure 7.6: Scatter plot of the minimal $cost(f_{\boldsymbol{\xi}})$ in terms of CPU seconds against $\mu$ using data from table 5.3.

## 7.7 Maple algorithms

### 7.7.1 Main algorithm for execution of original or segmented collocation method

```
# Load packages
with(LinearAlgebra):
with(plots):

# ODE and display parameters
Digits   :=20                        :  #Maple machine accuracy
mu       :=10.0                      :  #ODE parameter
```

```
plotsper := 1                          :  #Number of plotdots in each interval.
gphboo   := 0                          :  #Boolean value to control whether plots are drawn (time intensive)

# Numerical solution parameters
b        :=40.                         :  #Defines solution range [a=0,b]
L        :=3125                        :  #Total number of intervals per for PP, i.e. used to create break sequence.
iter_max :=500                         :  #Max number of newton iterations allowed per segment
tolz     := 0.0001                     :  #Tolerance to stop window iterations evaluated at x=w*b.
windows  :=625                         :  #Number of segments to use in the segmented colloc method. (1 = original method)
m        :=2                           :  #Order of ODE
iorder   :=11                          :  #Order of the pp
k        :=iorder-m                    :  #Implied number of interpolation points per pp
b:=b/windows : L:=L/windows            :  #b and L rescaled according to segmented colloc method
f0       :=unapply(0.5*(1-tanh(x-3)),x):  #Supposed initial solution to be converted to PP
ic0      :=1.                          :  #Initial value x(0) supplied.
ic1      :=0.                          :  #Initial value Dx(0) supplied.
n        :=(iorder-m)*L+m              :  #the dimension on the spline is implied and is used by knotseq(),tauseq() to construct.
rho      :=rhoseq(k)                   :  #interpolation site spread, of Dimension(k)
timer    :=time()                      :  #Start CPU timer

# Initialize plot and variable arrays before any iterations
plotdots := L*windows*plotsper + 1:
xval     :=Vector(plotdots):
yval     :=Vector(plotdots):
DEErr    :=Vector(plotdots):
DmB      :=Array(1..iorder,1..(m+1),0..(k*L)):
iter_max_w  :=Vector(windows):
tolv     :=Vector(3):

# Begin window iteration....
for w from 1 to windows do:

  # Design the break seq, tau seq, t seq and initial guess solution.
  bbreak    :=seq(evalf(b*(i-1)/L+(w-1)*b),i=1..(L+1)): bbreak:=Vector(L+1,[bbreak]):
  tau       :=tauseq(bbreak,L,k,rho)      :  tau:=Vector(k*L,tau):
  t,n       :=knotseq(bbreak,L,iorder,m)  :  t:=Vector(n+iorder,t):
  Bcoefs:=Bcoef_ini(m,n,iorder,t,tau)     :  #calculate PP from given initial solution, f0
  tolv[1]:=10: tolv[2]:=100: tolv[3]:=1000:

  # Create array of B-Spline derivative values for each value, tau, outside of Newton iteration.
  # DmB(i,m,j) contains the i'th spline value at tau_j of the (m-1)th derivative
  side:=t[1]:
  left,iflag:=interv(t,n+iorder,side,left1,iflag1):
  DmB_temp:=bsplvd(t,iorder,side,left,A1,dbiatx1,m+1):
  for i from 1 to iorder do:
    for mm from 1 to m+1 do:
      DmB[i,mm,0] := DmB_temp[i,mm]:
    end do:
  end do:

  for blok from 1 to L do:
    for j from 1 to k do:
      xx:=tau[(blok-1)*k+j]:
      left,iflag:=interv(t,n+iorder,xx,left1,iflag1):
      DmB_temp:=bsplvd(t,iorder,xx,left,A1,dbiatx1,m+1):
      for i from 1 to iorder do:
        for mm from 1 to m+1 do:
          DmB[i,mm,(blok-1)*k+j]:=DmB_temp[i,mm]:
        end do:
      end do:
    end do:
  end do:

  # Begin Newton iteration....
  for iter from 1 to iter_max do:

    # Build a block matrix system Ax=B
    A:=Matrix(n,n):
    B:=Vector(n):

    # Begin construction of A. ICs dictate top 2 rows of A.
    side:=t[1]:
    for i from 1 to iorder do:
      A[1,i]:=0.*DmB[i,2,0]+1.*DmB[i,1,0]:
      B[1]:= ic0:
      A[2,i]:=1.*DmB[i,2,0]+0.*DmB[i,1,0]:
      B[2]:= ic1:
    end do:

    # Continue to construct the remaining kL rows of A.
    for blok from 1 to L do:
      for j from 1 to k do:
        xx:=tau[(blok-1)*k+j]:
        v:=diffequ(m,xx,t,Bcoefs,n,iorder):
        # blok number and interpolation site identified and diff equ established
        # calculate spline derivatives
        for i from 1 to iorder do:
          A[2+(blok-1)*k+j,(blok-1)*k+i]:=v[3]*DmB[i,3,(blok-1)*k+j]+v[2]*DmB[i,2,(blok-1)*k+j]+v[1]*DmB[i,1,(blok-1)*k+j]:
        end do:
        B[2+(blok-1)*k+j]:=v[4]:
      end do:
    end do:

    # Aquire solution in terms of B-Spline PP, and in PP-Form
    Bcoefs:=LinearSolve(A,B,method=SparseDirect):

    #Terminate iterations after tolerance reached
    tolv[3]:=tolv[2]: tolv[2]:=tolv[1]:
    tolv[1]:=bvalue(t,Bcoefs,n,iorder,b*w,0):
    if abs(tolv[1]-tolv[2])<tolz and abs(tolv[1]-tolv[3])<tolz then:
      iter_max_w[w]:=iter:
      break:
    end if:

  # next Newton iteration...
  end do:

  # Collate data for plot
  if gphboo > 0 then:
    dx := evalf((b*windows)/(plotdots-1)):
    for i from ((w-1)*plotsper*L+1) to (w*plotsper*L) do:
      xval[i] := (i-1)*dx:
      yval[i] := bvalue(t,Bcoefs,n,iorder,xval[i],0):
      DEerr[i] := bvalue(t,Bcoefs,n,iorder,xval[i],2)-mu*(1-yval[i]**2)*bvalue(t,Bcoefs,n,iorder,xval[i],1)+yval[i]:
    end do:
  end if:
```

```
  #Record num iterations
  if iter_max_w[w]=0 then:
    iter_max_w[w]:=iter_max:
  end if:

  # Initialize function and initial conditions for next window
  ic0:=bvalue(t,Bcoefs,n,iorder,b*w,0):
  ic1:=bvalue(t,Bcoefs,n,iorder,b*w,1):
  f0:=unapply((ic0+ic1*(x-b*w))*0.5*(1-tanh(x-b*w-3)),x):

# next window iteration...
end do:

# configure final plot point
xval[plotdots] := b*windows:
yval[plotdots] := bvalue(t,Bcoefs,n,iorder,xval[plotdots],0):
DEerr[plotdots] := bvalue(t,Bcoefs,n,iorder,xval[plotdots],2)-mu*(1-yval[plotdots]**2)*bvalue(t,Bcoefs,n,iorder,xval[plotdots],1)+yval[plotdots]:

#Print cost analysis of method
print("CPU time: ",time()-timer);
total_iter:=0:
for w from 1 to windows do:
  total_iter := total_iter + iter_max_w[w]:
end do:
print("total iterations ", total_iter);

# Design plots
Plot1:=pointplot(xval,yval, style=line,color=colorvec[12]):
Plot2:=pointplot(xval,DEerr, style=line, color=colorvec[12]):

# Finally, display results
display(Plot1,view=[0 .. b*windows,-2.5 .. 2.5]);
display(Plot2,view=[0 .. b*windows,-2 .. 2]);
```

### 7.7.2 Auxiliary algorithms

```
#######################################################
knotseq:=proc(bbreak,L::integer,iorder::integer,m::integer)
# knotseq creates a sequence of knots for collocation diff eqn solving
local j,i,t,n:
#left side knots
for i from 1 to iorder do:
   t[i]:=bbreak(1):
end do:
#intermediate knots
for j from 2 to L do:
   for i from 1 to (iorder-m) do:
      t[iorder+(j-2)*(iorder-m)+i]:=bbreak(j):
   end do:
end do:
#right side knots
for i from 1 to iorder do:
   t[iorder+(L-1)*(iorder-m)+i]:=bbreak(L+1):
end do:
n:=(iorder-m)*L+m:
return eval(t),eval(n):
end proc:
########################################################
rhoseq:=proc(n::integer)
# rhoseq creates a sequence of zeroes for the Legendre polynomial of degree, n>0.
local P, rhos, rhovec,i :
P[n]:=1/(n!*2**n)*diff((z**2-1)**n,[z$n]):
P[n]:=collect(P[n],z):
rhos:=fsolve(P[n]):
rhovec:=Vector(n):
for i from 1 to n do:
rhovec[i]:=rhos[i]:
end do:
return rhovec:
end proc:
#######################################################
tauseq:=proc(bbreak,L::integer,k::integer,rho)
# tauseq creates a repeated sequence of interpolation sites amongst the breaks given a spread vector rho
local j,i,xm,dx,tau:
for j from 1 to L do:
  xm:=(bbreak(j+1)+bbreak(j))*0.5:
  dx:=(bbreak(j+1)-bbreak(j))*0.5:
  for i from 1 to k do:
     tau[(j-1)*k+i]:=xm+dx*rho[i]:
  end do:
end do:
return eval(tau):
end proc:
#######################################################
diffequ:=proc(m::integer,xx,t,bcoef,n::integer,iorder::integer)
# Requires globally defined parameter, mu.
# diffequ returns vector v[i] for values of the equation in form,
# FORM:  v[m+1]D**m + v[m]D**(m-1) + ... + v[1]D**0 = v[m+2]
# t,bcoef,n,iorder are passed to the procedure to calculate values D**0f(xx), D**1f(xx), etc.
# as required by the differential eqn.
local i,f0xx,f1xx,v:
f0xx:=bvalue(t,bcoef,n,iorder,xx,0):
f1xx:=bvalue(t,bcoef,n,iorder,xx,1):
#IMPLICIT:
v[1]:=1.+2*mu*f0xx*f1xx:
v[2]:=-(1.-f0xx*f0xx)*mu:
v[3]:=1.:
v[4]:=2*mu*f0xx*f0xx*f1xx:
return eval(v):
end proc:
########################################################
Bcoef_ini:=proc(m::integer,n::integer,iorder::integer,t,tau)
#Requires globally defined function, f0(x).
local Q,iflag,i,tauf,ftau,Bcoefs:
for i from 2 to (n-m+1) do   :
  tauf[i] :=tau[i-1]         :
  ftau[i] :=f0(tauf[i])      :
end do                       :
  tauf[1] :=t[1]             :
  ftau[1] :=f0(t[1])         :
  tauf[n] :=t[n+iorder]      :
```

```
  ftau[n] :=f0(t[n+iorder])  :
#Bcoefs are established for above soln guess using splint.
Q,Bcoefs,iflag:=splint(tauf,ftau,t,n,iorder,Q1,Bcoefs1,iflag1): Bcoefs:=Vector(n,Bcoefs) :
return eval(Bcoefs):
end proc:
########################################################
colorvec:=["DeepSkyBlue","DarkTurquoise","Aquamarine","SpringGreen","GreenYellow","Yellow","Gold","Orange","Coral"_
_,"IndianRed","Maroon","Red"]:
#creation of a vector to assist with colored plots
########################################################
```

### 7.7.3 Cited algorithms

The above procedures also make calls to the algorithms of De Boor[1] translated from Fortran to Maple.

# Bibliography

[1] de Boor C., *A practical guide to splines, revised edition.* ISBN 0-387-95366-3.

[2] de Boor C., *On calculating with B-splines.* J. Approx. Theory, 1972.

[3] de Boor C. and Swartz B., *Collocation at Gaussian points*, J. Numerical Analysis, 1973.

[4] Casciola G. and Valori G., *An inductive proof of the derivative B-spline recursion formula.* http://citeseerx.ist.psu.edu/viewdoc/summary?doi=10.1.1.8.6081

[5] Schoenberg I. and Whitney A., *The positivity of translation determinants with an application to the interpolation problem by spline curves*, Trans. Amer. Math. Soc., 1953. http://www.ams.org/journals/tran/1953-074-02/S0002-9947-1953-0053177-X/

[6] Schoenberg I. and Curry H., *The fundamental spline functions and their limits*, J. Analyse Math, 1966.

[7] Liénard, A., *Etude des oscillations entretenues*, Revue générale de l'électricité, 1928.

[8] Panayotounakos D., Panayotounakou D., and Vanakis A., *On the lack of analytic solutions of the Van der Pol oscillator*, J. App. Math. & Mech., 2002.

[9] Loscalzo F. and Talbot T., *Spline function approximations for solutions of ordinary differential equations*, J. Numerical Analysis, 1967.

[10] Multiple authors, *Recent advances in sparse direct solvers*, Conference on structural mechanics in reactor tech., 2013. http://weisbecker.perso.enseeiht.fr/documents/smirt22.pdf

[11] Cartwright M., *Balthazar Van der Pol*, J. London Math. Soc., 1960.

[12] Nagumo J., Arimoto S., and Yoshizawa S., *An active pulse transmission line simulating nerve axon*, Proc. IRE, 1962.

[13] Cartwright J., Eguiluz V., Hernandez-Garcia E., and Piro O., *Dynamics of elastic excitable media*, Internat. J. Bifur. Chaos Appl. Sci. Engrg., 1999.